\documentclass[11pt]{article}
\usepackage{amsmath,amsthm,amsfonts,latexsym,amsopn,verbatim,amscd,amssymb}
\usepackage{pgfplots}
\usepackage{caption}
\usepackage{geometry}
\theoremstyle{plain}
\newtheorem{theorem}{Theorem}[section]

\newtheorem{Question}[theorem]{Question}

\newtheorem{Conjecture}[theorem]{Conjecture}

\newcommand{\bnum}{\begin{enumerate}}
\newcommand{\enum}{\end{enumerate}}

\numberwithin{equation}{section}

\DeclareMathOperator{\spec}{Spec}
\DeclareMathOperator{\L-spec}{L-spec}
\DeclareMathOperator{\Q-spec}{Q-spec}

\begin{document}
\title{\textbf{Signless Laplacian energies of non-commuting graphs of finite groups and related results}}
\author{Monalisha Sharma and Rajat Kanti Nath\footnote{Corresponding author}}
\date{}
\maketitle
\begin{center}\small{\it
Department of Mathematical Sciences, Tezpur University,\\ Napaam-784028, Sonitpur, Assam, India.\\

Emails:\, monalishasharma2013@gmail.com and rajatkantinath@yahoo.com}
\end{center}
\begin{abstract}
The non-commuting graph of a non-abelian group $G$ with center $Z(G)$ is a simple undirected graph whose vertex set is $G\setminus Z(G)$ and two vertices $x, y$ are adjacent if  $xy \ne yx$. In this study, we compute Signless Laplacian spectrum and Signless Laplacian energy of non-commuting graphs of finite groups. We obtain several conditions such that the non-commuting graph of  $G$ is Q-integral and observe  relations between energy, Signless Laplacian energy and Laplacian energy. In addition, we look into the energetic hyper- and hypo-properties of non-commuting graphs of finite groups. We also assess whether  the same graphs are Q-hyperenergetic and L-hyperenergetic.
\end{abstract}


\section{Introduction} \label{S:intro}
Recall that for a graph ${\mathcal{G}}$, its spectrum denoted by $\spec({\mathcal{G}})$ is  $\{(\lambda_1)^{k_1}, (\lambda_2)^{k_2},$ $\ldots, (\lambda_n)^{k_n}\}$, where $\lambda_i$'s are the eigenvalues of the adjacency matrix of $\mathcal{G}$ with multiplicities $k_i$ for $1 \leq i \leq n$, respectively. Let $A({\mathcal{G}})$ and $D({\mathcal{G}})$  be the adjacency matrix and degree matrix of ${\mathcal{G}}$, respectively. Consequently, $L({\mathcal{G}})  = D({\mathcal{G}}) - A({\mathcal{G}})$ gives the Laplacian matrix of ${\mathcal{G}}$.  The set $\{(\beta_1)^{b_1}, (\beta_2)^{b_2}, \ldots, (\beta_m)^{b_m}\}$,  where $\beta_j$'s are the eigenvalues of  $L({\mathcal{G}})$ with multiplicities $b_j$ for $1 \leq j \leq m$, respectively is referred to as the Laplacian spectrum of ${\mathcal{G}}$ and is denoted by $\L-spec({\mathcal{G}})$. While, $Q({\mathcal{G}})  = A({\mathcal{G}}) - D({\mathcal{G}})$ provides the Signless Laplacian matrix of ${\mathcal{G}}$.  The set $\{(\gamma_1)^{d_1}, (\gamma_2)^{d_2}, \ldots, (\gamma_l)^{d_l}\}$, where $\gamma_r$'s are the eigenvalues of  $Q({\mathcal{G}})$ with multiplicities $d_r$ for $1 \leq r \leq l$, respectively is known as the Signless Laplacian spectrum of ${\mathcal{G}}$ and is denoted by $\Q-spec({\mathcal{G}})$. A graph ${\mathcal{G}}$ is called integral, L-integral and Q-integral if eigenvalues of $A({\mathcal{G}})$, $L({\mathcal{G}})$ and $Q({\mathcal{G}})$ ate integers.

The energy, Laplacian energy and Signless Laplacian energy of   ${\mathcal{G}}$ are given by

~~$E(\mathcal{G}) = \sum_{\lambda \in \spec({\mathcal{G}})}|\lambda|$, $LE(\mathcal{G}) = \sum_{\beta \in \L-spec({\mathcal{G}})}\left|\beta - \frac{2|e(\mathcal{G})|}{|v(\mathcal{G})|}\right|$ and $LE^+(\mathcal{G}) = \sum_{\gamma \in \Q-spec({\mathcal{G}})}\left|\gamma - \frac{2|e(\mathcal{G})|}{|v(\mathcal{G})|}\right|$ respectively,
 where $v(\mathcal{G})$ and $e(\mathcal{G})$ are the set of vertices and edges of $\mathcal{G}$ respectively. A finite graph $\mathcal{G}$  is hyperenergetic if  $E(\mathcal{G}) > E(K_{|v(\mathcal{G})|})$ 
and is hypoenergetic if $E(\mathcal{G}) < |v(\mathcal{G})|$. Likewise, the graph $\mathcal{G}$ is L-hyperenergetic if $LE(\mathcal{G}) > LE(K_{|v(\mathcal{G})|})$ and is Q-hyperenergetic if $LE^+(\mathcal{G}) > LE^+(K_{|v(\mathcal{G})|})$. Gutman \cite{Gutman-99} and Walikar et. al \cite{Walikar-99} pioneered the research of hyperenergetic graphs in 1999. Gutman \cite{Gutman-07} first presented the hypoenergetic graph in 2007. L-hyperenergetic and Q-hyperenergetic  graphs were considered in \cite{fns20}.

Assume that $G$ is a finite non-abelian group and that $Z(G) = \{z \in G: zx=xz, \forall x \in G\}$ is the  centre of $G$. The non-commuting graph of $G$, represented by $\Gamma_G$, is a basic undirected network connecting two unique vertices $x$ and $y$ whenever $xy \neq yx$ with $G \setminus Z(G)$ as the vertex set. Erd{\H o}s and Neumann's study \cite{EN76} considered this graph for the first time in 1976. \cite{AAM06, BEKN14, DBBM10, GET14, MSZZ05, TEJ14} are some references on non-commuting graphs of finite groups.  Spectral aspects of this graph have also grabbed attention of numerous mathematicians. In \cite{GGAZ17}, Ghorbani et al. calculated  spectrum of $\Gamma_G$ for certain groups. The energy of the same was later explored by Ghorbani and Gharavi--Alkhansari in \cite{GGA17}. 
 In \cite{DDN18}, Dutta et al. computed the Laplacian spectrum of $\Gamma_G$  following which Dutta and Nath \cite{DN18} worked on Laplacian energy of the same. Signless Laplacian spectrum and  energy of $\Gamma_G$ are not yet computed. At present, \cite{AEN17} is the only paper where Abdussakir et al. investigated the Signless Laplacian spectrum of non-commuting graphs of dihedral groups. However, various spectra and energies (including Signless Laplacian spectrum and  energy) of the complement of $\Gamma_G$, known as commuting graph of $G$, are already computed in \cite{DBN-KJM-2020,DN-17,B2,DN-KJM-2018,DN-IJPAM-2021,FWNT21,SND-PJM-2022}. Various spectra and energies of  non-commuting conjugacy class graphs of several families of finite groups are computed in \cite{Nathetal-23}. 
 In this paper, we consider the following questions and answer them up to some extent.
 \begin{Question}\label{Q1-Q-integral}
 Which finite non-abelian groups give Q-integral non-commuting graphs?  
\end{Question}
 \begin{Question} \label{Q2-hyperenergetic}
 Are there any finite non-abelian group $G$ such that $\Gamma_G$ is hypoenergetic, hyperenergetic, L-hyperenergetic and Q-hyperenergetic?
\end{Question}
 \begin{Question}\label{Q3-energies}
 Which finite non-abelian groups satisfy the following inequalities?
 \begin{enumerate}
 	\item $E(\Gamma_G) \leq LE(\Gamma_G)$.
 	\item $LE(\Gamma_G) \leq LE^+(\Gamma_G)$.
 \end{enumerate}
\end{Question}
  We  compute $\Q-spec(\Gamma_G)$ and determine several conditions such that  $\Gamma_G$ is Q-integral.  We also compute $LE^+(\Gamma_G)$ for many families of finite non-abelian groups and compare various energies of $\Gamma_G$. These comparisons are interesting because mathematicians want to know about the graphs $\mathcal{G}$ such that $E(\mathcal{G}) \leq LE(\mathcal{G})$ (see \cite{DBN-KJM-2020,GAVBR-08,Liu-09,SSM-09}) and $LE(\mathcal{G}) \leq LE^+(\mathcal{G})$ (see \cite{BN-21,DBN-KJM-2020}).
Using  energies of $\Gamma_G$ for various classes of finite groups we determine whether they are hyperenergetic or hypoenergetic. In \cite{Gutman-78}, Gutman posed Conjecture \ref{con-hyper} which was disproved by different mathematicians providing counter examples (see \cite{Gutman-2011}).
\begin{Conjecture}\label{con-hyper}
Any finite graph $\mathcal{G} \ncong K_{|v(\mathcal{G})|}$ is non-hyperenergetic. 
\end{Conjecture} 
It is worth mentioning that we also disprove Conjecture \ref{con-hyper} by considering non-commuting graphs of finite groups (see Theorem \ref{com}(b)). We shall also determine finite groups such that their non-commuting graphs are Q-hyperenergetic and L-hyperenergetic.
 
We conclude this section with the following results which will help us to compute spectrum, Laplacian spectrum and Signless Laplacian spectrum of $\Gamma_G$  in the succeeding sections.


\begin{theorem}[\protect{\cite[Corollary 2.3]{DDN18}}]\label{L1}
Let $\mathcal{G}$ be a graph and $\mathcal{G}=l_1K_{m_1}\cup l_2K_{m_2}\cup \cdots \cup l_kK_{m_k}$, where $l_iK_{m_i}$ denotes the disjoint union of $l_i$ copies of $K_{m_i}$ for $1 \leq i \leq k$ and $m_1 < m_2 <\cdots < m_k$. Then
\begin{align*}
\L-spec(\mathcal{G}^c)=&\left\{ (0)^1,\left({\displaystyle \sum_{i=1}^k}l_im_i-m_k\right)^{l_k(m_k-1)},\left({\displaystyle \sum_{i=1}^k}l_im_i-m_{k-1}\right)^{l_{k-1}(m_{k-1}-1)},\ldots,\right.\\
& \left.\left({\displaystyle \sum_{i=1}^k}l_im_i-m_{1}\right)^{l_{1}(m_{1}-1)},\left({\displaystyle \sum_{i=1}^k}l_im_i\right)^{{\displaystyle \sum_{i=1}^k}l_i-1}\right \}.
\end{align*}
\end{theorem}

\begin{theorem}[\protect{\cite[Corollary 2.3]{WLH04} and \cite[Corollary 2.2]{ZWL13}}]\label{R1}
Let $\mathcal{G}$ be the complete $r$-partite graph\\ $K_{\underbrace{p_1,  \dots, p_1}_{a_1\text{-times}} \,, \, \underbrace{p_2,  \dots, p_2}_{a_2\text{-times}} \,,\, \ldots, \,\underbrace{p_s,  \dots, p_s}_{a_s\text{-times}}}:= K_{a_1.p_1,a_2.p_2,\ldots,a_s.p_s}$ on $n$ vertices, where $r = a_1 + a_2 + \cdots + a_s$. Then
\begin{enumerate}
\item the characteristic polynomial of $A({\mathcal{G}})$ is
\[P_\mathcal{G}(x):=x^{n-r}{\displaystyle \prod_{i=1}^s}(x+p_i)^{a_i-1}\left({\displaystyle \prod_{i=1}^s}(x+p_i)-{\displaystyle \sum_{j=1}^s}a_jp_j{\displaystyle \prod_{i=1,i \neq j}^s}(x+p_i)\right).\]
\item the characteristic polynomial of $Q({\mathcal{G}})$ ($Q$-polynomial) is
\[Q_\mathcal{G}(x):={\displaystyle \prod_{i=1}^s} (x-n+p_i)^{a_i(p_i-1)}{\displaystyle \prod_{i=1}^s}(x-n+2p_i)^{a_i}\left(1-{\displaystyle \sum_{i=1}^s}\frac{a_ip_i}{x-n+2p_i}\right).\]
\end{enumerate}
\end{theorem}
%

\section{$\frac{G}{Z(G)}$ is isomorphic to $D_{2m}$}

 Here, we primarily compute the Signless Laplacian spectrum and Signless Laplacian energy of $\Gamma_G$, where $G$ is isomorphic to $D_{2m}$, $QD_{2^n}$, $M_{2rs}$, $Q_{4n}$ and  $U_{6n}$. Further we compare  different energies of $\Gamma_G$ and look into the hyper- and hypo-energetic properties of $\Gamma_G$  for each of the above-mentioned groups. The  energy and Laplacian energy of $\Gamma_G$  for each of the aforementioned groups have already been determined and it is crucial to note.

\subsection{The dihedral groups, $D_{2m}$}
We consider $D_{2m}:= \langle a, b : a^m = b^2 = 1; bab^{-1} = a^{-1} \rangle$, the dihedral groups of order $2m$ (where $m > 2$). Results regarding different energies of non-commuting graphs of $D_{2m}$ are given below.
\begin{theorem}[\protect{\cite[Corollary 4.1.7 and (4.3.e)]{FWNT21}}]\label{d}  Let $G$ be isomorphic to $D_{2m}$.
\begin{enumerate}
\item If $m$ is odd then
\[E(\Gamma_{D_{2m}}) = (m-1)+\sqrt{(m-1)(5m-1)} \, \text{ and } \,
LE(\Gamma_{D_{2m}}) = \frac{2m(m-1)(m-2)+2m(2m-1)}{2m-1}.\]
\item If $m$ is even then
\[E(\Gamma_{D_{2m}}) = (m-2)+\sqrt{(m-2)(5m-2)} \, \text{ and } \,
LE(\Gamma_{D_{2m}}) =  \frac{m(m-2)(m-4)+2m(m-1)}{m-1}.\]
\end{enumerate}
\end{theorem}

\begin{theorem}\label{Dihedral1}
Let $G$ be isomorphic to $D_{2m}$, where $m$ is odd. Then 
\[\Q-spec(\Gamma_{D_{2m}})\!=\!\left\{(m)^{m-2}, (2m-3)^{m-1}, \left(\frac{4m-3+\sqrt{8m^2-16m+9}}{2}\right)^1, \left(\frac{4m-3-\sqrt{8m^2-16m+9}}{2}\right)^1\right\}\]
 and
\begin{align*}
LE^+(\Gamma_{D_{2m}}) = \begin{cases}
	 \frac{9}{5}+\sqrt{33}, & \mbox{if $m=3$}\vspace{.2cm}\\ \frac{2m^3-10m^2+12m-3}{2m-1}+\sqrt{8m^2-16m+9}, & \mbox{if $m\geq 5$.}
\end{cases}   
\end{align*}
\end{theorem}
\begin{proof}
If $G \cong D_{2m}$ and $m$ is odd then $|v(\Gamma_{D_{2m}})|=2m-1$ and $\Gamma_{D_{2m}}=K_{m.1,1.(m-1)}$. Using Theorem \ref{R1}(b), we have 
\begin{align*}
Q_{\Gamma_{D_{2m}}}(x)=&{\displaystyle \prod_{i=1}^2} (x-(2m-1)+p_i)^{a_i(p_i-1)}{\displaystyle \prod_{i=1}^2}(x-(2m-1)+2p_i)^{a_i}\left(1-{\displaystyle \sum_{i=1}^2}\frac{a_ip_i}{x-(2m-1)+2p_i}\right)\\
  = \,& (x-(2m-2))^{0}(x-m)^{m-2}(x-2m+3)^{m}(x-1)\left(1-\frac{m}{x-2m+3}-\frac{m-1}{x-1}\right)\\
   =\,& (x-m)^{m-2}(x-(2m-3))^{m-1}(x^2-(4m-3)x+2m^2-2m).
\end{align*}
Thus, $\Q-spec(\Gamma_{D_{2m}})=\left\{(m)^{m-2}, (2m-3)^{m-1}, \left(\frac{4m-3+\sqrt{8m^2-16m+9}}{2}\right)^1,\left(\frac{4m-3-\sqrt{8m^2-16m+9}}{2}\right)^1\right\}$.

 Number of edges in $\Gamma_{D_{2m}}^c$ is $\frac{m^2-3m+2}{2}$. Thus, $|e(\Gamma_{D_{2m}})| = \frac{(2m-1)(2m-1-1)}{2}-\frac{m^2-3m+2}{2}=\frac{3m(m-1)}{2}$.
  Now
\[
\left|m - \frac{2|e(\Gamma_{D_{2m}})|}{|v(\Gamma_{D_{2m}})|}\right| = \left|\frac{-m(m-2)}{2m-1}\right|= \frac{m(m-2)}{2m-1}, 
\] 
\[\left|2m-3 - \frac{2|e(\Gamma_{D_{2m}})|}{|v(\Gamma_{D_{2m}})|}\right| = \left|\frac{m^2-5m+3}{2m-1}\right|= \begin{cases}
	\frac{3}{5}, & \mbox{if $m=3$}\vspace{.2cm}\\
	\frac{m^2-5m+3}{2m-1}, & \mbox{if $m \geq 5$,}
\end{cases}
\]
\begin{align*}
\left|\frac{1}{2}\left(4m-3+\sqrt{8m^2-16m+9}\right) - \frac{2|e(\Gamma_{D_{2m}})|}{|v(\Gamma_{D_{2m}})|}\right| = &\left|\frac{1}{2}\left(\sqrt{8m^2-16m+9}+m-\frac{3}{2}+\frac{3}{4m-2}\right)\right|\\
 =&\frac{1}{2}\left(\sqrt{8m^2-16m+9}+m-\frac{3}{2}+\frac{3}{4m-2}\right)
\end{align*}
 and
\begin{align*}
\left|\frac{1}{2}\left(4m-3-\sqrt{8m^2-16m+9}\right) - \frac{2|e(\Gamma_{D_{2m}})|}{|v(\Gamma_{D_{2m}})|}\right| =& \left|\frac{1}{2}\left(-\sqrt{8m^2-16m+9}+m-\frac{3}{2}+\frac{3}{4m-2}\right)\right|\\
 =&\frac{1}{2}\left(\sqrt{8m^2-16m+9}-m+\frac{3}{2}-\frac{3}{4m-2}\right).
\end{align*}
Therefore, for $m=3$  we have 
$LE^+(\Gamma_{D_{2m}})= \frac{9}{5}+\sqrt{33}$.
For $m \geq 5$  we have
\begin{align*}
LE^+(\Gamma_{D_{2m}}) = &(m-2)\times \frac{m(m-2)}{2m-1}+(m-1)\times \frac{m^2-5m+3}{2m-1}+\\
&\frac{1}{2}\left(\sqrt{8m^2-16m+9}+m-\frac{3}{2}+\frac{3}{4m-2}
\right) +\frac{1}{2}\left(\sqrt{8m^2-16m+9}-m+\frac{3}{2}-\frac{3}{4m-2}\right)
\end{align*}
and the result follows on simplification.
\end{proof}

\begin{theorem}\label{Dihedral2}
Let $G$ be isomorphic to $D_{2m}$, where  $m$ is even. Then 
\begin{align*}
\Q-spec&(\Gamma_{D_{2m}})\\
&=\left\lbrace (2m-4)^\frac{m}{2}, (m)^{m-3}, (2m-6)^{\frac{m}{2}-1}, \left(2m-3+\sqrt{2m^2-8m+9}\right)^1, \left(2m-3-\sqrt{2m^2-8m+9}\right)^1\right\rbrace
\end{align*} and 
\begin{align*}
LE^+(\Gamma_{D_{2m}}) = \begin{cases}
\frac{m^3-4m^2+12}{2m-2}+2\sqrt{2m^2-8m+9}, & \mbox{if $4 \leq m \leq 8$}\vspace{.2cm}\\
\frac{m^3-8m^2+16m-6}{m-1}+2\sqrt{2m^2-8m+9}, & \mbox{if $m \geq 10$.}
\end{cases}   
\end{align*}
\end{theorem}
\begin{proof}
If $G \cong D_{2m}$ and $m$ is even then $|v(\Gamma_{D_{2m}})|=2m-2$ and $\Gamma_{D_{2m}}=K_{\frac{m}{2}.2,1.(m-2)}$. Using Theorem \ref{R1}(b), we have 
\begin{align*}
Q_{\Gamma_{D_{2m}}}(x)=&{\displaystyle \prod_{i=1}^2} (x-(2m-2)+p_i)^{a_i(p_i-1)}{\displaystyle \prod_{i=1}^2}(x-(2m-2)+2p_i)^{a_i}\left(1-{\displaystyle \sum_{i=1}^2}\frac{a_ip_i}{x-(2m-2)+2p_i}\right)\\
  = \,& (x-(2m-4 ))^{\frac{n}{2}}(x-m)^{m-3}(x-2m+6)^{\frac{m}{2}}(x-2)\left(1-\frac{m}{x-2m+6}-\frac{m-2}{x-2}\right)\\
   =\,& (x-(2m-4))^{\frac{m}{2}}(x-m)^{m-3}(x-(2m-6))^{\frac{m}{2}-1}(x^2-(4m-6)x+2m^2-4m).
\end{align*}
Therefore
\begin{align*}
\Q-spec&(\Gamma_{D_{2m}})\\
=&\left\lbrace (2m-4)^\frac{m}{2}, (m)^{m-3}, (2m-6)^{\frac{m}{2}-1}, \left(2m-3+\sqrt{2m^2-8m+9}\right)^1, \left(2m-3-\sqrt{2m^2-8m+9}\right)^1\right\rbrace.
\end{align*}

 Number of edges in $\Gamma_{D_{2m}}^c$ is $\frac{m^2-4m+6}{2}$ and so $|e(\Gamma_{D_{2m}})| = \frac{(2m-2)(2m-2-1)}{2}-\frac{m^2-4m+6}{2}=\frac{3m(m-2)}{2}$.   Now
\[
\left|2m-4 - \frac{2|e(\Gamma_{D_{2m}})|}{|v(\Gamma_{D_{2m}})|}\right| = \left|\frac{(m-2)(m-4)}{2m-2}\right| =\frac{(m-2)(m-4)}{2m-2}, 
\]
\[
\left|m - \frac{2|e(\Gamma_{D_{2m}})|}{|v(\Gamma_{D_{2m}})|}\right| = \left|\frac{(-m^2+4m)}{2m-2}\right|= \frac{(m^2-4m)}{2m-2},
\]
\[
\left|2m-6 - \frac{2|e(\Gamma_{D_{2m}})|}{|v(\Gamma_{D_{2m}})|}\right| = \left|\frac{(m^2-10m+12)}{2m-2}\right|= \begin{cases}
\frac{(-m^2+10m-12)}{2m-2}, & \mbox{if $m \leq 8$}\vspace{.2cm}\\
\frac{(m^2-10m+12)}{2m-2}, & \mbox{if $m \geq 10$,}
\end{cases}
\]
\begin{align*}
\left|2m-3+\sqrt{2m^2-8m+9} - \frac{2|e(\Gamma_{D_{2m}})|}{|v(\Gamma_{D_{2m}})|}\right| = & \left|\sqrt{2m^2-8m+9}+\frac{m}{2}-\frac{3}{2}+\frac{3}{2m-2}\right|\\
= & \sqrt{2m^2-8m+9}+\frac{m}{2}-\frac{3}{2}+\frac{3}{2m-2}
\end{align*}
and
\begin{align*}
\left|2m-3-\sqrt{2m^2-8m+9} - \frac{2|e(\Gamma_{D_{2m}})|}{|v(\Gamma_{D_{2m}})|}\right| = & \left|-\sqrt{2m^2-8m+9}+\frac{m}{2}-\frac{3}{2}+\frac{3}{2m-2}\right|\\
= &\sqrt{2m^2-8m+9}-\frac{m}{2}+\frac{3}{2}-\frac{3}{2m-2}.
\end{align*}
Therefore, for $4 \leq m \leq 8$, we have  
\begin{align*}
	LE^+(\Gamma_{D_{2m}}) = &\frac{m}{2}\times \frac{(m-2)(m-4)}{2m-2}+(m-3)\times \frac{(m^2-4m)}{2m-2}+ \left(\frac{m}{2}-1\right) \times \frac{-(m^2-10m+12)}{2m-2} \\
	&+ \sqrt{2m^2-8m+9}+\frac{m}{2}-\frac{3}{2}+\frac{3}{2m-2}+\sqrt{2m^2-8m+9}-\frac{m}{2}+\frac{3}{2}-\frac{3}{2m-2}   
\end{align*}
and for $m \geq 10$, we have 
\begin{align*}
LE^+(\Gamma_{D_{2m}}) = & \frac{m}{2}\times \frac{(m-2)(m-4)}{2m-2}+(m-3)\times\frac{(m^2-4m)}{2m-2}+ \left(\frac{m}{2}-1\right) \times \frac{(m^2-10m+12)}{2m-2} \\
&+ \sqrt{2m^2-8m+9}+\frac{m}{2}-\frac{3}{2}+\frac{3}{2m-2}+\sqrt{2m^2-8m+9}-\frac{m}{2}+\frac{3}{2}-\frac{3}{2m-2}.    
\end{align*}
The required expressions for $LE^+(\Gamma_{D_{2m}})$ can be obtained on simplification.
\end{proof}

\begin{theorem}\label{D_{2m}}
If $G$ is isomorphic to $D_{2m}$ then
\begin{enumerate}
\item $E(\Gamma_{D_{2m}}) \leq LE^+(\Gamma_{D_{2m}}) \leq LE(\Gamma_{D_{2m}})$, equality holds if and only if $G \cong D_8$.
\item $\Gamma_{D_{2m}}$ is non-hypoenergetic as well as non-hyperenergetic. 
\item $\Gamma_{D_6}$ is L-hyperenergetic but not Q-hyperenergetic. $\Gamma_{D_8}$ is not L-hyperenergetic and  not Q-hyperenergetic. If $m \ne 3, 4$ then $\Gamma_{D_{2m}}$ is Q-hyperenergetic and L-hyperenergetic.
\end{enumerate}

\end{theorem}
\begin{proof} 
\noindent (a)	\textbf{Case 1:} $m$ is odd

For $m = 3$, using Theorems \ref{d} and \ref{Dihedral1}, we have
$E(\Gamma_{D_6})=2+2\sqrt{7}$, $LE(\Gamma_{D_6})=\frac{42}{5}$ and $LE^+(\Gamma_{D_6})=\frac{9}{5}+\sqrt{33}$. Clearly, $E(\Gamma(D_{6}))< LE^+(\Gamma(D_{6})) < LE(\Gamma(D_{6}))$. 


 For $m \geq 5$, using Theorems \ref{d} and \ref{Dihedral1}, we have 
\begin{equation}\label{E1}
LE(\Gamma_{D_{2m}})-LE^+(\Gamma_{D_{2m}})=\frac{8m^2-10m+3}{2m-1}-\sqrt{8m^2-16m+9}     
\end{equation}
and
\begin{equation}\label{E2}
LE^+(\Gamma_{D_{2m}})-E(\Gamma_{D_{2m}})=\frac{2m^2(m-6)+15m-4}{2m-1}+\sqrt{8m^2-16m+9}-\sqrt{5m^2-6m+1}.
\end{equation}
Since $8m^2-10m+3 > 0$, $(2m-1)\sqrt{8m^2-16m+9} > 0$ and $(8m^2-10m+3)^2-\left(\sqrt{8m^2-16m+9}\right)^2(2m-1)^2=32m^3(m-2)+8m(5m-1)>0$
we have $8m^2-10m+3 - (2m-1)\sqrt{8m^2-16m+9} > 0$. Therefore, by \eqref{E1}, $(2m-1)(LE(\Gamma_{D_{2m}})-LE^+(\Gamma_{D_{2m}})) > 0$. Hence, $LE(\Gamma_{D_{2m}}) > LE^+(\Gamma_{D_{2m}})$.  

Again, $\sqrt{8m^2-16m+9} > 0, \sqrt{5m^2-6m+1} > 0$ and $\left(\sqrt{8m^2-16m+9}\right)^2-\left(\sqrt{5m^2-6m+1}\right)^2=m(3m-10)+8>0$. Therefore, 
$\sqrt{8m^2-16m+9}-\sqrt{5m^2-6m+1}>0$. Since $2m^2(m-6)+15m-4 > 0$ we have $\frac{2m^2(m-6)+15m-4}{2m-1}+\sqrt{8m^2-16m+9}-\sqrt{5m^2-6m+1}>0$. Therefore, by \eqref{E2}, $LE^+(\Gamma_{D_{2m}}) >  E(\Gamma_{D_{2m}})$. Hence, 
$E(\Gamma_{D_{2m}}) < LE^+(\Gamma_{D_{2m}}) < LE(\Gamma_{D_{2m}})$.
%

\vspace{.5cm}

\noindent \textbf{Case 2:} $m$ is even

For $4 \leq m \leq 8$, using Theorems \ref{d} and \ref{Dihedral2}, we have 
\begin{equation}\label{E3}
	LE(\Gamma_{D_{2m}})-LE^+(\Gamma_{D_{2m}})=\frac{m^3-4m^2+12m-12}{2m-2}-2\sqrt{2m^2-8m+9}   
\end{equation}
and
\begin{equation}\label{E4}
	LE^+(\Gamma_{D_{2m}})-E(\Gamma_{D_{2m}})=\frac{(m-4)(m^2-2m-2)}{m-1}+2\sqrt{2m^2-8m+9}-\sqrt{5m^2-12m+4}.
\end{equation} 

Since $m^3-4m^2+12m-12 > 0$, $2(2m-2)\sqrt{2m^2-8m+9} > 0$ and
\begin{center}
 $(m^3-4m^2+12m-12)^2-\left(2\sqrt{2m^2-8m+9}\right)^2(2m-2)^2=m(m-4)^2(m-2)(m^2+2m-4) \geq 0$
\end{center}
(equality holds if and only if $m = 4$).
It follows that $m^3-4m^2+12m-12 - 2(2m-2)\sqrt{2m^2-8m+9} \geq 0$. Therefore, by \eqref{E3}, $(2m-2)\left(LE(\Gamma_{D_{2m}})-LE^+(\Gamma_{D_{2m}})\right) \geq 0$. Hence, $LE(\Gamma_{D_{2m}}) \geq LE^+(\Gamma_{D_{2m}})$ equality holds if and only if $G \cong D_8$.  

Again, $2\sqrt{2m^2-8m+9} > 0, \sqrt{5m^2-12m+4} > 0$ and $\left(2\sqrt{2m^2-8m+9}\right)^2-\left(\sqrt{5m^2-12m+4}\right)^2=$ $(m - 4)(3m - 8) \geq 0$. Therefore, 
$2\sqrt{2m^2-8m+9}-\sqrt{5m^2-12m+4} \geq 0$ (equality holds if and only if $m = 4$). Since $(m-4)(m^2-2m-2) \geq 0$ we have $\frac{(m-4)(m^2-2m-2)}{m-1}+2\sqrt{2m^2-8m+9}-\sqrt{5m^2-12m+4} \geq 0$ (equality holds if and only if $m = 4$). Therefore, by \eqref{E4}, $LE^+(\Gamma_{D_{2m}}) \geq  E(\Gamma_{D_{2m}})$. Hence, 
$E(\Gamma_{D_{2m}}) \leq LE^+(\Gamma_{D_{2m}}) \leq LE(\Gamma_{D_{2m}})$ equality holds if and only if $G \cong D_8$. 



For $m \geq 10$, using Theorems \ref{d} and \ref{Dihedral2}, we have 
\begin{equation}\label{E5}
	LE(\Gamma_{D_{2m}})-LE^+(\Gamma_{D_{2m}})=\frac{4m^2-10m+6}{m-1}-2\sqrt{2m^2-8m+9}      
\end{equation}
and
\begin{equation}\label{E6}
	LE^+(\Gamma_{D_{2m}})-E(\Gamma_{D_{2m}})=\frac{m^3-9m^2+19m-8}{m-1}+2\sqrt{2m^2-8m+9}-\sqrt{5m^2-12m+4}.
\end{equation} 
Since $4m^2-10m+6 > 0$, $2(m-1)\sqrt{2m^2-8m+9} > 0$ and $(4m^2-10m+6)^2-\left(2\sqrt{2m^2-8m+9}\right)^2(m-1)^2=8m^3(m-4)+8m(5m-2)>0$
 we have $4m^2-10m+6 - 2(m-1)\sqrt{2m^2-8m+9} > 0$. Therefore, by \eqref{E5}, $(m-1)(LE(\Gamma_{D_{2m}})-LE^+(\Gamma_{D_{2m}})) > 0$. Hence, $LE(\Gamma_{D_{2m}}) > LE^+(\Gamma_{D_{2m}})$. 

Again, $2\sqrt{2m^2-8m+9} > 0, \sqrt{5m^2-12m+4} > 0$ and $\left(2\sqrt{2m^2-8m+9}\right)^2-\left(\sqrt{5m^2-12m+4}\right)^2=m(3m-10)+8 > 0$. Therefore, 
$2\sqrt{2m^2-8m+9}-\sqrt{5m^2-12m+4} > 0$. Since $m^3-9m^2+19m-8 > 0$ we have $\frac{m^3-9m^2+19m-8}{m-1}+2\sqrt{2m^2-8m+9}-\sqrt{5m^2-12m+4} > 0$. Therefore, by \eqref{E6}, $LE^+(\Gamma_{D_{2m}}) >  E(\Gamma_{D_{2m}})$. Hence, 
$E(\Gamma_{D_{2m}}) < LE^+(\Gamma_{D_{2m}}) < LE(\Gamma_{D_{2m}})$. 




\vspace{.5cm}

\noindent (b)  \textbf{Case 1:} $m$ is odd

Here, $|v(\Gamma_{D_{2m}})|=2m-1$ and $E(K_{|v(\Gamma_{D_{2m}})|})=LE(K_{|v(\Gamma_{D_{2m}})|})=LE^+(K_{|v(\Gamma_{D_{2m}})|})=4m-4$. Using Theorem \ref{d}, we have
\begin{equation}\label{E7}
E(\Gamma_{D_{2m}})-|v(\Gamma_{D_{2m}})|= \sqrt{(m-1)(5m-1)}-m   
\end{equation}
and
\begin{equation}\label{E8}
E(K_{|v(\Gamma_{D_{2m}})|})-E(\Gamma_{D_{2m}})= 3(m-1)-\sqrt{(m-1)(5m-1)}.    
\end{equation} 
Since $\sqrt{(m-1)(5m-1)} > 0$, $m > 0$ and $\left(\sqrt{(m-1)(5m-1)}\right)^2-m^2=4m^2-6m+1>0$
we have $ \sqrt{(m-1)(5m-1)}-m> 0$. Therefore, by \eqref{E7}, $E(\Gamma_{D_{2m}}) > |v(\Gamma_{D_{2m}})|$.

Again,  $\sqrt{(m-1)(5m-1)} > 0$, $3(m-1) > 0$ and $(3(m-1))^2-\left(\sqrt{(m-1)(5m-1)}\right)^2=4(m^2-3m+2)>0$
and so $3(m-1) - \sqrt{(m-1)(5m-1)}> 0$. Therefore, by \eqref{E8}, $E(K_{|v(\Gamma_{D_{2m}})|}) > E(\Gamma_{D_{2m}})$.


\vspace{.5cm}

\noindent \textbf{Case 2:} $m$ is even

Here, $|v(\Gamma_{D_{2m}})|=2m-2$ and $E(K_{|v(\Gamma_{D_{2m}})|})=LE(K_{|v(\Gamma_{D_{2m}})|})=LE^+(K_{|v(\Gamma_{D_{2m}})|})=4m-6$. Using Theorem \ref{d}, we have
\begin{equation}\label{E9}
	E(\Gamma_{D_{2m}})-|v(\Gamma_{D_{2m}})|= \sqrt{(m-2)(5m-2)}-m   
\end{equation}
and
\begin{equation}\label{E10}
	E(K_{|v(\Gamma_{D_{2m}})|})-E(\Gamma_{D_{2m}})= 3(m-2)+2-\sqrt{(m-2)(5m-2)}.    
\end{equation}
Since $\sqrt{(m-2)(5m-2)} > 0$, $m > 0$ and $\left(\sqrt{(m-2)(5m-2)}\right)^2-m^2=4(m^2-3m+1)>0$
we have $ \sqrt{(m-2)(5m-2)}-m> 0$. Therefore, by \eqref{E9}, $E(\Gamma_{D_{2m}}) > |v(\Gamma_{D_{2m}})|$.

Again,  $\sqrt{(m-2)(5m-2)} > 0$, $3(m-2)+2 > 0$ and 
$$
(3(m-2)+2)^2-\left(\sqrt{(m-2)(5m-2)}\right)^2=4(m^2-3m+3)>0
$$
and so $3(m-2)+2 - \sqrt{(m-2)(5m-2)}> 0$. Therefore, by \eqref{E10}, $E(K_{|v(\Gamma_{D_{2m}})|}) > E(\Gamma_{D_{2m}})$.


\vspace{0.5cm}

\noindent (c) \textbf{Case 1:} $m$ is odd

For $m=3$, using Theorems \ref{d} and \ref{Dihedral1}, $LE(\Gamma(D_6))=\frac{42}{5}$, $LE^+(\Gamma(D_6))=\frac{9}{5}+\sqrt{33}$ and $LE^+(K_{|v(\Gamma(D_6))}|)$ $= LE(K_{|v(\Gamma(D_6))}|) =8$. Clearly, 
$$
LE^+(\Gamma(D_6)) < LE^+(K_{|v(\Gamma(D_6))}|) = LE(K_{|v(\Gamma(D_6))}|) < LE(\Gamma(D_6)).
$$

 For $m\geq 5$, using Theorem \ref{Dihedral1}, we have
$$
LE^+(\Gamma_{D_{2m}})-LE^+(K_{|v(\Gamma_{D_{2m}})|})= \frac{2m^2(m-9)+24m-7}{2m-1}+\sqrt{8m^2-16m+9}>0.
$$
 Therefore, $LE^+(\Gamma_{D_{2m}})>LE^+(K_{|v(\Gamma_{D_{2m}})|})$ which implies $\Gamma_{D_{2m}}$ is Q-hyperenergetic and consequently part (a) implies $\Gamma_{D_{2m}}$ is L-hyperenergetic. 

\vspace{0.5cm}

\noindent \textbf{Case 2:} $m$ is even

For $m=4$, using Theorem \ref{d}, we have $LE(\Gamma(D_8))=8$ and $LE(K_{|v(\Gamma(D_8))}|)=10$. Clearly, $LE(\Gamma(D_8))$ $ < LE(K_{|v(\Gamma(D_8))}|)$. Therefore, $\Gamma_{D_8}$ is not L-hyperenergetic and  not Q-hyperenergetic.

 Using Theorem \ref{Dihedral2}, for $m = 6$ and $8$, we have   
 $$
 LE^+(\Gamma_{D_{2m}})-LE^+(K_{|v(\Gamma_{D_{2m}})|})= \frac{m^2(m-12)+20m}{2m-2}+2\sqrt{2m^2-8m+9}>0
 $$  
 and for $m\geq 10$,
$LE^+(\Gamma_{D_{2m}})-LE^+(K_{|v(\Gamma_{D_{2m}})|})= \frac{m^2(m-12)+26m-9}{m-1}+2\sqrt{2m^2-8m+9}>0$.
 Therefore, $LE^+(\Gamma_{D_{2m}})>LE^+(K_{|v(\Gamma_{D_{2m}})|})$ which implies $\Gamma_{D_{2m}}$ is Q-hyperenergetic and consequently part (a) implies $\Gamma_{D_{2m}}$ is L-hyperenergetic.
\end{proof}


In Theorem \ref{D_{2m}}, we compare $E(\Gamma_{D_{2m}})$, $LE(\Gamma_{D_{2m}})$ and $LE^+(\Gamma_{D_{2m}})$. However, in the following figures, we show how close are they.

\vspace{.3cm}

\begin{minipage}[t]{.5\linewidth}
\begin{tikzpicture}
\begin{axis}
[
xlabel={$m$ $\rightarrow$},
ylabel={Energies of $\Gamma_{D_{2m}} \rightarrow$},
xmin=2, xmax=21,
ymin=0, ymax=438,
grid = both,
minor tick num = 1,
major grid style = {lightgray},
minor grid style = {lightgray!25},
width=.7\textwidth,
height=.7\textwidth,
legend style={legend pos=north west},
 ]
\addplot[domain=3:22,samples at={3,5,7,9,11,13,15,17,19,21},mark=*,green, samples=18, mark size=.8pt]{(x-1)+(5*x*x-6*x+1)^(1/2)};
\tiny
\addlegendentry{$E$}
\addplot[domain=3:22,samples at={3,5,7,9,11,13,15,17,19,21},mark=triangle*,blue,mark size=.8pt, samples=18]{(2*x*(x-1)*(x-2)+2*x*(2*x-1))/(2*x-1)};
\tiny
\addlegendentry{$LE$}
\addplot[domain=5:22,samples at={5,7,9,11,13,15,17,19,21},mark=square*, red, mark size=.8pt, samples=18]{(2*x*x*x-10*x*x+12*x-3)/(2*x-1)+(8*x*x-16*x+9)^(1/2)};
\tiny
\addlegendentry{$LE^+$}
\addplot[domain=3:4,samples at={3},mark=square*, red, mark size=.8pt, samples=18]{(9)/(5)+(33)^(1/2)};
\tiny
\end{axis}
\end{tikzpicture}
\vspace{-.2 cm}
\captionsetup{font=footnotesize}

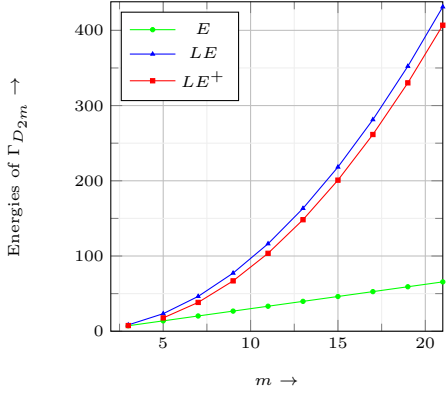
\captionof{figure}{Energies of $\Gamma_{D_{2m}}$, $m$ is odd} 
\end{minipage}
\hspace{0.05cm}
\begin{minipage}[t]{.5\linewidth}
\begin{tikzpicture}
\begin{axis}
[
xlabel={$m$ $\rightarrow$},
ylabel={Energies of $\Gamma_{D_{2m}} \rightarrow$},
xmin=2, xmax=20,
ymin=0, ymax=370,
grid = both,
minor tick num = 1,
major grid style = {lightgray},
minor grid style = {lightgray!25},
width=.7\textwidth,
height=.7\textwidth,
legend style={legend pos=north west},
 ]
\addplot[domain=10:50,samples at={4,6,8,10,12,14,16,18,20},mark=*,green, samples=20, mark size=.8pt]{(x-2)+(5*x*x-12*x+4)^(1/2)};
\tiny
\addlegendentry{$E$}
\addplot[domain=10:50,samples at={4,6,8,10,12,14,16,18,20},mark=triangle*,blue,mark size=.8pt, samples=20]{(x*(x-2)*(x-4)+2*x*(x-1))/(x-1)};
\tiny
\addlegendentry{$LE$}
\addplot[domain=10:50,samples at={10,12,14,16,18,20},mark=square*, red, mark size=.8pt, samples=20]{(x*x*x-8*x*x+16*x-3)/(x-1)+2*(2*x*x-8*x+9)^(1/2)};
\tiny
\addlegendentry{$LE^+$}
\addplot[domain=4:8,samples at={4,6,8},mark=square*, red, mark size=.8pt, samples=20]{(x*x*x-4*x*x+12)/(2*x-2)+2*(2*x*x-8*x+9)^(1/2)};
\tiny
\end{axis}
\end{tikzpicture}
\vspace{-.2 cm}
\captionsetup{font=footnotesize}

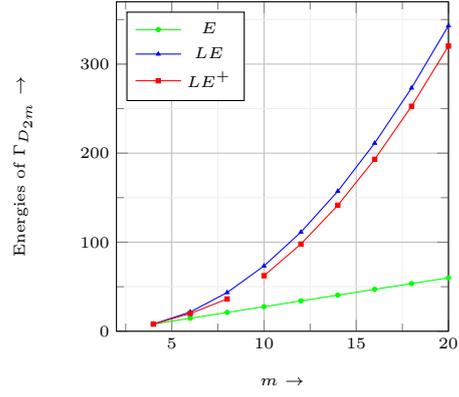
\captionof{figure}{Energies of $\Gamma_{D_{2m}}$, $m$ is even}
\end{minipage}

\subsection{The Quasidihedral groups, $QD_{2^n}$}
We consider $QD_{2^n}:=\langle a,b : a^{2^{n-1}} = b^2 = 1, bab^{-1} = a^{2^{n-2}-1} \rangle$, the quasidihedral groups of order $2^n$ (where $n \geq 4$). Results regarding different energies of non-commuting graphs of $QD_{2^n}$ are given below.

\begin{theorem} [\protect{\cite[Result 1.2.16 (d)]{FWNT21}and \cite[Proposition 3.2]{DN18}}]\label{QD}
Let $G$ be isomorphic to $QD_{2^n}$. Then
\[E(\Gamma_{QD_{2^n}})= (2^{n-1}-2)+2\sqrt{(5\times 2^{n-2}-1)(2^{n-2}-1)} \, \text{ and } \,
LE(\Gamma_{QD_{2^n}}) = \frac{2^{3n-3}-2^{2n}+3\times 2^n}{2^{n-1}-1}.\]
\end{theorem}

\begin{theorem}\label{Quasidihedral}
Let $G$ be isomorphic to $QD_{2^n}$. Then 
\begin{align*}
\Q-spec(\Gamma_{QD_{2^n}})=&\left\lbrace(2^n-4)^{2^{n-2}}, (2^n-2^{n-1})^{2^{n-1}-3},(2^n-6)^{2^{n-2}-1},\left(2^n-3+\sqrt{2^{2n-1}-2^{n+2}+9}\right)^1,\right.\\
&\left.\left(2^n-3-\sqrt{2^{2n-1}-2^{n+2}+9}\right)^1\right\rbrace
\end{align*}
and
\begin{align*}
LE^+(\Gamma_{QD_{2^n}}) = \begin{cases}
\frac{134}{7}+2\sqrt{73}, & \mbox{if $n=4$}\vspace{.2cm}\\
 \frac{2^{3n-2}+2^{n+4}-2^{2n+2}-12}{2^n-2}+2\sqrt{2^{2n-1}-2^{n+2}+9}, & \mbox{if $n \geq 5.$}
\end{cases}   
\end{align*}
\end{theorem}
\begin{proof}
If $G\cong QD_{2^n}$ then $|v(\Gamma_{QD_{2^n}})|=2^n-2$ and $\Gamma_{QD_{2^n}}=K_{2^{n-2}.2,1.(2^{n-1}-2)}$. Using Theorem \ref{R1}(b), we have 
\begin{align*}
Q_{\Gamma_{QD_{2^n}}}(x)=&{\displaystyle \prod_{i=1}^2} (x-(2^n-2)+p_i)^{a_i(p_i-1)}{\displaystyle \prod_{i=1}^2}(x-(2^n-2)+2p_i)^{a_i}\left(1-{\displaystyle \sum_{i=1}^2}\frac{a_ip_i}{x-(2^n-2)+2p_i}\right)\\
  = \,& (x-2^n+4)^{2^{n-2}}(x-2^n+2^{n-1})^{2^{n-1}-3}(x-2^n+6)^{2^{n-2}}(x-2) \left(1-\frac{2^{n-1}}{x-2^n+6}-\frac{2^{n-1}-2}{x-2}\right)\\
   =\,& (x-(2^n-4))^{2^{n-2}}(x-(2^n-2^{n-1}))^{2^{n-1}-3}(x-(2^n-6))^{2^{n-2}-1} (x^2-(2^{n+1}-6)x+2^{2n-1}-2^{n+1}).
\end{align*}
Thus, \begin{align*}
\Q-spec(\Gamma_{QD_{2^n}})=&\left\lbrace(2^n-4)^{2^{n-2}}, (2^n-2^{n-1})^{2^{n-1}-3},(2^n-6)^{2^{n-2}-1},\left(2^n-3+\sqrt{2^{2n-1}-2^{n+2}+9}\right)^1,\right.\\
&\left.\left(2^n-3-\sqrt{2^{2n-1}-2^{n+2}+9}\right)^1\right\rbrace.
\end{align*}

Number of edges of $\Gamma_{QD_{2^n}}^c$ is $2^{2n-3}-2^n+3$. Thus, $|e(\Gamma_{QD_{2^n}})| = \frac{(2^n-2)(2^n-2-1)}{2}-(2^{2n-3}-2^n+3)$ $=3(2^{2n-3}-2^{n-1})$. Now

\[\left|2^n-4 - \frac{2|e(\Gamma_{QD_{2^n}})|}{|v(\Gamma_{QD_{2^n}})|}\right| = \left|\frac{8+2^{2n-2}-3\times 2^{n}}{2^n-2}\right|= \frac{8+2^{2n-2}-3\times 2^{n}}{2^n-2},\] 

\[\left|2^n-2^{n-1} - \frac{2|e(\Gamma_{QD_{2^n}})|}{|v(\Gamma_{QD_{2^n}})|}\right| = \left|\frac{2^{n+1}-2^{2n-2}}{2^n-2}\right|  =\frac{2^{2n-2}-2^{n+1}}{2^n-2},\]

\[\left|2^n-6 - \frac{2|e(\Gamma_{QD_{2^n}})|}{|v(\Gamma_{QD_{2^n}})|}\right| = \left|\frac{12+2^{2n-2}-5\times 2^{n}}{2^n-2}\right|= \begin{cases}\frac{2}{7}, & \mbox{if $n=4$}\vspace{.2cm}\\
\frac{12+2^{2n-2}-5\times 2^{n}}{2^n-2}, & \mbox{if $n \geq 5$,}
 \end{cases}\]

 \begin{align*}
\left|2^n-3+\sqrt{2^{2n-1}-2^{n+2}+9} - \frac{2|e(\Gamma_{QD_{2^n}})|}{|v(\Gamma_{QD_{2^n}})|}\right|
= &\left|\sqrt{2^{2n-1}-2^{n+2}+9}+\frac{2^{2n-2}-2^{n+1}+6}{2^n-2}\right|\\
=&\sqrt{2^{2n-1}-2^{n+2}+9}+\frac{2^{2n-2}-2^{n+1}+6}{2^n-2}\end{align*}

and
\begin{align*}
\left|2^n-3-\sqrt{2^{2n-1}-2^{n+2}+9} - \frac{2|e(\Gamma_{QD_{2^n}})|}{|v(\Gamma_{QD_{2^n}})|}\right|
= &\left|-\sqrt{2^{2n-1}-2^{n+2}+9}+\frac{2^{2n-2}-2^{n+1}+6}{2^n-2}\right|\\
=&\sqrt{2^{2n-1}-2^{n+2}+9}-\frac{2^{2n-2}-2^{n+1}+6}{2^n-2}.\end{align*}
Therefore, for $n=4$ we have  $LE^+(\Gamma_{QD_{2^n}})= \frac{134}{7}+2\sqrt{73}$. For $n \geq 5$ we have 
\begin{align*}
LE^+(\Gamma_{QD_{2^n}}) = &(2^{n-2}) \times \frac{8+2^n(2^{n-2}-3)}{2^n-2} + (2^{n-1}-3) \times  \frac{2^{n+1}(2^{n-3}-1)}{2^n-2}+ (2^{n-2}-1) \times \frac{12+2^n(2^{n-2}-5)}{2^n-2}\\
&+\sqrt{2^{2n-1}-2^{n+2}+9}+\frac{2^{2n-2}-2^{n+1}+6}{2^n-2}+\sqrt{2^{2n-1}-2^{n+2}+9}-\frac{2^{2n-2}-2^{n+1}+6}{2^n-2}
\end{align*} and the result follows on simplification.
\end{proof}

\begin{theorem}\label{QDn}
If $G$ is isomorphic to $QD_{2^n}$ then 
\begin{enumerate}
\item $E(\Gamma_{QD_{2^n}})< LE^+(\Gamma_{QD_{2^n}}) < LE(\Gamma_{QD_{2^n}})$.
\item $\Gamma_{QD_{2^n}}$ is non-hypoenergetic as well as non-hyperenergetic.
\item $\Gamma_{QD_{2^n}}$ is Q-hyperenergetic and L-hyperenergetic.
\end{enumerate}
\end{theorem}

\begin{proof}
(a) For $n = 4$, using Theorems \ref{QD} and \ref{Quasidihedral}, we have
$E(\Gamma_{QD_{2^n}})=6+2\sqrt{57}$, $LE(\Gamma_{QD_{2^n}})=\frac{304}{7}$ and $LE^+(\Gamma_{QD_{2^n}})=\frac{134}{7}+2\sqrt{73}$. Clearly, $E(\Gamma_{QD_{16}})< LE^+(\Gamma_{QD_{16}}) < LE(\Gamma_{QD_{16}})$. 

For $n \geq 5$, using Theorems \ref{QD} and \ref{Quasidihedral},  we have
\begin{equation}\label{E30}
LE(\Gamma_{QD_{2^n}})-LE^+(\Gamma_{QD_{2^n}})= \frac{12+2^{2n+1}-5\times 2^{n+1}}{2^n-2}-2\sqrt{2^{2n-1}-2^{n+2}+9}
\end{equation}
and
\begin{equation}\label{E31}
LE^+(\Gamma_{QD_{2^n}})-E(\Gamma_{QD_{2^n}})= \frac{2^{2n-2}(2^n-18)+ 19\times 2^n-16}{2^n-2}+2\sqrt{2^{2n-1}-2^{n+2}+9}-2\sqrt{5\times 2^{2n-4}-3\times 2^{n-1}+1}.
\end{equation} 
Since $12+2^{2n+1}-5 \times 2^{n+1} > 0$, $2\sqrt{2^{2n-1}-2^{n+2}+9}(2^n-2) > 0$ and
\begin{center}
 $(12+2^{2n+1}-5 \times 2^{n+1})^2-\left(2\sqrt{2^{2n-1}-2^{n+2}+9}\right)^2(2^n-2)^2=2^{3n+1}(2^n-8)+2^{n+3}(5\times 2^n-4)>0$
\end{center}
we have $12+2^{2n+1}- 5 \times 2^{n+1}-2(2^n-2)\sqrt{2^{2n-1}-2^{n+2}+9}>0$. Therefore, by \eqref{E30}, $(2^n-2)(LE(\Gamma_{QD_{2^n}})-LE^+(\Gamma_{QD_{2^n}})) > 0$. Hence, $LE(\Gamma_{QD_{2^n}})>LE^+(\Gamma_{QD_{2^n}})$.  

Again, $\sqrt{2^{2n-1}-2^{n+2}+9} > 0, \sqrt{5\times 2^{2n-4}-3\times 2^{n-1}+1} > 0$ and
\begin{center}
 $\left(\sqrt{2^{2n-1}-2^{n+2}+9}\right)^2-\left(\sqrt{5\times2^{2n-4}-3\times 2^{n-1}+1}\right)^2=2^{n-4}(3\times 2^n-40)+8>0$.
\end{center}
 Therefore, 
 we have  $\sqrt{2^{2n-1}-2^{n+2}+9}-\sqrt{5\times 2^{2n-4}-3\times 2^{n-1}+1} > 0$. Since $2^{2n-2}(2^n-18)+ 19 \times 2^n-16 > 0$ we have $\frac{2^{2n-2}(2^n-18)+19 \times 2^n-16}{2^n-2}+2\sqrt{2^{2n-1}-2^{n+2}+9}-2\sqrt{5\times 2^{2n-4}-3\times2^{n-1}+1} > 0$. Therefore, by \eqref{E31}, $LE^+(\Gamma_{QD_{2^n}}) \geq  E(\Gamma_{QD_{2^n}})$. Hence, 
$E(\Gamma_{QD_{2^n}}) < LE^+(\Gamma_{QD_{2^n}}) < LE(\Gamma_{QD_{2^n}})$. 


\vspace{.5cm}
\noindent (b) Here, $|v(\Gamma_{QD_{2^n}})|=2^n-2$ and $E(K_{|v(\Gamma_{QD_{2^n}})|})= LE(K_{|v(\Gamma_{QD_{2^n}})|})= LE^+(K_{|v(\Gamma_{QD_{2^n}})|})= 2^{n+1}-6$. Using Theorem \ref{QD}, we have
\begin{equation}\label{E32}
E(\Gamma_{QD_{2^n}})-|v(\Gamma_{QD_{2^n}})|= 2\left(\sqrt{(5\times 2^{n-2}-1)(2^{n-2}-1)}-(2^{n-1}-2^{n-2})\right) 
\end{equation} 
and
\begin{equation}\label{E33}
E(K_{|v(\Gamma_{QD_{2^n}})|})-E(\Gamma_{QD_{2^n}})= 2\left(3\times 2^{n-2}-2-\sqrt{(5\times 2^{n-2}-1)(2^{n-2}-1)}\right).
\end{equation} 
Since $\sqrt{(5\times 2^{n-2}-1)(2^{n-2}-1)} > 0$, $2^{n-1}-2^{n-2} > 0$ and $\left(\sqrt{(5\times 2^{n-2}-1)(2^{n-2}-1)}\right)^2-(2^{n-1}-2^{n-2})^2=2^{n-2}(4\times 2^{n-2}-6)+1>0$
we have $ \sqrt{(5\times 2^{n-2}-1)(2^{n-2}-1)}-(2^{n-1}-2^{n-2})> 0$. Therefore, by \eqref{E32}, $E(\Gamma_{QD_{2^n}}) > |v(\Gamma_{QD_{2^n}})|$.

Again,  $\sqrt{(5\times 2^{n-2}-1)(2^{n-2}-1)} > 0$, $3\times 2^{n-2}-2 > 0$ and $(3\times 2^{n-2}-2)^2-\left(\sqrt{(5\times 2^{n-2}-1)(2^{n-2}-1)}\right)^2=(2^{n-2}-3)(4\times 2^{n-2}+6)+21>0$
and so $3\times 2^{n-2}-2 - \sqrt{(5\times 2^{n-2}-1)(2^{n-2}-1)}> 0$. Therefore, by \eqref{E33}, $E(K_{|v(\Gamma_{QD_{2^n}})|}) > E(\Gamma_{QD_{2^n}})$.


\vspace{.5cm}
\noindent (c) For $n = 4$, using Theorem \ref{Quasidihedral},
$LE^+(\Gamma_{QD_{16}})-LE^+(K_{|v(\Gamma_{QD_{16}})|})= 2\sqrt{73}-\frac{48}{7}>0.$  Therefore, $LE^+(\Gamma_{QD_{16}})>LE^+(K_{|v(\Gamma_{QD_{16}})|})$ which implies $\Gamma_{QD_{16}}$ is Q-hyperenergetic and consequently part (a) implies $\Gamma_{QD_{16}}$ is L-hyperenergetic.

For $n \geq 5$, using Theorem \ref{Quasidihedral}, 
\begin{center}
	$LE^+(\Gamma_{QD_{2^n}})-LE^+(K_{|v(\Gamma_{QD_{2^n}})|})= \frac{2^{2n-2}(2^n-24)+2(13 \times 2^n-12)}{2^n-2}+2\sqrt{2^{2n-1}-2^{n+2}+9}>0.$
\end{center}
 Therefore, $LE^+(\Gamma_{QD_{2^n}})>LE^+(K_{|v(\Gamma_{QD_{2^n}})|})$ which implies $\Gamma_{QD_{2^n}}$ is Q-hyperenergetic and consequently part (a) implies $\Gamma_{QD_{2^n}}$ is L-hyperenergetic.
\end{proof}

In Theorem \ref{QDn}, we compare $E(\Gamma_{QD_{2^n}})$, $LE(\Gamma_{QD_{2^n}})$ and $LE^+(\Gamma_{QD_{2^n}})$. However, in the following figures, we show how close are they.

\vspace{.3cm}

\begin{minipage}[t]{.5\linewidth}
\begin{tikzpicture}
\begin{axis}
[
xlabel={$n$ $\rightarrow$},
ylabel={Energies  of $\Gamma_{QD_{2^n}}$ $\rightarrow$},
xmin=3, xmax=9,
ymin=0, ymax=20000,
grid = both,
minor tick num = 1,
major grid style = {lightgray},
minor grid style = {lightgray!25},
width=.7\textwidth,
height=.7\textwidth,
legend style={legend pos=north west},
 ]
\addplot[domain=4:9,samples at={4,5,6,7,8,9},mark=*,green, samples=6, mark size=.8pt]{(2^(x-1)-2)+2*(5*2^(2*x-4)-3*2^(x-1)+1)^(1/2)};
\tiny
\addlegendentry{$E$}
\addplot[domain=4:9,samples at={4,5,6,7,8,9},mark=triangle*,blue,mark size=.8pt, samples=6]{(2^(3*x-3)-2^(2*x)+3*2^x)/(2^(x-1)-1)};
\tiny
\addlegendentry{$LE$}
\addplot[domain=5:9,samples at={5,6,7,8,9},mark=square*, red, mark size=.8pt, samples=6]{(2^(3*x-2)+2^(x+4)-2^(2*x+2)-12)/(2^x-2)+2*(2^(2*x-1)-2^(x+2)+9)^(1/2)};
\tiny
\addlegendentry{$LE^+$}
\addplot[domain=4:4.5,samples at={4},mark=square*, red, mark size=.8pt, samples=6]{(134)/(7)+2*(73)^(1/2)};
\tiny
\end{axis}
\end{tikzpicture}
\vspace{-.2 cm}
\captionsetup{font=footnotesize}

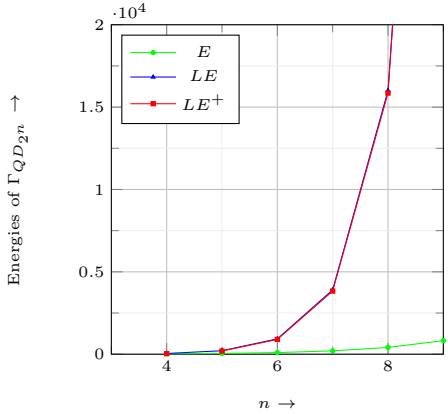
\captionof{figure}{Energies of $\Gamma_{QD_{2^n}}$}
\end{minipage}
\hspace{0.05cm}
\begin{minipage}[t]{.5\linewidth}
\includegraphics[width=6cm]{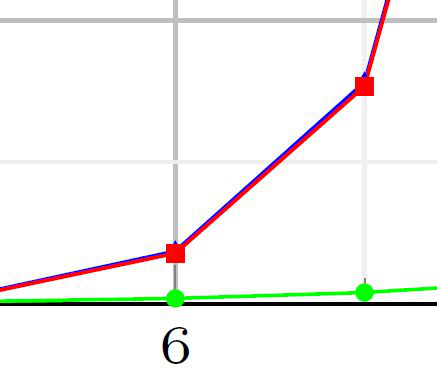}
\vspace{-.2 cm}
\captionsetup{font=footnotesize}
{\captionof{figure}{A close up view of Figure 3}}
\end{minipage}

\subsection{The groups $M_{2rs}$}
We consider $M_{2rs}:=\langle a,b:a^r=b^{2s}=1, bab^{-1}=a^{-1}\rangle$, the groups of order $2rs$ (where $r \geq 3$ and $s \geq 1$). Results regarding different energies of non-commuting graphs of $M_{2rs}$ are given below.
\begin{theorem}[\protect{\cite[Corollary 4.1.6  and (4.3.d)]{FWNT21}}]\label{M}
Let $G$ be isomorphic to $M_{2rs}$.
\begin{enumerate}
\item If $m$ is odd then
\[E(\Gamma_{M_{2rs}}) = s(r-1)+s\sqrt{(r-1)(5r-1)} \, \text{ and } \,
LE(\Gamma_{M_{2rs}}) = \frac{s}{2r-1}\left(2r^3s-6r^2s+4rs+4r^2-2r\right).\]
\item If $m$ is even then
\[E(\Gamma_{M_{2rs}}) = s(r-2)+s\sqrt{(r-2)(5r-2)} \, \text{ and } \,
LE(\Gamma_{M_{2rs}}) = \frac{s}{r-1}\left(r^3s-6r^2s+8rs+2r^2-2r\right).\]
\end{enumerate}
\end{theorem}

\begin{theorem}\label{M_2rs1}
Let $G$ be isomorphic to $M_{2rs}$, where  $r$ is odd. Then 
\begin{align*}
\Q-spec(\Gamma_{M_{2rs}})=&\left\{(2s(r-1))^{r(s-1)},(rs)^{(r-1)s-1},((2r-3)s)^{r-1}, \left(\frac{s\left(4r-3+\sqrt{8r^2-16r+9}\right)}{2}\right)^1,\right.\\ &\left.\left(\frac{s\left(4r-3+\sqrt{8r^2-16r+9}\right)}{2}\right)^1\right\}   
\end{align*} and 
\begin{align*}
LE^+(\Gamma_{M_{2rs}}) =&\begin{cases}
\frac{3s(4s-1)}{5}+s\sqrt{33}, & \mbox{if $r=3$}\vspace{.2cm}\\
s\left(\frac{(2r(r-1)(r-2))s}{2r-1}-(2r-3)+\sqrt{8r^2-16r+9}\right), & \mbox{if $r \geq 5.$}
\end{cases}   
\end{align*}
\end{theorem}

\begin{proof}
If $G\cong M_{2rs}$, where $r$ is odd, then $|v(\Gamma_{M_{2rs}})|=(2r-1)s$ and $\Gamma_{M_{2rs}}=K_{r.s,1.((r-1)s)}$. Using Theorem \ref{R1}(b), we have 
\begin{align*}
Q_{\Gamma_{M_{2rs}}}(x)=&{\displaystyle \prod_{i=1}^2} (x-(2rs-s)+p_i)^{a_i(p_i-1)}{\displaystyle \prod_{i=1}^2}(x-(2rs-s)+2p_i)^{a_i}\left(1-{\displaystyle \sum_{i=1}^2}\frac{a_ip_i}{x-(2rs-s)+2p_i}\right)\\
= \,& (x-2s(r-1))^{r(s-1)}(x-rs)^{(r-1)s-1}(x-(2r-3)s)^{r}(x-s)\left(1-\frac{rs}{x-(2r-3)s}-\frac{(r-1)s}{x-s}\right)\\
= \,& (x-2s(r-1))^{r(s-1)}(x-rs)^{(r-1)s-1}(x-(2r-3)s)^{r-1}(x^2-(4r-3)sx+(2r^2-2r)s^2).
\end{align*}
Thus, $\Q-spec(\Gamma_{M_{2rs}})=\left\{(2s(r-1))^{r(s-1)},(rs)^{(r-1)s-1},((2r-3)s)^{r-1}, \left(\frac{s(4r-3+\sqrt{8r^2-16r+9})}{2}\right)^1,\right.\\ 
 ~~~~~~~~~~~~~~~~~~~~~~~~~~~~~~~~~\left.\left(\frac{s(4r-3-\sqrt{8r^2-16r+9})}{2}\right)^1\right\}$.
 
 Number of edges of $\Gamma_{M_{2rs}}^c$ is $\frac{(r^2-r+1)s^2-(2r-1)s}{2}$. Therefore, \begin{center}
 	$|e(\Gamma_{M_{2rs}})| = \frac{(2r-1)^2s^2-(2r-1)s}{2}-\frac{(r^2-r+1)s^2-(2r-1)s}{2}=\frac{3r(r-1)s^2}{2}$.
 \end{center} 
Now
\[\left|(2r-2)s- \frac{2|e(\Gamma_{M_{2rs}})|}{|v(\Gamma_{M_{2rs}})|}\right| = \left|\frac{(r-1)(r-2)s}{2r-1}\right| = \frac{(r-1)(r-2)s}{2r-1},\]

\[\left|rs - \frac{2|e(\Gamma_{M_{2rs}})|}{|v(\Gamma_{M_{2rs}})|}\right| = \left|\frac{-r(r-2)s}{2r-1}\right| = \frac{r(r-2)s}{2r-1},\]

\[\left|(2r-3)s - \frac{2|e(\Gamma_{M_{2rs}})|}{|v(\Gamma_{M_{2rs}})|}\right| = 
\left|\frac{(r^2-5r+3)s}{2r-1}\right|= \begin{cases} \frac{3s}{5}, & \mbox{if $r = 3$}\vspace{.2cm}\\
 \frac{(r^2-5r+3)s}{2r-1}, & \mbox{if $r \geq 5,$}\end{cases}\]

\begin{align*}
\left|\frac{s}{2}\left(4r-3+\sqrt{8r^2-16r+9}\right) - \frac{2|e(\Gamma_{M_{2rs}})|}{|v(\Gamma_{M_{2rs}})|}\right|= & \left|\frac{s}{2}\left(\sqrt{8r^2-16r+9}+r-\frac{3}{2}+\frac{3}{4r-2}\right)\right|\\=&\frac{s}{2}\left(\sqrt{8r^2-16r+9}+r-\frac{3}{2}+\frac{3}{4r-2}\right)\end{align*} 
and
\begin{align*}
\left|\frac{s}{2}\left(4r-3-\sqrt{8r^2-16r+9}\right) - \frac{2|e(\Gamma_{M_{2rs}})|}{|v(\Gamma_{M_{2rs}})|}\right|= & \left|\frac{s}{2}\left(-\sqrt{8r^2-16r+9}+r-\frac{3}{2}+\frac{3}{4r-2}\right)\right|\\=& \frac{s}{2}\left(\sqrt{8r^2-16r+9}-r+\frac{3}{2}-\frac{3}{4r-2}\right).
\end{align*}


Therefore, for $n=3$, we have  $LE^+(\Gamma_{M_{2rs}})= \frac{3s(4s-1)}{5}+s\sqrt{33}$. For $r \geq 5$, we have
\begin{align*}
LE^+(\Gamma_{M_{2rs}}) =& r(s-1)\times \frac{(r-1)(r-2)s}{2r-1}+((r-1)s-1)\times \frac{r(r-2)s}{2r-1}+ (r-1)\times \frac{(r^2-5r+3)s}{2r-1}+\\
&\frac{s}{2}\left(\sqrt{8r^2-16r+9}+r-\frac{3}{2}+\frac{3}{4r-2}\right)+\frac{s}{2}\left(\sqrt{8r^2-16r+9}-r+\frac{3}{2}-\frac{3}{4r-2}\right)
\end{align*}
and the result follows on simplification.
\end{proof}

\begin{theorem}\label{M_2rs2}
Let $G$ be isomorphic to $M_{2rs}$, where $r$ is even. Then 
\begin{align*}
\Q-spec(\Gamma_{M_{2rs}})=&\left\lbrace(2s(r-2))^{rs-\frac{r}{2}},(rs)^{rs-2s-1},(2s(r-3))^{\frac{r}{2}-1}, \left(4rs-6s+2s\sqrt{2r^2-8r+9}\right)^1, \right.\\ &\left.\left(4rs-6s-2s\sqrt{2r^2-8r+9}\right)^1\right\rbrace  
\end{align*}
 and 
\begin{align*}
LE^+(\Gamma_{M_{2rs}}) = \begin{cases}\frac{s(r^3s-6r^2s+8rs-\frac{r^3}{2}+4r^2-8r+6)}{r-1}+2s\sqrt{2r^2-8r+9}, &  \mbox{if $4 \leq r \leq 8$}\vspace{.2cm}\\
\frac{s(r^3s-6r^2s+8rs-2r^2+8r-6)}{r-1}+2s\sqrt{2r^2-8r+9}, & \mbox{if $r \geq 10.$}
\end{cases}   
\end{align*}
\end{theorem}

\begin{proof}
If $G\cong M_{2rs}$ and $r$ is even then $|v(\Gamma_{M_{2rs}})|=2s(r-1)$ and $\Gamma_{M_{2rs}}=K_{\frac{r}{2}.(2s),1.((\frac{r}{2}-1)2s)}$. Using Theorem \ref{R1}(b), we have 
\begin{align*}
Q_{\Gamma_{M_{2rs}}}(x)=&{\displaystyle \prod_{i=1}^2} (x-2(rs-s)+p_i)^{a_i(p_i-1)}{\displaystyle \prod_{i=1}^2}(x-2(rs-s)+2p_i)^{a_i}\left(1-{\displaystyle \sum_{i=1}^2}\frac{a_ip_i}{x-2(rs-s)+2p_i}\right)\\
  = \,& (x-2s(r-2))^{rs-\frac{r}{2})}(x-rs)^{rs-2s-1}(x-2s(r-3))^{\frac{r}{2}}(x-2s)\left(1-\frac{rs}{x-2s(r-3)}-\frac{rs-2s}{x-2s}\right)\\
   =\,& (x-2s(r-2))^{rs-\frac{r}{2}}(x-rs)^{rs-2s-1}(x-2s(r-3))^{\frac{r}{2}-1}(x^2-(4r-6)sx + rs^2(2r-4)).
\end{align*}
Thus, \begin{align*}
\Q-spec(\Gamma_{M_{2rs}})=&\left\lbrace(2s(r-2))^{rs-\frac{r}{2}},(rs)^{rs-2s-1},(2s(r-3))^{\frac{r}{2}-1}, \left(4rs-6s+2s\sqrt{2r^2-8r+9}\right)^1, \right.\\ &\left.\left(4rs-6s-2s\sqrt{2r^2-8r+9}\right)^1\right\rbrace. 
\end{align*}

\noindent Number of edges of $\Gamma_{M_{2rs}}^c$ is $\frac{(r^2-2r+4)s^2-2(r-1)s}{2}$. Thus, $|e(\Gamma_{M_{2rs}})| = \frac{2(r-1)s(2(r-1)s-1)}{2}-\frac{(r^2-r+1)s^2-2(r-1)s}{2}$ $=\frac{3r(r-2)s^2}{2}$. Now 
\[\left|2s(r-2)- \frac{2|e(\Gamma_{M_{2rs}})|}{|v(\Gamma_{M_{2rs}})|}\right| = \left|\frac{(r-2)(r-4)s}{2r-2}\right| = \frac{(r-2)(r-4)s}{2r-2},\] 

\[\left|rs - \frac{2|e(\Gamma_{M_{2rs}})|}{|v(\Gamma_{M_{2rs}})|}\right| = \left|\frac{-r(r-4)s}{2r-2}\right|=
\frac{r(r-4)s}{2r-2},\]

\[\left|(2r-3)s - \frac{2|e(\Gamma_{M_{2rs}})|}{|v(\Gamma_{M_{2rs}})|}\right| = \left|\frac{(r^2-10r+12)s}{2r-2}\right|= \begin{cases}\frac{(-r^2+10r-12)s}{2r-2}, & \mbox{if $4\leq r \leq 8$}\vspace{.2cm}\\
\frac{(r^2-10r+12)s}{2r-2}, & \mbox{if $r \geq 10,$}
 \end{cases}\]
\begin{align*}\left|4rs-6s+2s\sqrt{2r^2-8r+9} - \frac{2|e(\Gamma_{M_{2rs}})|}{|v(\Gamma_{M_{2rs}})|}\right| = &\left|\frac{rs}{2}-\frac{3s}{2}+\frac{3s}{2r-2}+s\sqrt{2r^2-8r+9}\right|\\ 
=& \frac{rs}{2}-\frac{3s}{2}+\frac{3s}{2r-2}+s\sqrt{2r^2-8r+9} 
\end{align*}
and
\begin{align*}\left|4rs-6s-2s\sqrt{2r^2-8r+9} - \frac{2|e(\Gamma_{M_{2rs}})|}{|v(\Gamma_{M_{2rs}})|}\right| = &\left|\frac{rs}{2}-\frac{3s}{2}+\frac{3s}{2r-2}-s\sqrt{2r^2-8r+9}\right| \\ = & -\frac{rs}{2}+\frac{3s}{2}-\frac{3s}{2r-2}+s\sqrt{2r^2-8r+9}.\end{align*}

Therefore, for $4\leq r \leq 8$, we have 
\begin{align*}
LE^+(\Gamma_{M_{2rs}}) = &\left(rs-\frac{r}{2}\right)\times \frac{(r-2)(r-4)s}{2r-2}+(rs-2s-1)\times \frac{r(r-4)s}{2r-2}+ \left(\frac{r}{2}-1\right)\times \frac{(-r^2+10r-12)s}{2r-2}\\
&+\frac{rs}{2}-\frac{3s}{2}+\frac{3s}{2r-2}+s\sqrt{2r^2-8r+9}-\frac{rs}{2}+\frac{3s}{2}-\frac{3s}{2r-2}+s\sqrt{2r^2-8r+9}    \end{align*}
 
and for $r \geq 10$, we have
\begin{align*}
 LE^+(\Gamma_{M_{2rs}}) =& \left(rs-\frac{r}{2}\right)\times \frac{(r-2)(r-4)s}{2r-2}+(rs-2s-1)\times \frac{r(r-4)s}{2r-2}+\left(\frac{r}{2}-1\right)\times \frac{(r^2-10r+12)s}{2r-2}\\
 &+\frac{rs}{2}-\frac{3s}{2}+\frac{3s}{2r-2}+s\sqrt{2r^2-8r+9}-\frac{rs}{2}+\frac{3s}{2}-\frac{3s}{2r-2}+s\sqrt{2r^2-8r+9}. \end{align*}
Hence, the results follow on simplification.
\end{proof}

\begin{theorem}\label{m}
If  $G$ is isomorphic to $M_{2rs}$ then
\begin{enumerate}
\item $E(\Gamma_{M_{2rs}}) \leq LE^+(\Gamma_{M_{2rs}}) \leq LE(\Gamma_{M_{2rs}})$, equality holds if and only if $G \cong M_{8s}$.
\item $\Gamma_{M_{2rs}}$ is non-hypoenergetic as well as non-hyperenergetic.
\item $\Gamma_{M_{6}}$ is L-hyperenergetic but not Q-hyperenergetic. $\Gamma_{M_{8s}}$ is not L-hyperenergetic and not Q-hyperenergetic. If $2rs \ne 6$ and $8s$ then $\Gamma_{M_{2rs}}$ is Q-hyperenergetic and L-hyperenergetic.
\end{enumerate}
\end{theorem}

\begin{proof} (a) \textbf{Case 1:} $r$ is odd

Using Theorems \ref{M} and \ref{M_2rs1}, for $r=3$,  $LE(\Gamma_{M_{2rs}})-LE^+(\Gamma_{M_{2rs}})=\frac{33s}{5}-s\sqrt{33}>0$ and $LE^+(\Gamma_{M_{2rs}})-E(\Gamma_{M_{2rs}})= \frac{12s^2-13s}{5} + (\sqrt{33}-2\sqrt{7})s>0$. 
For $r \geq 5$, using Theorems \ref{M} and \ref{M_2rs1}, we have 
\begin{equation}\label{m11}
LE(\Gamma_{M_{2rs}})-LE^+(\Gamma_{M_{2rs}})=\frac{s(8r^2-10r+3)}{2r-1}-s\sqrt{8r^2-16r+9}
\end{equation}
and
\begin{equation}\label{m22}
LE^+(\Gamma_{M_{2rs}})-E(\Gamma_{M_{2rs}})=s\left(\frac{(2r^3-6r^2+4r)s-6r^2+11r-4}{2r-1}+\sqrt{8r^2-16r+9}-\sqrt{5r^2-6r+1}\right).
\end{equation}

Since $8r^2-10r+3 > 0$, $(2r-1)\sqrt{8r^2-16r+9} > 0$ and $(8r^2-10r+3)^2-(2r-1)^2(8r^2-16r+9)=32r^4-64r^3+40r^2-8r>0$
we have $8r^2-10r+3 - (2r-1)\sqrt{8r^2-16r+9} > 0$. Therefore, by \eqref{m11}, $(2r-1)(LE(\Gamma_{M_{2rs}})-LE^+(\Gamma_{M_{2rs}})) > 0$. Hence, $LE(\Gamma_{M_{2rs}}) > LE^+(\Gamma_{M_{2rs}})$.  

Again, $\sqrt{8r^2-16r+9} > 0, \sqrt{5r^2-6r+1} > 0$ and $(\sqrt{8r^2-16r+9})^2-(\sqrt{5r^2-6r+1})^2$ $=r(3r-10)+8>0$. Therefore, 
$\sqrt{8r^2-16r+9}-\sqrt{5r^2-6r+1}>0$. Since $2r^3-6r^2+4r > 6r^2-11r+4$ we have $\frac{(2r^3-6r^2+4r)s-6r^2+11r-4}{2r-1}+\sqrt{8r^2-16r+9}-\sqrt{5r^2-6r+1}>0$. Therefore, by \eqref{m22}, $LE^+(\Gamma_{M_{2rs}}) >  E(\Gamma_{M_{2rs}})$. Hence, 
$E(\Gamma_{M_{2rs}}) < LE^+(\Gamma_{M_{2rs}}) < LE(\Gamma_{M_{2rs}})$.



\vspace{.5cm}
\noindent \textbf{Case 2:} $r$ is even

For $4 \leq r \leq 8$, using Theorems \ref{M} and \ref{M_2rs2}, we have 
\begin{equation}\label{m3}
LE(\Gamma_{M_{2rs}})-LE^+(\Gamma_{M_{2rs}})=\frac{s}{r-1}\left(\frac{r^3}{2}-2r^2+6r-6\right)-2s\sqrt{2r^2-8r+9}
\end{equation}
and
\begin{equation}\label{m4}
LE^+(\Gamma_{M_{2rs}})-E(\Gamma_{M_{2rs}})=\frac{(r^3-6r^2+8r)s^2-\frac{r^3s}{2}+3r^2s-5rs+4s}{r-1}+2s\sqrt{2r^2-8r+9}-s\sqrt{5r^2-12r+4}.
\end{equation}

Since $\frac{r^3}{2}-2r^2+6r-6 > 0$, $2(r-1)\sqrt{2r^2-8r+9} > 0$ and $(\frac{r^3}{2}-2r^2+6r-6)^2-4(r-1)^2(2r^2-8r+9)= \frac{r^5(r-8)}{4}+2r^4+6r^2(3r-8)+32r \geq 0$ (equality holds if and only if $r = 4$)
we have $\frac{r^3}{2}-2r^2+6r-6 - 2(r-1)\sqrt{2r^2-8r+9} \geq 0$. Therefore, by \eqref{m3}, $(r-1)(LE(\Gamma_{M_{2rs}})-LE^+(\Gamma_{M_{2rs}})) \geq 0$. Hence, $LE(\Gamma_{M_{2rs}}) \geq LE^+(\Gamma_{M_{2rs}})$ equality holds if and only if $G \cong M_{8s}$.  

Again, $2\sqrt{2r^2-8r+9} > 0, \sqrt{5r^2-12r+4} > 0$ and $(2\sqrt{2r^2-8r+9})^2-(\sqrt{5r^2-12r+4})^2$ $=(r-4)(3r-8) \geq 0$ (equality holds if and only if $r = 4$). Therefore, 
$2\sqrt{2r^2-8r+9}-\sqrt{5r^2-12r+4} \geq 0$. Since $r^3-6r^2+8r \geq  \frac{r^3}{2}-3r^2+5r-4$ we have $\frac{(r^3-6r^2+8r)s^2-\frac{r^3s}{2}+3r^2s-5rs+4s}{r-1}+2s\sqrt{2r^2-8r+9}-s\sqrt{5r^2-12r+4} \geq 0$. Therefore, by \eqref{m4}, $LE^+(\Gamma_{M_{2rs}}) \geq  E(\Gamma_{M_{2rs}})$. Hence, 
$E(\Gamma_{M_{2rs}}) \leq LE^+(\Gamma_{M_{2rs}}) \leq LE(\Gamma_{M_{2rs}})$ equality holds if and only if $G \cong M_{8s}$.


For $r \geq 10$, using Theorems \ref{M} and \ref{M_2rs2}, we have 
\begin{equation}\label{m1}
LE(\Gamma_{M_{2rs}})-LE^+(\Gamma_{M_{2rs}})=2s\left(\frac{2r^2-5r+3}{r-1}-\sqrt{2r^2-8r+9}\right)  
\end{equation}
and
\begin{equation}\label{m2}
LE^+(\Gamma_{M_{2rs}})-E(\Gamma_{M_{2rs}})=s\left(\frac{(r^3-6r^2+8r)s-3r^2+11r-8}{r-1}+2\sqrt{2r^2-8r+9}-\sqrt{5r^2-12r+4}\right).
\end{equation}

Since $2r^2-5r+3 > 0$, $(r-1)\sqrt{2r^2-8r+9} > 0$ and $(2r^2-5r+3)^2-(r-1)^2(2r^2-8r+9)=2r(r-2)(r-1)^2>0$
we have $2r^2-5r+3 - (r-1)\sqrt{2r^2-8r+9} > 0$. Therefore, by \eqref{m1}, $(r-1)(LE(\Gamma_{M_{2rs}})-LE^+(\Gamma_{M_{2rs}})) > 0$. Hence, $LE(\Gamma_{M_{2rs}}) > LE^+(\Gamma_{M_{2rs}})$.  

Again, $2\sqrt{2r^2-8r+9} > 0, \sqrt{5r^2-12r+4} > 0$ and $(2\sqrt{2r^2-8r+9})^2-(\sqrt{5r^2-12r+4})^2$ $=(r-4)(3r-8)>0$. Therefore, 
$2\sqrt{2r^2-8r+9}-\sqrt{5r^2-12r+4}>0$. Since $r^3-6r^2+8r > 3r^2-11r+8$ we have $\frac{(r^3-6r^2+8r)s-3r^2+11r-8}{r-1}+2\sqrt{2r^2-8r+9}-\sqrt{5r^2-12r+4}>0$. Therefore, by \eqref{m2}, $LE^+(\Gamma_{M_{2rs}}) >  E(\Gamma_{M_{2rs}})$. Hence, 
$E(\Gamma_{M_{2rs}}) < LE^+(\Gamma_{M_{2rs}}) < LE(\Gamma_{M_{2rs}})$.


\vspace{.5cm}
\noindent (b) \noindent \textbf{Case 1:} $r$ is odd

Here, $|v(\Gamma_{M_{2rs}})|=2rs-s$ and $E(K_{|v(\Gamma_{M_{2rs}})|})=LE(K_{|v(\Gamma_{M_{2rs}})|})= LE^+(K_{|v(\Gamma_{M_{2rs}})|})=4rs-2s-2$. Using Theorem \ref{M}, we have
\begin{equation}\label{E11}
E(\Gamma_{M_{2rs}})-|v(\Gamma_{M_{2rs}})| = s(\sqrt{(r-1)(5r-1)}-r) 
\end{equation}
and
\begin{equation}\label{E12}
E(K_{|v(\Gamma_{M_{2rs}})|})-E(\Gamma_{M_{2rs}}) = 3rs-s-2-s\sqrt{(r-1)(5r-1)}.  
\end{equation}
Since $\sqrt{(r-1)(5r-1)} > 0$, $r > 0$ and $\left(\sqrt{(r-1)(5r-1)}\right)^2-(r)^2=2r(2r-3)+1>0$
we have $ \sqrt{(r-1)(5r-1)}-r> 0$. Therefore, by \eqref{E11}, $E(\Gamma_{M_{2rs}}) > |v(\Gamma_{M_{2rs}})|$.

Again,  $s\sqrt{(r-1)(5r-1)} > 0$, $3rs-s-2 > 0$ and $(3rs-s-2)^2-\left(s\sqrt{(r-1)(5r-1)}\right)^2=4rs(rs-3)+4(s+1)>0$
and so $3rs-s-2 - s\sqrt{(r-1)(5r-1)}> 0$. Therefore, by \eqref{E12}, $E(K_{|v(\Gamma_{M_{2rs}})|}) > E(\Gamma_{M_{2rs}})$. 


\vspace{.5cm}
\noindent \textbf{Case 2:} $r$ is even

Here, $|v(\Gamma_{M_{2rs}})|=2rs-2s$ and $E(K_{|v(\Gamma_{M_{2rs}})|})=LE(K_{|v(\Gamma_{M_{2rs}})|})=LE^+(K_{|v(\Gamma_{M_{2rs}})|})=4rs-4s-2.$ Using Theorem \ref{M}, we have
\begin{equation}\label{E14}
E(\Gamma_{M_{2rs}})-|v(\Gamma_{M_{2rs}})| = s(\sqrt{(r-2)(5r-2)}-r) 
\end{equation}
and
\begin{equation}\label{E15}
E(K_{|v(\Gamma_{M_{2rs}})|})-E(\Gamma_{M_{2rs}}) = 3rs-2s-2-s\sqrt{(r-2)(5r-2)}.
\end{equation}
Since $\sqrt{(r-2)(5r-2)} > 0$, $r > 0$ and $\left(\sqrt{(r-2)(5r-2)}\right)^2-r^2= 4(r(r-3)+1)>0$
we have $ \sqrt{(r-2)(5r-2)}-r> 0$. Therefore, by \eqref{E14}, $E(\Gamma_{M_{2rs}}) > |v(\Gamma_{M_{2rs}})|$.

Again,  $s\sqrt{(r-2)(5r-2)} > 0$, $3rs-2s-2 > 0$ and $(3rs-2s-2)^2-\left(s\sqrt{(r-2)(5r-2)}\right)^2=4rs(rs-3)+4(2s+1)>0$
and so $3rs-2s-2 - s\sqrt{(r-2)(5r-2)}> 0$. Therefore, by \eqref{E15}, $E(K_{|v(\Gamma_{M_{2rs}})|}) > E(\Gamma_{M_{2rs}})$.


\vspace{.5cm}
\noindent (c) \textbf{Case 1:} $r$ is odd  

For $r=3$, using Theorem \ref{M_2rs1}, we have $LE^+(\Gamma_{M_{2rs}})-LE^+(K_{|v(\Gamma_{M_{2rs}})|})=\frac{12s^2-53s}{5}+2+s\sqrt{33}>0$ for all $s \neq 1$. Therefore, for $r=3$ and $s \neq 1$, $LE^+(\Gamma_{M_{2rs}})>LE^+(K_{|v(\Gamma_{M_{2rs}})|})$ which implies $\Gamma_{M_{2rs}}$ is Q-hyperenergetic and consequently part (a) implies $\Gamma_{M_{2rs}}$ is L-hyperenergetic. If $r=3$ and $s=1$, then $G \cong D_6$ so result follows from Theorem \ref{D_{2m}}(c).

For $r \geq 5$, using Theorem \ref{M_2rs1}, we have
$LE^+(\Gamma_{M_{2rs}}) - LE^+(K_{|v(\Gamma_{M_{2rs}})|})= \frac{(2r^3-6r^2+4r)s^2-(12r^2-16r+5)s}{2r-1}+s\sqrt{8r^2-16r+9}+2>0.$ Therefore, $LE^+(\Gamma_{M_{2rs}})>LE^+(K_{|v(\Gamma_{M_{2rs}})|})$ which implies $\Gamma_{M_{2rs}}$ is Q-hyperenergetic and consequently part (a) implies $\Gamma_{M_{2rs}}$ is L-hyperenergetic.

\vspace{.5cm}
\noindent \textbf{Case 2:} $r$ is even

For $r=4$ and $s \neq 1$, using Theorem \ref{M}, we have $LE(K_{|v(\Gamma_{M_{2rs}})|})-LE(\Gamma_{M_{2rs}})=\frac{12s}{7}-2>0$. Therefore, $\Gamma_{M_{8s}}$,  is not L-hyperenergetic and consequently part (a) implies that it is not Q-hyperenergetic. If $r=4$ and $s=1$, then $G \cong D_8$ so result follows from Theorem \ref{D_{2m}}(c).

Using Theorem \ref{M_2rs2}, for $4 < r \leq 8$, we get
\begin{center}
 $LE^+(\Gamma_{M_{2rs}}) - LE^+(K_{|v(\Gamma_{M_{2rs}})|})= \frac{(r^3-6r^2+8r)s^2-(\frac{r^3}{2}-2)s}{r-1}+2s\sqrt{2r^2-8r+9}+2>0$.
\end{center} 
  Therefore,  $LE^+(\Gamma_{M_{2rs}})>LE^+(K_{|v(\Gamma_{M_{2rs}})|})$ which implies $\Gamma_{M_{2rs}}$ is Q-hyperenergetic and consequently part (a) implies $\Gamma_{M_{2rs}}$ is L-hyperenergetic. 

Using Theorem \ref{M_2rs2}, for $r \geq 10$, we get
\begin{center}
$LE^+(\Gamma_{M_{2rs}}) - LE^+(K_{|v(\Gamma_{M_{2rs}})|})= \frac{(r^3-6r^2+8r)s^2-(6r^2-16r+10)s}{r-1}+2s\sqrt{2r^2-8r+9}+2>0.$
\end{center}
	 Therefore, $LE^+(\Gamma_{M_{2rs}})>LE^+(K_{|v(\Gamma_{M_{2rs}})|})$ which implies $\Gamma_{M_{2rs}}$ is Q-hyperenergetic and consequently part (a) implies $\Gamma_{M_{2rs}}$ is L-hyperenergetic.
 \end{proof}

In Theorem \ref{m}, we compare $E(\Gamma_{M_{2rs}})$, $LE(\Gamma_{M_{2rs}})$ and $LE^+(\Gamma_{M_{2rs}})$. However, in the following figures, we show how close are they.


\begin{minipage}[t]{.5\linewidth}
\begin{tikzpicture}
\begin{axis}
[
xlabel={$r$ $\rightarrow$},
ylabel={Energies of $\Gamma_{M_{8r}}$ $\rightarrow$},
xmin=1, xmax=15,
ymin=0, ymax=3200,
grid = both,
minor tick num = 1,
major grid style = {lightgray},
minor grid style = {lightgray!25},
width=.7\textwidth,
height=.7\textwidth,
legend style={legend pos=north west},
 ]
\addplot[domain=5:25,samples at={3,5,7,9,11,13,15,17,19,21,23,25},mark=*,green, samples=11, mark size=.8pt]{4*(x-1)+4*(5*x*x-6*x+1)^(1/2)};
\tiny
\addlegendentry{$E$}
\addplot[domain=5:25,samples at={3,5,7,9,11,13,15,17,19,21,23,25},mark=triangle*,blue,mark size=.8pt, samples=11]{4*(8*x*x*x-24*x*x+16*x+4*x*x-2*x)/(2*x-1)};
\tiny
\addlegendentry{$LE$}
\addplot[domain=5:25,samples at={5,7,9,11,13,15,17,19,21,23,25},mark=square*, red, mark size=.8pt, samples=11]{4*((8*x*x*x-24*x*x+16*x-4*x*x+8*x-3)/(2*x-1)+(8*x*x-16*x+9)^(1/2))}; 
\tiny
\addlegendentry{$LE^+$}
\addplot[domain=3:4,samples at={3},mark=square*, red, mark size=.8pt, samples=18]{36+4*(33)^(1/2)};
\tiny
\end{axis}
\end{tikzpicture}
\vspace{-.2 cm}
\captionsetup{font=footnotesize}

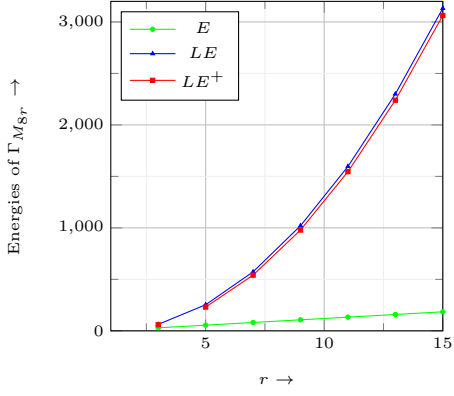
\captionof{figure}{Energies of $\Gamma_{M_{8r}}$}
\end{minipage}
\hspace{0.05cm}
\begin{minipage}[t]{.5\linewidth}
\begin{tikzpicture}
\begin{axis}
[
xlabel={$r$ $\rightarrow$},
ylabel={Energies $\Gamma_{M_{10r}}$ $\rightarrow$},
xmin=2, xmax=16,
ymin=0, ymax=4700,
grid = both,
minor tick num = 1,
major grid style = {lightgray},
minor grid style = {lightgray!25},
width=.7\textwidth,
height=.7\textwidth,
legend style={legend pos=north west},
 ]
\addplot[domain=10:30,samples at={4,6,8,10,12,14,16},mark=*,green, samples=10, mark size=.8pt]{5*(x-2)+5*(5*x*x-12*x+4)^(1/2)};
\tiny
\addlegendentry{$E$}
\addplot[domain=10:30,samples at={4,6,8,10,12,14,16},mark=triangle*,blue,mark size=.8pt, samples=10]{5*(5*x*x*x-30*x*x+40*x+2*x*x-2*x)/(x-1)};
\tiny
\addlegendentry{$LE$}
\addplot[domain=10:30,samples at={10,12,14,16},mark=square*, red, mark size=.8pt, samples=10]{5*(5*x*x*x-30*x*x+40*x-2*x*x+8*x-6)/(x-1)+10*(2*x*x-8*x+9)^(1/2)};
\tiny
\addlegendentry{$LE^+$}
\addplot[domain=4:8,samples at={4,6,8},mark=square*, red, mark size=.8pt, samples=10]{5*(5*x*x*x-30*x*x+40*x-(x*x*x)/(2)+4*x*x-8*x+6)/(x-1)+10*(2*x*x-8*x+9)^(1/2)};
\tiny
\end{axis}
\end{tikzpicture}
\vspace{-.2 cm}
\captionsetup{font=footnotesize}

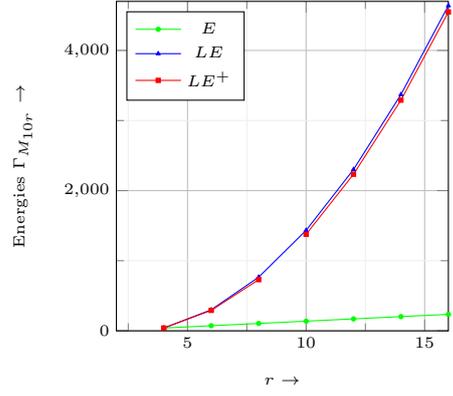
\captionof{figure}{Energies of $\Gamma_{M_{10r}}$}
\end{minipage}

\subsection{The dicyclic groups, $Q_{4n}$}
We consider $Q_{4n}: = \langle x, y : x^{2n} = y^4=1; x^n=y^2; y^{-1}xy = x^{-1} \rangle$, the dicyclic groups of order $4n$ (where $n \geq 2$). Results regarding different energies of non-commuting graphs of $Q_{4n}$ are given below.
\begin{theorem}[\protect{\cite[Corollary 4.1.8 and (4.3.f)]{FWNT21}}]\label{Q4n}
Let $G$ be isomorphic to $Q_{4n}$. Then 
\[E(\Gamma_{Q_{4n}}) = 2\left((n-1)+\sqrt{(n-1)(5n-1)}\right) \, \text{ and } \,
LE(\Gamma_{Q_{4n}}) = \frac{8n(n-1)(n-2)+4n(2n-1)}{2n-1}.\]
\end{theorem}

\begin{theorem}\label{Quarternion}
Let $G$ be isomorphic to $Q_{4n}$. Then 
\[
\Q-spec(\Gamma_{Q_{4n}})=\left\lbrace(4n-4)^n, (2n)^{2n-3}, (4n-6)^{n-1}, (4n-3+\sqrt{8n^2-16n+9})^1, \left(4n-3-\sqrt{8n^2-16n+9}\right)^1\right\rbrace
\] and 
\begin{align*}
LE^+(\Gamma_{Q_{4n}}) = \begin{cases}
 \frac{4n^3-8n^2+6}{2n-1}+2\sqrt{8n^2-16n+9}, & \mbox{if $n \leq 4$}\vspace{.2cm}\\
\frac{8n^3-32n^2+32n-6}{2n-1}+2\sqrt{8n^2-16n+9}, & \mbox{if $n\geq 5$.}
\end{cases}   
\end{align*}
\end{theorem}
\begin{proof}
If $G\cong Q_{4n}$ then $|v(\Gamma_{Q_{4n}})|=4n-2$ and $\Gamma_{Q_{4n}}=K_{n.2,1.(2n-2)}$. Using Theorem \ref{R1}(b), we have 
\begin{align*}
Q_{\Gamma_{Q_{4n}}}(x)=&{\displaystyle \prod_{i=1}^2} (x-(4n-2)+p_i)^{a_i(p_i-1)}{\displaystyle \prod_{i=1}^2}(x-(4n-2)+2p_i)^{a_i}\left(1-{\displaystyle \sum_{i=1}^2}\frac{a_ip_i}{x-(4n-2)+2p_i}\right)\\
  = \,& (x-(4n-4))^{n}(x-2n)^{2n-3}(x-4n+6)^{n}(x-2)\left(1-\frac{2n}{x-4n+6}-\frac{2n-2}{x-2}\right)\\
   =\,& (x-(4n-4))^{n}(x-2n)^{2n-3}(x-(4n-6))^{n-1}(x^2-(8n-6)x+8n^2-8n).
\end{align*}
Thus, 
$$\Q-spec(\Gamma_{Q_{4n}})=\left\lbrace(4n-4)^n, (2n)^{2n-3}, (4n-6)^{n-1}, \left(4n-3+\sqrt{8n^2-16n+9}\right)^1, (4n-3-\sqrt{8n^2-16n+9})^1\right\rbrace.$$

 Number of edges of $\Gamma_{Q_{4n}}^c$ is $2n^2-4n+3$. Thus, $|e(\Gamma_{Q_{4n}})| = \frac{(4n-2)(4n-2-1)}{2}-(2n^2-4n+3)=6n(n-1)$. Now
\[\left|4n-4 - \frac{2|e(\Gamma_{Q_{4n}})|}{|v(\Gamma_{Q_{4n}})|}\right| = \left|\frac{(2n-4)(n-1)}{2n-1}\right|= \frac{(2n-4)(n-1)}{2n-1},\] 

\[\left|2n - \frac{2|e(\Gamma_{Q_{4n}})|}{|v(\Gamma_{Q_{4n}})|}\right| = \left|\frac{-2n(n-2)}{2n-1}\right|= \frac{2n(n-2)}{2n-1},\]

\[\left|4n-6 - \frac{2|e(\Gamma_{Q_{4n}})|}{|v(\Gamma_{Q_{4n}})|}\right| = \left|\frac{2(n^2-5n+3)}{2n-1}\right|= \begin{cases}\frac{-2(n^2-5n+3)}{2n-1}, & \mbox{if $n \leq 4$}\vspace{.2cm}\\
\frac{2(n^2-5n+3)}{2n-1}, & \mbox{if $n \geq 5,$}
 \end{cases}\] 
\begin{align*}
\left|4n-3+\sqrt{8n^2-16n+9} - \frac{2|e(\Gamma_{Q_{4n}})|}{|v(\Gamma_{Q_{4n}})|}\right|=& \left|\sqrt{8n^2-16n+9}+n-\frac{3}{2}+\frac{3}{4n-2}\right|\\
=&\sqrt{8n^2-16n+9}+n-\frac{3}{2}+\frac{3}{4n-2}\end{align*}
 and
\begin{align*} \left|4n-3-\sqrt{8n^2-16n+9} - \frac{2|e(\Gamma_{Q_{4n}})|}{|v(\Gamma_{Q_{4n}})|}\right| = &\left|-\sqrt{8n^2-16n+9}+n-\frac{3}{2}+\frac{3}{4n-2}\right|\\ =& \sqrt{8n^2-16n+9}-n+\frac{3}{2}-\frac{3}{4n-2}.\end{align*}


Therefore, for $n \leq 4$ we have
\begin{align*}
 LE^+(\Gamma_{Q_{4n}}) = &n\times \frac{(2n-4)(n-1)}{2n-1}+(2n-3)\times  \frac{2n(n-2)}{2n-1}+ (n-1) \times \frac{-2(n^2-5n+3)}{2n-1}+\\ &\sqrt{8n^2-16n+9}+n-\frac{3}{2}+\frac{3}{4n-2}+\sqrt{8n^2-16n+9}-n+\frac{3}{2}-\frac{3}{4n-2} 
 \end{align*} 
and for $n \geq 5$ we have
\begin{align*}
LE^+(\Gamma_{Q_{4n}}) = &n\times \frac{(2n-4)(n-1)}{2n-1}+(2n-3)\times  \frac{2n(n-2)}{2n-1}+ (n-1) \times \frac{2(n^2-5n+3)}{2n-1}+\\ &\sqrt{8n^2-16n+9}+n-\frac{3}{2}+\frac{3}{4n-2}+\sqrt{8n^2-16n+9}-n+\frac{3}{2}-\frac{3}{4n-2}   
\end{align*}
Hence, the results follow on simplification.
\end{proof}

\begin{theorem}\label{Q_{4m}}
If $G$ is isomorphic to $Q_{4n}$ then
\begin{enumerate}
\item $E(\Gamma_{Q_{4n}}) \leq  LE^+(\Gamma_{Q_{4n}}) \leq LE(\Gamma_{Q_{4n}})$, equality holds if and only if $G \cong Q_8$.
\item $\Gamma_{Q_{4n}}$ is non-hypoenergetic as well as non-hyperenergetic. 
\item $\Gamma_{Q_8}$ is not L-hyperenergetic and not Q-hyperenergetic. If $n \neq 2$ then $\Gamma_{Q_{4n}}$ is Q-hyperenergetic and L-hyperenergetic. 
\end{enumerate}
\end{theorem}
\begin{proof}
(a) For $n \leq 4$, using Theorems \ref{Q4n} and \ref{Quarternion}, we have 
\begin{equation}\label{E19}
LE(\Gamma_{Q_{4n}})-LE^+(\Gamma_{Q_{4n}})=2\left(\frac{2n^3-4n^2+6n-3}{2n-1}-\sqrt{8n^2-16n+9}\right)    
\end{equation} 
and
\begin{equation}\label{E20}
LE^+(\Gamma_{Q_{4n}})-E(\Gamma_{Q_{4n}})=\frac{2(n-2)(2n^2-2n-1)}{2n-1}+2\sqrt{8n^2-16n+9}-2\sqrt{5n^2-6n+1}.   
\end{equation}
Since $2n^3-4n^2+6n-3 > 0$, $(2n-1)\sqrt{8n^2-16n+9} > 0$ and $(2n^3-4n^2+6n-3)^2-\left(\sqrt{8n^2-16n+9}\right)^2(2n-1)^2=8n^5(n-4)+16n^4+24n^2(3n-4)+32n \geq 0$ (equality holds if and only if $n = 2$)
we have $2n^3-4n^2+6n-3 - (2n-1)\sqrt{8n^2-16n+9} \geq 0$. Therefore, by \eqref{E19}, $(2n-1)(LE(\Gamma_{Q_{4n}})-LE^+(\Gamma_{Q_{4n}})) \geq 0$. Hence, $LE(\Gamma_{Q_{4n}}) \geq LE^+(\Gamma_{Q_{4n}})$ equality holds if and only if $G \cong Q_8$.  

Again, $\sqrt{8n^2-16n+9} > 0, \sqrt{5n^2-6n+1} > 0$ and $\left(\sqrt{8n^2-16n+9}\right)^2-\left(\sqrt{5n^2-6n+1}\right)^2=n(3n-10)+8 \geq 0$ (equality holds if and only if $n = 2$). Therefore, 
$\sqrt{8n^2-16n+9}-\sqrt{5n^2-6n+1} \geq 0$. Since $(n-2)(2n^2-2n-1) \geq 0$ we have $\frac{2(n-2)(2n^2-2n-1)}{2n-1}+2\sqrt{8n^2-16n+9}-2\sqrt{5n^2-6n+1} \geq 0$ (equality holds if and only if $n = 2$). Therefore, by \eqref{E20}, $LE^+(\Gamma_{Q_{4n}}) \geq  E(\Gamma_{Q_{4n}})$. Hence, 
$E(\Gamma_{Q_{4n}}) \leq LE^+(\Gamma_{Q_{4n}}) \leq LE(\Gamma_{Q_{4n}})$ equality holds if and only if $G \cong Q_8$.


For $n \geq 5$, using Theorems \ref{Q4n} and \ref{Quarternion}, we have 
\begin{equation}\label{E17}
LE(\Gamma_{Q_{4n}})-LE^+(\Gamma_{Q_{4n}})=2\left(\frac{8n^2-10n+3}{2n-1}-\sqrt{8n^2-16n+9}\right)     
\end{equation}
and
\begin{equation}\label{E18}
LE^+(\Gamma_{Q_{4n}})-E(\Gamma_{Q_{4n}})=\frac{2(2n^2(4n-9)+19n-4)}{2n-1}+2\sqrt{8n^2-16n+9}-2\sqrt{5n^2-6n+1}.
\end{equation} 
Since $8n^2-10n+3 > 0$, $(2n-1)\sqrt{8n^2-16n+9} > 0$ and $(8n^2-10n+3)^2-\left(\sqrt{8n^2-16n+9}\right)^2(2n-1)^2=32n^3(n-2)+8n(5n-1)>0$
we have $8n^2-10n+3 - (2n-1)\sqrt{8n^2-16n+9} > 0$. Therefore, by \eqref{E17}, $(2n-1)(LE(\Gamma_{Q_{4n}})-LE^+(\Gamma_{Q_{4n}})) > 0$. Hence, $LE(\Gamma_{Q_{4n}}) > LE^+(\Gamma_{Q_{4n}})$.  

Again, $\sqrt{8n^2-16m+9} > 0, \sqrt{5n^2-6n+1} > 0$ and $\left(\sqrt{8n^2-16n+9}\right)^2-\left(\sqrt{5n^2-6n+1}\right)^2=n(3n-10)+8>0$. Therefore, 
$\sqrt{8n^2-16n+9}-\sqrt{5n^2-6n+1}>0$. Since $2n^2(4n-9)+19n-4> 0$ we have $\frac{2(2n^2(4n-9)+19n-4)}{2n-1}+2\sqrt{8n^2-16n+9}-2\sqrt{5n^2-6n+1}>0$. Therefore, by \eqref{E18}, $LE^+(\Gamma_{Q_{4n}}) >  E(\Gamma_{Q_{4n}})$. Hence, 
$E(\Gamma_{Q_{4n}}) < LE^+(\Gamma_{Q_{4n}}) < LE(\Gamma_{Q_{4n}})$.


\vspace{.5cm}

\noindent (b) Here, $|v(\Gamma_{Q_{4n}})|=4n-2=2(2n-1)$ and $E(K_{|v(\Gamma_{Q_{4n}})|})=LE(K_{|v(\Gamma_{Q_{4n}})|})= LE^+(K_{|v(\Gamma_{Q_{4n}})|})=8n-6$. Using Theorem \ref{Q4n},
\begin{equation}\label{E21}
E(\Gamma_{Q_{4n}})-|v(\Gamma_{Q_{4n}})| = 2(\sqrt{(n-1)(5n-1)}-n)      
\end{equation} and
\begin{equation}\label{E22}
E(K_{|v(\Gamma_{Q_{4n}})|})-E(\Gamma_{Q_{4n}})= 2(3(n-1)+1-\sqrt{(n-1)(5n-1)}).
\end{equation}
Since $\sqrt{(n-1)(5n-1)} > 0$, $n > 0$ and $\left(\sqrt{(n-1)(5n-1)}\right)^2-n^2=2n(2n-3)+1>0$
 we have $ \sqrt{(n-1)(5n-1)}-n> 0$. Therefore, by \eqref{E21}, $E(\Gamma_{Q_{4n}}) > |v(\Gamma_{Q_{4n}})|$.

Again,  $\sqrt{(n-1)(5n-1)} > 0$, $3(n-1)+1 > 0$ and $(3(n-1)+1)^2-\left(\sqrt{(n-1)(5n-1)}\right)^2=2n(2n-3)+3>0$
and so $3(n-1)+1 - \sqrt{(n-1)(5n-1)}> 0$. Therefore, by \eqref{E22}, $E(K_{|v(\Gamma_{Q_{4n}})|}) > E(\Gamma_{Q_{4n}})$.


\vspace{.5cm}
\noindent (c) For $n=2$, using Theorem \ref{Q4n}, $LE(\Gamma_{Q_8})=8$ and $LE(K_{|v(\Gamma_{Q_8})}|)=10$. Clearly, $LE(\Gamma_{Q_8}) < LE(K_{|v(\Gamma_{Q_8})}|)$. Therefore, $\Gamma_{Q_8}$ is not L-hyperenergetic and consequently part (a) implies that $\Gamma_{Q_8}$ is not Q-hyperenergetic.

 Using Theorem \ref{Quarternion}, for $2 < n \leq 4$, 
 $$
 LE^+(\Gamma_{Q_{4n}})-LE^+(K_{|v(\Gamma_{Q_{4n}})|})= \frac{4n(n-1)(n-5)}{2n-1}+2\sqrt{8n^2-16n+9}>0.
 $$
  Also, for $n \geq 5$, $LE^+(\Gamma_{Q_{4n}})-LE^+(K_{|v(\Gamma_{Q_{4n}})|})= \frac{8n^2(n-6)+52n-12}{2n-1}+2\sqrt{8n^2-16n+9}>0$. Therefore, $LE^+(\Gamma_{Q_{4n}})>LE^+(K_{|v(\Gamma_{Q_{4n}})|})$ which implies $\Gamma_{Q_{4n}}$ is Q-hyperenergetic and consequently part (a) implies that $\Gamma_{Q_{4n}}$ is L-hyperenergetic. Hence, the result holds. 
\end{proof}

\subsection{The groups $U_{6n}$}
We consider $U_{6n}: = \langle x, y : x^{2n} = y^3=1; x^{-1}yx = y^{-1} \rangle$, the groups of order $6n$. Results regarding different energies of non-commuting graphs of $U_{6n}$ are given below.

\begin{theorem}[\protect{\cite[Corollary 4.1.9 and putting $m=3$ in (4.3.c)]{FWNT21}}]\label{U}
Let $G$ be isomorphic to $U_{6n}$. Then $E(\Gamma_{U_{6n}})) = 2n(1+\sqrt{7})\text{ and }$ $LE(\Gamma_{U_{6n}})) = \frac{12n^2+30n}{5}$.
\end{theorem}

\begin{theorem}\label{U6n}
Let $G$ be isomorphic to $U_{6n}$. Then
\[\Q-spec(\Gamma_{U_{6n}})=\left\{(3n)^{2n+1}, (4n)^{3n-3},\left(\frac{(9+\sqrt{33})n}{2}\right)^1, \left(\frac{(9-\sqrt{33})n}{2}\right)^1\right\}\] and \[LE^+(\Gamma_{U_{6n}}) = \frac{12n^2-3n}{5}+\sqrt{33}n.\]
\end{theorem}
\begin{proof}
If $G\cong U_{6n}$ then $|v(\Gamma_{U_{6n}})|=5n$ and $\Gamma_{U_{6n}}=K_{1.2n,3.n}$. Using Theorem \ref{R1}(b), we have 
\begin{align*}
Q_{\Gamma_{U_{6n}}}(x)=&{\displaystyle \prod_{i=1}^2} (x-5n+p_i)^{a_i(p_i-1)}{\displaystyle \prod_{i=1}^2}(x-5n+2p_i)^{a_i}\left(1-{\displaystyle \sum_{i=1}^2}\frac{a_ip_i}{x-5n+2p_i}\right)\\
= \,& (x-3n)^{2n-1}(x-4n)^{3n-3}(x-n)(x-3n)^{3}\left(1-\frac{2n}{x-n}-\frac{3n}{x-3n}\right)\\
   = \,& (x-3n)^{2n+1}(x-4n)^{3n-3}(x^2-9nx+12n^2).
\end{align*}
Thus, $\Q-spec(\Gamma_{U_{6n}})=\left\{(3n)^{2n+1}, (4n)^{3n-3},\left(\frac{(9+\sqrt{33})n}{2}\right)^1, \left(\frac{(9-\sqrt{33})n}{2}\right)^1\right\}$.

 Number of edges of $\Gamma_{U_{6n}}^c$ is $\frac{7n^2-5n}{2}$. Thus, $|e(\Gamma_{U_{6n}})| = \frac{5n(5n-1)}{2}-\frac{7n^2-5n}{2}=\frac{18n^2}{2}$. Now
\[\left|3n- \frac{2|e(\Gamma_{U_{6n}})|}{|v(\Gamma_{U_{6n}})|}\right| = \left|\frac{-3n}{5}\right| =\frac{3n}{5}, \quad
\left|4n- \frac{2|e(\Gamma_{U_{6n}})|}{|v(\Gamma_{U_{6n}})|}\right| = \left|\frac{2n}{5}\right| =\frac{2n}{5},\] \[\left|\frac{(9+\sqrt{33})n}{5}- \frac{2|e(\Gamma_{U_{6n}})|}{|v(\Gamma_{U_{6n}})|}\right| = \left|\frac{(9+5\sqrt{33})n}{10}\right| =\frac{(9+5\sqrt{33})n}{10}\] and \[\left|\frac{(9-\sqrt{33})n}{5}- \frac{2|e(\Gamma_{U_{6n}})|}{|v(\Gamma_{U_{6n}})|}\right| = \left|\frac{(9-5\sqrt{33})n}{10}\right| =\frac{(5\sqrt{33}-9)n}{10}.
\]
Therefore, 
 $LE^+(\Gamma_{U_{6n}}) = (2n+1)\times \frac{3n}{5}+(3n-3)\times \frac{2n}{5}+ \frac{(9+5\sqrt{33})n}{10}+\frac{(5\sqrt{33}-9)n}{10}$ and the result follows on simplification.
\end{proof}

\begin{theorem}\label{U1}
If $G$ is isomorphic to $U_{6n}$ then
\begin{enumerate}
\item $E(\Gamma_{U_{6n}})<LE^+(\Gamma_{U_{6n}})<LE(\Gamma_{U_{6n}})$.
\item $\Gamma_{U_{6n}}$ is non-hypoenergetic as well as non-hyperenergetic. 
\item $\Gamma_{U_{6n}}$ is Q-hyperenergetic and L-hyperenergetic.
\end{enumerate}
\end{theorem}

\begin{proof}
(a) Using Theorems \ref{U} and \ref{U6n}, we have \begin{equation}\label{E23}
LE(\Gamma_{U_{6n}})-LE^+(\Gamma_{U_{6n}})= \frac{33n}{5}-\sqrt{33}n>0
\end{equation}
and
\begin{equation}\label{E24}
LE^+(\Gamma_{U_{6n}})-E(\Gamma_{U_{6n}})= \frac{12n^2-13n}{5}+(\sqrt{33}-2\sqrt{7})n>0.
\end{equation}
Thus, the conclusion is drawn from equations \eqref{E23} and \eqref{E24}.

\vspace{.5cm}
\noindent (b) Here, $|v(\Gamma_{U_{6n}}))|=5n$ and $E(K_{|v(\Gamma_{U_{6n}}))|}) = LE(K_{|v(\Gamma_{U_{6n}}))|}) = LE^+(K_{|v(\Gamma_{U_{6n}}))|})= 10n-2$. Using Theorem \ref{U}, we have
\begin{equation}\label{E25}
E(\Gamma_{U_{6n}})-|v(\Gamma_{U_{6n}})|= (2\sqrt{7}-3)n>0
\end{equation} 
and
\begin{equation}\label{E26}
E(K_{|v(\Gamma_{U_{6n}})|})-E(\Gamma_{U_{6n}})= (8-2\sqrt{7})n-2>0.
\end{equation} 
Thus, the conclusion is drawn from equations \eqref{E25} and \eqref{E26}.

\vspace{.5cm}
\noindent (c) Using Theorem \ref{U6n}, we have $LE^+(\Gamma_{U_{6n}})-LE^+(K_{|v(\Gamma_{U_{6n}})|})= \frac{12n^2-53n+10}{5}+n\sqrt{33}>0.$ Therefore, $LE^+(\Gamma_{U_{6n}})>LE^+(K_{|v(\Gamma_{U_{6n}})|})$ which implies $\Gamma_{U_{6n}}$ is Q-hyperenergetic and consequently part (a) implies $\Gamma_{U_{6n}}$ is L-hyperenergetic. 
\end{proof}

 In Theorems \ref{Q_{4m}} and \ref{U1}, we compare $E(\Gamma_G)$, $LE(\Gamma_G)$ and $LE^+(\Gamma_G)$ if   $G \cong Q_{4n}$ and $U_{6n}$ respectively. However, in the following figures, we show how close are they for both the groups.
\vspace{.3cm} 

\begin{minipage}[t]{.5\linewidth}
\begin{tikzpicture}
\begin{axis}
[
xlabel={$n$ $\rightarrow$},
ylabel={Energies of $\Gamma_{Q_{4n}}$ $\rightarrow$},
xmin=1, xmax=20,
ymin=0, ymax=1500,
grid = both,
minor tick num = 1,
major grid style = {lightgray},
minor grid style = {lightgray!25},
width=.7\textwidth,
height=.7\textwidth,
legend style={legend pos=north west},
 ]
\addplot[domain=2:20,samples at={2,3,4,5,...,18,19,20},mark=*,green, samples=39, mark size=.8pt]{2*(x-1)+2*(5*x*x-6*x+1)^(1/2)};
\tiny
\addlegendentry{$E$}
\addplot[domain=2:20,samples at={2,3,4,5,...,18,19,20},mark=triangle*,blue,mark size=.8pt, samples=39]{(8*x*(x-1)*(x-2)+4*x*(2*x-1))/(2*x-1)};
\tiny
\addlegendentry{$LE$}
\addplot[domain=5:20,samples at={5,6,7,...,18,19,20},mark=square*, red, mark size=.8pt, samples=39]{(8*x*x*x-32*x*x+32*x-6)/(2*x-1)+2*(8*x*x-16*x+9)^(1/2)};
\tiny
\addlegendentry{$LE^+$}
\addplot[domain=2:4,samples at={2,3,4},mark=square*, red, mark size=.8pt, samples=39]{(4*x*x*x-8*x*x+6)/(2*x-1)+2*(8*x*x-16*x+9)^(1/2)};
\tiny
\end{axis}
\end{tikzpicture}
\vspace{-.2 cm}
\captionsetup{font=footnotesize}

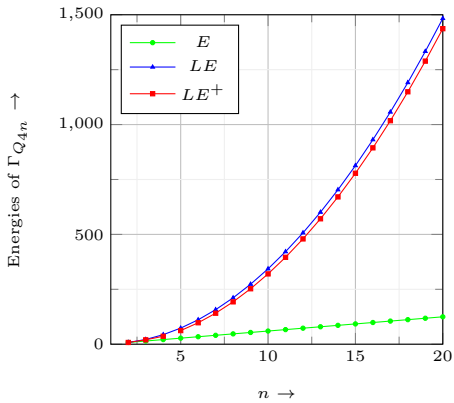
\captionof{figure}{Energies of $\Gamma_{Q_{4n}}$}
\end{minipage}
\hspace{0.05cm}
\begin{minipage}[t]{.5\linewidth}
\begin{tikzpicture}
\begin{axis}
[
xlabel={$n$ $\rightarrow$},
ylabel={Energies of $\Gamma_{U_{6n}}$ $\rightarrow$},
xmin=0, xmax=12,
ymin=0, ymax=450,
grid = both,
minor tick num = 1,
major grid style = {lightgray},
minor grid style = {lightgray!25},
width=.7\textwidth,
height=.7\textwidth,
legend style={legend pos=north west},
 ]
\addplot[domain=2:40,samples at={2,3,4,5,...,38,39,40},mark=*,green, samples=39, mark size=.8pt]{2*x*(1+(7)^(1/2))};
\tiny
\addlegendentry{$E$}
\addplot[domain=2:40,samples at={2,3,4,5,...,38,39,40},mark=triangle*,blue,mark size=.8pt, samples=39]{(12*x*x+30*x)/5};
\tiny
\addlegendentry{$LE$}
\addplot[domain=2:40,samples at={2,3,4,5,...,38,39,40},mark=square*, red, mark size=.8pt, samples=39]{(12*x*x-3*x)/5+(33)^(1/2)*x};
\tiny
\addlegendentry{$LE^+$}
\end{axis}
\end{tikzpicture}
\vspace{-.2 cm}
\captionsetup{font=footnotesize}

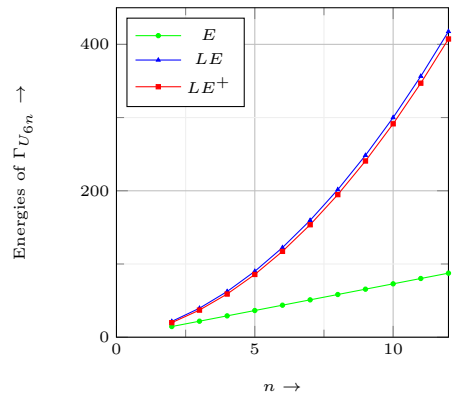
\captionof{figure}{Energies of $\Gamma_{U_{6n}}$}
\end{minipage}
\vspace{0.5cm}

It can be seen that if $G$ is isomorphic to $D_{2m}$, $QD_{2^n}$, $M_{2rs}$, $Q_{4n}$ or $U_{6n}$, then the central quotient of $G$ is also isomorphic to some dihedral group. Therefore, we conclude this section with the following theorems for the non-commuting graphs of the groups $G$ such that $\frac{G}{Z(G)} \cong D_{2m}$.

\begin{theorem}[\protect{\cite[Theorem 4.1.5]{FWNT21}}]\label{D}
Let $\frac{G}{Z(G)} $ be isomorphic to $ D_{2m}$ ($m \geq 3$) and $|Z(G)|=n$. Then
\[E(\Gamma_G) = n\left((m-1)+\sqrt{(m-1)(5m-1)}\right) \, \text{ and } \,
LE(\Gamma_G)=\frac{n}{2m-1}\left((2m^3-6m^2+4m)n+4m^2-2m\right).\]
\end{theorem}

\begin{theorem}\label{D2m}
Let $\frac{G}{Z(G)} $ be isomorphic to $D_{2m}$  ($m \geq 3$) and $|Z(G)|=n$. Then
\begin{align*}
 \Q-spec(\Gamma_G)=&\left\{((2m-2)n)^{m(n-1)}, (mn)^{(m-1)n-1}, ((2m-3)n)^{m-1}, \left(\frac{n(4m-3+\sqrt{8m^2-16m+9})}{2}\right)^1,\right.\\ & \left.\left(\frac{n(4m-3-\sqrt{8m^2-16m+9})}{2}\right)^1\right\}   
\end{align*}
and
\begin{align*}
LE^+(\Gamma_G) = \begin{cases}
\frac{12n^2-3n}{5}+n\sqrt{33}, & \mbox{if $m=3$}\vspace{.2cm}\\
\frac{48n^2-29n}{7}+n\sqrt{73}, & \mbox{if $m=4$}\vspace{.2cm}\\
\frac{(2m^3-6m^2+4m)n^2-(4m^2-8m+3)n}{2m-1}+n\sqrt{8m^2-16m+9}, & \mbox{if $m\geq 5$.}
\end{cases}   
\end{align*}
\end{theorem}
\begin{proof}
If $\frac{G}{Z(G)}\cong D_{2m}$ then $|v(\Gamma_G)|=(2m-1)n$ and $\Gamma_G=K_{m.n,1.((m-1)n)}$. Using Theorem \ref{R1}(b), we have 
\begin{align*}
Q_{\Gamma_G}(x)=&{\displaystyle \prod_{i=1}^2} (x-(2m-1)n+p_i)^{a_i(p_i-1)}{\displaystyle \prod_{i=1}^2}(x-(2m-1)n+2p_i)^{a_i} \left(1-{\displaystyle \sum_{i=1}^2}\frac{a_ip_i}{x-(2m-1)n+2p_i}\right)\\
  = \,& (x-(2m-2)n)^{m(n-1)}(x-mn)^{(m-1)n-1}(x-(2m-3)n)^{m}(x-n)\left(1-\frac{mn}{x-(2m-3)n}-\frac{(m-1)n}{x-n}\right)\\
   = \,& (x-(2m-2)n)^{m(n-1)}(x-mn)^{(m-1)n-1}(x-(2m-3)n)^{m-1}(x^2-(4m-3)nx+(2m^2-2m)n^2).
\end{align*}
Thus,\begin{align*}
 \Q-spec(\Gamma_G)=&\left\{((2m-2)n)^{m(n-1)}, (mn)^{(m-1)n-1}, ((2m-3)n)^{m-1}, \left(\frac{n(4m-3+\sqrt{8m^2-16m+9})}{2}\right)^1,\right.\\ & \left.\left(\frac{n(4m-3-\sqrt{8m^2-16m+9})}{2}\right)^1\right\}. 
\end{align*}

 Number of edges of $\Gamma_G^c$ is $\frac{(m^2-m+1)n^2-(2m-1)n}{2}$. Thus, $|e(\Gamma_G)| = \frac{(2m-1)^2n^2-(2m-1)n}{2}-\frac{(m^2-m+1)n^2-(2m-1)n}{2}$ $=\frac{3m(m-1)n^2}{2}$. Now
\[\left|(2m-2)n- \frac{2|e(\Gamma_G)|}{|v(\Gamma_G)|}\right| = \left|\frac{(m-1)(m-2)n}{2m-1}\right| =\frac{(m-1)(m-2)n}{2m-1},\] 
\[\left|mn - \frac{2|e(\Gamma_G)|}{|v(\Gamma_G)|}\right| = \left|\frac{-m(m-2)n}{2m-1}\right| = \frac{m(m-2)n}{2m-1},\]
\begin{align*}
\left|(2m-3)n - \frac{2|e(\Gamma_G)|}{|v(\Gamma_G)|}\right| = \left|\frac{(m^2-5m+3)n}{2m-1}\right|= \begin{cases}
\frac{(-m^2+5m-3)n}{2m-1}, & \mbox{if $m \leq 4$}\vspace{.2cm}\\
\frac{(m^2-5m+3)n}{2m-1}, & \mbox{if $m \geq 5,$}\end{cases}   
\end{align*} 
\begin{align*}
\left|\frac{n}{2}\left(4m-3+\sqrt{8m^2-16m+9}\right) - \frac{2|e(\Gamma_G)|}{|v(\Gamma_G)|}\right| =& \left|\frac{n}{2}\left(\sqrt{8m^2-16m+9}+m-\frac{3}{2}+\frac{3}{4m-2}\right)\right| \\=&\frac{n}{2}\left(\sqrt{8m^2-16m+9}+m-\frac{3}{2}+\frac{3}{4m-2}\right)    
\end{align*}
and
\begin{align*}
\left|\frac{n}{2}\left(4m-3-\sqrt{8m^2-16m+9}\right) - \frac{2|e(\Gamma_G)|}{|v(\Gamma_G)|}\right| =& \left|\frac{n}{2}\left(-\sqrt{8m^2-16m+9}+m-\frac{3}{2}+\frac{3}{4m-2}\right)\right|\\ =&\frac{n}{2}\left(\sqrt{8m^2-16m+9}-m+\frac{3}{2}-\frac{3}{4m-2}\right).   
\end{align*} 
Therefore, for $m \leq 4$ we have
\begin{align*}
	LE^+(\Gamma_G) =& m(n-1)\times \frac{(m-1)(m-2)n}{2m-1}+((m-1)n-1)\times \frac{m(m-2)n}{2m-1}+ (m-1)\times \frac{(-m^2+5m-3)n}{2m-1}+\\ &\frac{n}{2}\left(\sqrt{8m^2-16m+9}+m-\frac{3}{2}+\frac{3}{4m-2}\right)+\frac{n}{2}\left(\sqrt{8m^2-16m+9}-m+\frac{3}{2}-\frac{3}{4m-2}\right)
\end{align*}
and for $m \geq 5$ we have
\begin{align*}
LE^+(\Gamma_G) =& m(n-1)\times \frac{(m-1)(m-2)n}{2m-1}+((m-1)n-1)\times \frac{m(m-2)n}{2m-1}+ (m-1)\times \frac{(m^2-5m+3)n}{2m-1}+\\ &\frac{n}{2}\left(\sqrt{8m^2-16m+9}+m-\frac{3}{2}+\frac{3}{4m-2}\right)+\frac{n}{2}\left(\sqrt{8m^2-16m+9}-m+\frac{3}{2}-\frac{3}{4m-2}\right)   
\end{align*}
Hence,  the results follow on simplification.
\end{proof}

\begin{theorem}\label{G/ZGD2m}
If  $\frac{G}{Z(G)}$ is isomorphic to $D_{2m}$  ($m \geq 3$) and $|Z(G)|=n$ then
\begin{enumerate}
\item $E(\Gamma_G) < LE^+(\Gamma_G) < LE(\Gamma_G)$.
\item $\Gamma_G$ is non-hypoenergetic as well as non-hyperenergetic.
\item $\Gamma_G$ is L-hyperenergetic but not Q-hyperenergetic if $m = 3$ and $|Z(G)|=1$. 
For $m = 3, 4$ and $|Z(G)|\ne 1$ or $m \geq 5$ and $|Z(G)| \geq 1$,   $\Gamma_G$ is Q-hyperenergetic and L-hyperenergetic. 
\end{enumerate}
\end{theorem}

\begin{proof}
(a) For $m=3$, using Theorems \ref{D} and \ref{D2m}, we have $LE(\Gamma_G)-LE^+(\Gamma_G)=\frac{33n}{5}-n\sqrt{33}>0$ and $LE^+(\Gamma_G)-E(\Gamma_G)= \frac{12n^2-13n}{5} + (\sqrt{33}-2\sqrt{7})n>0$.

For $m=4$, using Theorems \ref{D} and \ref{D2m}, we have $LE(\Gamma_G)-LE^+(\Gamma_G)=\frac{85n}{7}-n\sqrt{73}>0$ and $LE^+(\Gamma_G)-E(\Gamma_G)= \frac{48n^2-50n}{5} + (\sqrt{73}-\sqrt{57})n>0$.

For $m \geq 5$, using Theorems \ref{D} and \ref{D2m}, we have
\begin{equation}\label{Dm11}
LE(\Gamma_G)-LE^+(\Gamma_G)=\frac{n(8m^2-10m+3)}{2m-1}-n\sqrt{8m^2-16m+9}
\end{equation}
and
\begin{equation}\label{Dm22}
LE^+(\Gamma_G)-E(\Gamma_G)=n\left(\frac{(2m^3-6m^2+4m)n-6m^2+11m-4}{2m-1}+\sqrt{8m^2-16m+9}-\sqrt{5m^2-6m+1}\right).
\end{equation}

Since $8m^2-10m+3 > 0$, $(2m-1)\sqrt{8m^2-16m+9} > 0$ and $(8m^2-10m+3)^2-(2m-1)^2(8m^2-16m+9)=32m^4-64m^3+40m^2-8m>0$
we have $8m^2-10m+3 - (2m-1)\sqrt{8m^2-16m+9} > 0$. Therefore, by \eqref{Dm11}, $(2m-1)(LE(\Gamma_G)-LE^+(\Gamma_G)) > 0$. Hence, $LE(\Gamma_G) > LE^+(\Gamma_G)$.  

Again, $\sqrt{8m^2-16m+9} > 0, \sqrt{5m^2-6m+1} > 0$ and $(\sqrt{8m^2-16m+9})^2-(\sqrt{5m^2-6m+1})^2$ $=m(3m-10)+8>0$. Therefore, 
$\sqrt{8m^2-16m+9}-\sqrt{5m^2-6m+1}>0$. Since $2m^3-6m^2+4m > 6m^2-11m+4$ we have $\frac{(2m^3-6m^2+4m)n-6m^2+11m-4}{2m-1}+\sqrt{8m^2-16m+9}-\sqrt{5m^2-6m+1}>0$. Therefore, by \eqref{Dm22}, $LE^+(\Gamma_G) >  E(\Gamma_G)$. Hence, 
$E(\Gamma_G) < LE^+(\Gamma_G) < LE(\Gamma_G)$.



\vspace{.5cm}


\noindent (b) Here, $|v(\Gamma_G)|=2mn-n$ and $E(K_{|v(\Gamma_G)|})=LE(K_{|v(\Gamma_G)|})= LE^+(K_{|v(\Gamma_G)|})= 4mn-2n-2.$ Using Theorem \ref{D}, we have
\begin{equation}\label{E27}
E(\Gamma_G)-|v(\Gamma_G)| = n(\sqrt{(m-1)(5m-1)}-m)  
\end{equation}
and
\begin{equation}\label{E28}
E(K_{|v(\Gamma_G)|})-E(\Gamma_G) = 3mn-n-2-n\sqrt{(m-1)(5m-1)}.  
\end{equation}
Since $\sqrt{(m-1)(5m-1)}>0, m>0$ and $(\sqrt{(m-1)(5m-1)})^2-m^2=2m(2m-3)+1>0$ we have $ \sqrt{(m-1)(5m-1)}-m> 0$. Therefore, by \eqref{E27}, $E(\Gamma_G) > |v(\Gamma_G)|$.

Again,  $n\sqrt{(m-1)(5m-1)} > 0$, $3mn-n-2 > 0$ and $(3mn-n-2)^2-\left(n\sqrt{(m-1)(5m-1)}\right)^2=4mn(mn-3)+4(n+1)>0$
and so $3mn-n-2 - n\sqrt{(m-1)(5m-1)}> 0$. Therefore, by \eqref{E28}, $E(K_{|v(\Gamma_G)|}) > E(\Gamma_G)$. 


\vspace{.5cm}
\noindent (c) Using Theorem \ref{D2m}, for $m=3$, we have $LE^+(\Gamma_G)-LE^+(K_{|v(\Gamma_G)|})=\frac{12n^2-53n}{7}+2+n\sqrt{33}>0$ for all $n \neq 1$. Therefore, for $m=3$ and $n \neq 1$, $LE^+(\Gamma_G))>LE^+(K_{|v(\Gamma_G))|})$ which implies $\Gamma_G$ is Q-hyperenergetic and consequently part (a) implies $\Gamma_G$ is L-hyperenergetic. If $m=3$ and $n=1$, then $G \cong D_6$ so result follows from Theorem \ref{D_{2m}}(c). Using Theorem \ref{D2m}, for $m=4$, we have $LE^+(\Gamma_G)-LE^+(K_{|v(\Gamma_G)|})=\frac{48n^2-127n}{7}+2+n\sqrt{73}>0$ for all $n \neq 1$. Therefore, for $m=4$ and $n \neq 1$, $LE^+(\Gamma_G))>LE^+(K_{|v(\Gamma_G))|})$ which implies $\Gamma_G$ is Q-hyperenergetic and consequently part (a) implies $\Gamma_G$ is L-hyperenergetic. The case $m=4$ and $n=1$ does not arise since $|Z(D_8)| = 2$.
 
 For $r \geq 5$, using Theorem \ref{D2m}, we have
$LE^+(\Gamma_G) - LE^+(K_{|v(\Gamma_G)|})= \frac{(2m^3-6m^2+4m)n^2-(12m^2-16m+5)n}{2m-1}+n\sqrt{8m^2-16m+9}+2>0.$ Therefore, $LE^+(\Gamma_G)>LE^+(K_{|v(\Gamma_G)|})$ which implies $\Gamma_G$ is Q-hyperenergetic and consequently part (a) implies $\Gamma_G$ is L-hyperenergetic.
 
\end{proof}

In Theorem \ref{G/ZGD2m}, we compare $E(\Gamma_G)$, $LE(\Gamma_G)$ and $LE^+(\Gamma_G)$. However, in the following figures, we show how close are they.

\begin{minipage}[t]{.5\linewidth}
\begin{tikzpicture}
\begin{axis}
[
xlabel={$m$ $\rightarrow$},
ylabel={Energies of $\Gamma_G$ $\rightarrow$},
xmin=1, xmax=15,
ymin=0, ymax=1800,
grid = both,
minor tick num = 1,
major grid style = {lightgray},
minor grid style = {lightgray!25},
width=.7\textwidth,
height=.7\textwidth,
legend style={legend pos=north west},
 ]
\addplot[domain=3:15,samples at={3,4,5,...,12,13,14,15},mark=*,green, samples=8, mark size=.8pt]{3*(x-1)+3*(5*x*x-6*x+1)^(1/2)};
\tiny
\addlegendentry{$E$}
\addplot[domain=3:15,samples at={3,4,5,...,12,13,14,15},mark=triangle*,blue,mark size=.8pt, samples=8]{3*(6*x*x*x-18*x*x+12*x+4*x*x-2*x)/(2*x-1)};
\tiny
\addlegendentry{$LE$}
\addplot[domain=5:15,samples at={5,...,12,13,14,15},mark=square*, red, mark size=.8pt, samples=8]{3*(6*x*x*x-18*x*x+12*x-4*x*x+8*x-3)/(2*x-1)+3*(8*x*x-16*x+9)^(1/2)};
\tiny
\addlegendentry{$LE^+$}
\addplot[domain=3:4,samples at={3,4},mark=square*, red, mark size=.8pt, samples=8]{3*(6*x*x*x-18*x*x+12*x-4*x*x+8*x-3)/(2*x-1)+3*(8*x*x-16*x+9)^(1/2)};
\tiny
\end{axis}
\end{tikzpicture}
\vspace{-.2 cm}
\captionsetup{font=footnotesize}

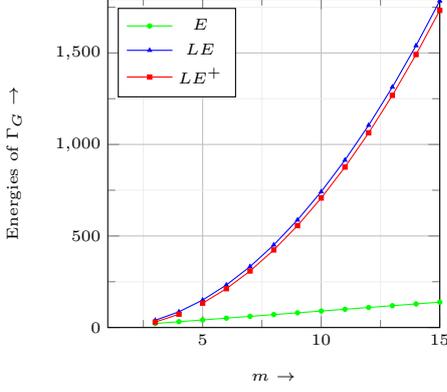
\captionof{figure}{Energies of $\Gamma_G$,  $\frac{G}{Z(G)} \cong D_{2m}$, $|Z(G)|=3$}
\end{minipage}
\hspace{0.05cm}
\begin{minipage}[t]{.5\linewidth}
\begin{tikzpicture}
\begin{axis}
[
xlabel={$m$ $\rightarrow$},
ylabel={Energies of $\Gamma_G$ $\rightarrow$},
xmin=1, xmax=15,
ymin=0, ymax=3200,
grid = both,
minor tick num = 1,
major grid style = {lightgray},
minor grid style = {lightgray!25},
width=.7\textwidth,
height=.7\textwidth,
legend style={legend pos=north west},
 ]
\addplot[domain=5:15,samples at={3,4,5,6,7,8,9,10,11,12,13,14,15},mark=*,green, samples=8, mark size=.8pt]{4*(x-1)+4*(5*x*x-6*x+1)^(1/2)};
\tiny
\addlegendentry{$E$}
\addplot[domain=5:15,samples at={3,4,5,6,7,8,9,10,11,12,13,14,15},mark=triangle*,blue,mark size=.8pt, samples=8]{4*(8*x*x*x-24*x*x+16*x+4*x*x-2*x)/(2*x-1)};
\tiny
\addlegendentry{$LE$}
\addplot[domain=5:15,samples at={5,6,7,8,9,10,11,12,13,14,15},mark=square*, red, mark size=.8pt, samples=8]{4*(8*x*x*x-24*x*x+16*x-4*x*x+8*x-3)/(2*x-1)+4*(8*x*x-16*x+9)^(1/2)};
\tiny
\addlegendentry{$LE^+$}
\addplot[domain=3:4,samples at={3,4},mark=square*, red, mark size=.8pt, samples=8]{4*(8*x*x*x-24*x*x+16*x-2*x*x*x+8*x*x-8*x+3)/(2*x-1)+4*(8*x*x-16*x+9)^(1/2)};
\tiny
\end{axis}
\end{tikzpicture}
\vspace{-.2 cm}
\captionsetup{font=footnotesize}

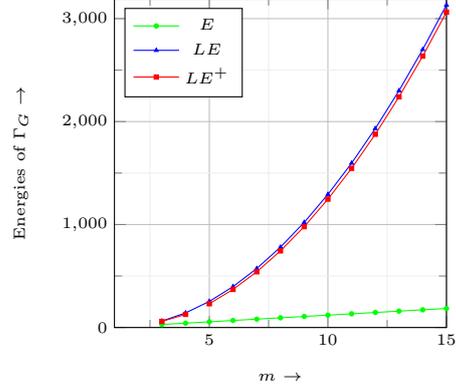
\captionof{figure}{Energies of $\Gamma_G$,  $\frac{G}{Z(G)} \cong D_{2m}$, $|Z(G)|=4$}
\end{minipage}

\section{$\frac{G}{Z(G)}$ is isomorphic to $\mathbb{Z}_p \times \mathbb{Z}_p$}
We compute the Signless Laplacian spectrum and Signless Laplacian energy of $\Gamma_G$ considering the group $G$ whose central quotient is isomorphic to $\mathbb{Z}_p \times \mathbb{Z}_p$, where $p$ is a prime. Further we compare energy, Laplacian energy and Signless Laplacian energy of $\Gamma_G$ and look into the hyper- and hypo-properties of of $\Gamma_G$.

\begin{theorem}[\protect{\cite[Theorem 2.2]{DN18} and \cite[Theorem 4.1.1 (c)]{FWNT21}}]\label{Z_p}
Let $\frac{G}{Z(G)}$ be isomorphic to $ \mathbb{Z}_p \times \mathbb{Z}_p$. Then
$E(\Gamma_G) = LE(\Gamma_G)= 2p(p-1)|Z(G)|$.
In particular, if $G$ is non-abelian and $|G|=p^3$ then $E(\Gamma_G) = LE(\Gamma_G) = 2p^2(p-1)$.
\end{theorem}

\begin{theorem}\label{Z_p*Z_p}
Let $\frac{G}{Z(G)}$ be isomorphic to $ \mathbb{Z}_p \times \mathbb{Z}_p$. Then 
\[\Q-spec(\Gamma_G)=\left\lbrace(pn(p-1))^{(p^2-1)n-(p+1)},\left(n(p-1)^2\right)^p,(2pn(p-1))^1\right\rbrace, \text{ where } |Z(G)| = n\]  and $LE^+(\Gamma_G)= 2p(p-1)|Z(G)|$. In particular, if $G$ is non-abelian and $|G|=p^3$ then $LE^+(\Gamma_G) = 2p^2(p-1)$.    
\end{theorem}

\begin{proof}
If $\frac{G}{Z(G)} \cong \mathbb{Z}_p \times \mathbb{Z}_p$  then $|v(\Gamma_G)|=(p^2-1)n$ and $\Gamma_G=K_{(p+1).(p-1)n}$, where $|Z(G)| = n$. Using Theorem \ref{R1}(b), we have 
\begin{align*}
Q_{\Gamma_G}&(x)\\
=& (x-(p^2-1)n+(p-1)n)^{(p+1)((p-1)n-1)}(x-(p^2-1)n+2(p-1)n)^{p+1}
 \left(1-\frac{(p^2-1)n}{x-(p^2-1)n+2(p-1)n}\right)\\
   =& (x-pn(p-1))^{(p^2-1)n-(p+1)}(x-n(p-1)^2)^{p}(x-2pn(p-1)).
\end{align*}
Thus, $\Q-spec(\Gamma_G)=\left\lbrace(pn(p-1))^{(p^2-1)n-(p+1)}, \left(n(p-1)^2\right)^p,(2pn(p-1))^1\right\rbrace$.

 Number of edges of $\Gamma_G^c$ is $\frac{n(p^2-1)(pn-n-1)}{2}$. Therefore, 
\begin{center}
 $|e(\Gamma_G)| = \frac{n^2(p^2-1)^2-n(p^2-1)}{2}-\frac{n(p^2-1)(pn-n-1)}{2} =$ $\frac{(p^2-p)(p^2-1)n^2}{2}$. 
\end{center}
Now
\[
\left|pn(p-1) - \frac{2|e(\Gamma_G)|}{|v(\Gamma_G)|}\right| = \left|pn(p-1)-(p^2-p)n\right|= 0, \quad \left|n(p-1)^2 - \frac{2|e(\Gamma_G)|}{|v(\Gamma_G)|}\right| = \left|n-pn\right|  =pn-n\] and
\[\left|2pn(p-1)- \frac{2|e(\Gamma_G)|}{|v(\Gamma_G)|}\right| = \left|p^2n-pn\right|=p^2n-pn.\]
 Therefore, $LE^+(\Gamma_G) = ((p^2-1)n-(p+1)) \times 0 + p \times (pn-n) + p^2n-pn= 2pn(p-1)$.
 In particular, if $G$ is non-abelian and $|G|=p^3$ then $n=p$. Therefore, $LE^+(\Gamma_G) = 2p^2(p-1)$.
\end{proof}

\begin{theorem}\label{Zp-energetic}
If $\frac{G}{Z(G)} \cong \mathbb{Z}_p \times \mathbb{Z}_p$ then 
\begin{enumerate}
\item $E(\Gamma_G)= LE(\Gamma_G)= LE^+(\Gamma_G)$.
\item $\Gamma_G$ is non-hypoenergetic, non-hyperenergetic, not L-hyperenergetic and not Q-hyperenergetic.
\end{enumerate} 

\noindent In particular, if $G$ is non-abelian and $|G|=p^3$ then $E(\Gamma_G)= LE(\Gamma_G)= LE^+(\Gamma_G)$ and $\Gamma_G$ is non-hyperenergetic, non-hypoenergetic, not L-hyperenergetic as well as not Q-hyperenergetic.
\end{theorem}

\begin{proof}
(a) For $\frac{G}{Z(G)} \cong \mathbb{Z}_p \times \mathbb{Z}_p$, from Theorems \ref{Z_p} and \ref{Z_p*Z_p}, we have $E(\Gamma_G) = LE(\Gamma_G) = LE^+(\Gamma_G)= 2p^{2n}(p^n-1)$ and hence follows.

\vspace{.5cm}

\noindent (b) Here, $|v(\Gamma_G)|=(p^2-1)|Z(G)|$. Thus, $E(K_{|v(\Gamma_G)|}) = LE(K_{|v(\Gamma_G)|})= LE^+(K_{|v(\Gamma_G)|})= 2p(p-1)|Z(G)|+2(p-1)|Z(G)|-2$.
Therefore, by Theorems \ref{Z_p} and \ref{Z_p*Z_p}, we have  $E(K_{|v(\Gamma_G)|})-E(\Gamma_G) = LE(K_{|v(\Gamma_G)|})-LE(\Gamma_G) = LE^+(K_{|v(\Gamma_G)|})-LE^+(\Gamma_G) = 2(p-1)|Z(G)|-2 > 0$. Also, $E(\Gamma_G)-|v(\Gamma_G)|=(p-1)^2|Z(G)|>0$. Hence, the results follow.

 In particular, if $G$ is non-abelian and $|G|=p^3$, then $|Z(G)|=p$ in the above cases so the results hold.
\end{proof}




\section{$\frac{G}{Z(G)}$ is isomorphic to $Sz(2)$}
We compute  spectrum, energy, Signless Laplacian Spectrum and Signless Laplacian energy of $\Gamma_G$ considering the group $G$ whose central quotient is isomorphic to the Suzuki group of order $20$ denoted by $Sz(2)$. Further we compare energy, Laplacian energy and Signless Laplacian energy of $\Gamma_G$ and look into the hyper- and hypo-properties of $\Gamma_G$.
\begin{theorem}[\protect{\cite[Theorem 2.1]{DN18}}]\label{Sz}
Let $\frac{G}{Z(G)} \cong Sz(2)$. Then $LE(\Gamma_G) = \left(\frac{120}{19}n+30\right)n$, where $|Z(G)| = n$.
\end{theorem}

\begin{theorem}\label{Suzuki}
Let $\frac{G}{Z(G)} \cong Sz(2)$. Then
\begin{enumerate}
\item $\spec(\Gamma_G)=\left\lbrace 0^{19n-6},(-3n)^4,\left(2n\left(3+2\sqrt{6}\right)\right)^1,\left(2n\left(3-2\sqrt{6}\right)\right)^1\right\rbrace$ and $E(\Gamma_G)=4n(3+2\sqrt{6})$, where $n=|Z(G)|$.
\item $\Q-spec(\Gamma_G)=\left\{(16n)^{15n-5},(15n)^{4n-1},(13n)^4,
\left(\frac{n(43+\sqrt{409})}{2}\right)^1,\left(\frac{n(43-\sqrt{409})}{2}\right)^1\right\}$ and\\
$LE^+(\Gamma_G) =\frac{120n^2+177n}{19}+\sqrt{409}n$, where $n=|Z(G)|$.
\end{enumerate}
\end{theorem}
\begin{proof}
If $\frac{G}{Z(G)} \cong Sz(2)$ and $|Z(G)| = n$ then $\Gamma_G=K_{5.3n,1.4n}$ and it is a complete 6-partite graph with $19n$ vertices.

\noindent (a) Using Theorem \ref{R1}(a), the characteristic polynomial of $\Gamma_G$ is
\[
P_{\Gamma_G}(x)=x^{19n-6}(x+3n)^4(x^2-12nx-60n^2).
\]
Therefore, $\spec(\Gamma_G)=\left\lbrace 0^{19n-6},(-3n)^4,\left(2n\left(3+2\sqrt{6}\right)\right)^1,\left(2n\left(3-2\sqrt{6}\right)\right)^1\right\rbrace$ and $E(\Gamma_G)=4n\left(3+2\sqrt{6}\right)$.

\vspace{.5cm}

\noindent (b) Using Theorem \ref{R1}(b), we have 
\begin{align*}
Q_{\Gamma_G}(x)=&{\displaystyle \prod_{i=1}^2} (x-19n+p_i)^{a_i(p_i-1)}{\displaystyle \prod_{i=1}^2}(x-19n+2p_i)^{a_i}\left(1-{\displaystyle \sum_{i=1}^2}\frac{a_ip_i}{x-19n+2p_i}\right)\\
  = \,& (x-19n+3n)^{5(3n-1)}(x-19n+4n)^{4n-1}(x-19n+2\times 3n)^{5}(x-19n+2\times 4n)^{1}\\
  \,& \times \left(1-\frac{5\times 3n}{x-19n+2\times 3n}-\frac{4n}{x-19n+2\times 4n}\right)\\
   =\,& (x-16n)^{15n-5}(x-15n)^{4n-1}(x-13n)^{4}(x^2-43nx+360n^2).
\end{align*}
Thus, $\Q-spec(\Gamma_G)=\left\{(16n)^{15n-5},(15n)^{4n-1},(13n)^4,\left(\frac{n(43+\sqrt{409})}{2}\right)^1,\left(\frac{n(43-\sqrt{409})}{2}\right)^1\right\}$.

 Number of edges of $\Gamma_G^c$ is $\frac{61n^2-19n}{2}$. Thus, $|e(\Gamma_G)| = \frac{19n(19n-1)}{2}-\frac{61n^2-19n}{2}=150n^2$. Now 
\[\left|16n - \frac{2|e(\Gamma_G)|}{|v(\Gamma_G)|}\right| = \left|\frac{4n}{19}\right|= \frac{4n}{19},\quad
\left|15n - \frac{2|e(\Gamma_G)|}{|v(\Gamma_G)|}\right| = \left|\frac{-15n}{19}\right|  =\frac{15n}{19},\]
\[\left|13n- \frac{2|e(\Gamma_G)|}{|v(\Gamma_G)|}\right| = \left|\frac{-53n}{19}\right|= \frac{53n}{19},\]
\[\left|\frac{(43+\sqrt{409})n}{2} - \frac{2|e(\Gamma_G)|}{|v(\Gamma_G)|}\right| =\left|\frac{(217+19\sqrt{409})n}{38}\right|=\frac{(217+19\sqrt{409})n}{38}\] and
\[\left|\frac{(43-\sqrt{409})n}{2} - \frac{2|e(\Gamma_G)|}{|v(\Gamma_G)|}\right|= \left|\frac{(217-19\sqrt{409})n}{38}\right| = \frac{(19\sqrt{409}-217)n}{38}.\]
Therefore, $LE^+(\Gamma_G) = (15n-5) \times \frac{4n}{19} + (4n-1) \times  \frac{15n}{19}+ 4 \times \frac{53n}{19} +\frac{(217+19\sqrt{409})n}{38}+\frac{(19\sqrt{409}-217)n}{38}$ and the result follows on simplification.
\end{proof}

\begin{theorem}{\label{SSzz}}
If $\frac{G}{Z(G)}$ is isomorphic to $Sz(2)$ then
\begin{enumerate}
\item $E(\Gamma_G) < LE^+(\Gamma_G) < LE(\Gamma_G)$.
\item $\Gamma_G$ is non-hypoenergetic as well as non-hyperenergetic.
\item $\Gamma_{Sz(2)}$ is not Q-hyperenergetic but is L-hyperenergetic. If $G \ncong Sz(2)$ then $\Gamma_G$ is Q-hyperenergetic and L-hyperenergetic. 
\end{enumerate}
\end{theorem}
\begin{proof} (a) Using Theorems \ref{Suzuki} and \ref{Sz}, we have $LE(\Gamma_G)-LE^+(\Gamma_G)=\left(\frac{393}{19}-\sqrt{409}\right)n > 0$ and $LE^+(\Gamma_G)-E(\Gamma_G)=\frac{3n(40n-17)}{19}+(\sqrt{409}-8\sqrt{6})n > 0$, where $n=|Z(G)|$. Hence, the result follows.
	
\vspace{.5cm}	

\noindent (b) Here, $|v(\Gamma_G)|=19n$, $n=|Z(G)|$ and $E(K_{|v(\Gamma_G)|})=LE(K_{|v(\Gamma_G)|})= LE^+(K_{|v(\Gamma_G)|}) =38n-2$. Using Theorem \ref{Suzuki}(a), we have $E(\Gamma_G)-|v(\Gamma_G)|=(8\sqrt{6}-7)n>0$ and also $E(K_{|v(\Gamma_G)|})-E(\Gamma_G) = 2(13-4\sqrt{6})n-2>0$. Hence, the result follows.

\vspace{.5cm}	

\noindent (c)  For $n = |Z(G)| = 1$, from Proposition 4.3.13 of \cite{FWNT21}, we have $LE(\Gamma_{Sz(2)})=\frac{690}{19}>36= LE(K_{|v(\Gamma_{Sz(2)})|})$. Hence, for $n = 1$, $\Gamma_G$ is L-hyperenergetic.

For $n = |Z(G)| > 1$, using Theorem \ref{Suzuki}(b), we have $LE^+(\Gamma_G)-LE^+(K_{|v(\Gamma_G)|}) = \frac{5n(24n-109)+38}{19}+n\sqrt{409}>0$. Therefore,   $LE^+(\Gamma_G)>LE^+(K_{|v(\Gamma_G)|})$ which implies $\Gamma_G$ is Q-hyperenergetic and consequently part (a) implies $\Gamma_G$ is L-hyperenergetic.
\end{proof}

  For $\frac{G}{Z(G)} \cong Sz(2)$, the following figures also demonstrate that among the three energies, $E(\Gamma_G)$ is the least and the fact that although $LE^+(\Gamma_G) < LE(\Gamma_G)$ but these two energies are very close to each other. 

\vspace{.3cm}

\begin{minipage}[t]{.5\linewidth}
\begin{tikzpicture}
\begin{axis}
[
xlabel={$n$ $\rightarrow$},
ylabel={Energies of $\Gamma_G$ $\rightarrow$},
xmin=2, xmax=10,
ymin=0, ymax=1000,
grid = both,
minor tick num = 1,
major grid style = {lightgray},
minor grid style = {lightgray!25},
width=.7\textwidth,
height=.7\textwidth,
legend style={legend pos=north west},
 ]
\addplot[domain=2:15,samples at={2,3,4,...,15},mark=*,green, samples=24, mark size=.8pt]{4*x*(3+2*(6)^(1/2))};
\tiny
\addlegendentry{$E$}
\addplot[domain=2:15,samples at={2,3,4,...,15},mark=triangle*,blue,mark size=.8pt, samples=24]{(120*x*x)/(19)+30*x};
\tiny
\addlegendentry{$LE$}
\addplot[domain=2:15,samples at={2,3,4,...,15},mark=square*, red, mark size=.8pt, samples=24]{(120*x*x+177*x)/(19)+(409)^(1/2)*x};
\tiny
\addlegendentry{$LE^+$}
\end{axis}
\end{tikzpicture}
\vspace{-.2 cm}
\captionsetup{font=footnotesize}

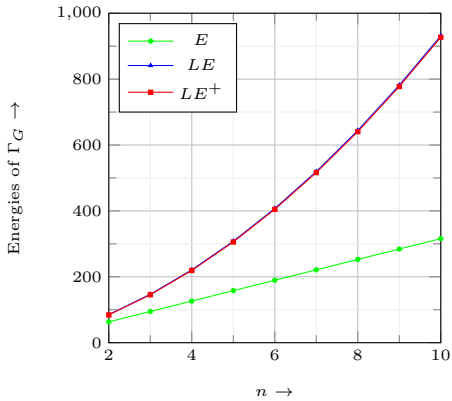
\captionof{figure}{Energies of $\Gamma_G$ where $\frac{G}{Z(G)} \cong Sz(2)$}
\end{minipage}
\hspace{0.01cm}
\begin{minipage}[t]{.5\linewidth}
\includegraphics[width=5cm]{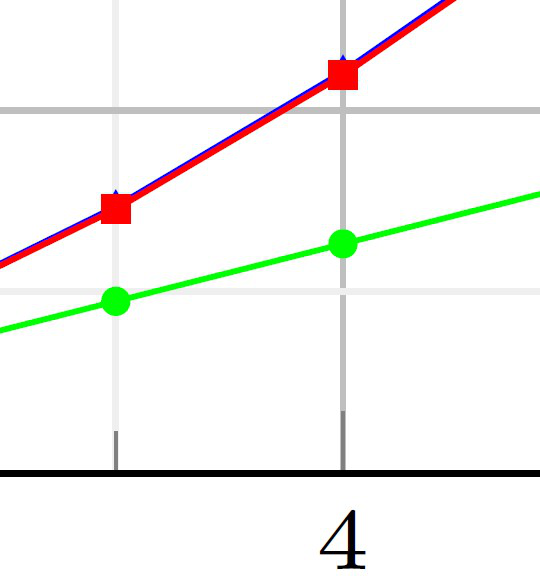}
\vspace{-.2 cm}
\captionsetup{font=footnotesize}
{\captionof{figure}{A close up view of Figure 11}}
\end{minipage}

\section{Some more classes of groups}
In this section we discuss   results on energy, Laplacian energy and Signless Laplacian energy of non-commuting graph of certain well-known classes of finite groups.

\subsection{The Hanaki groups}
We consider the Hanaki groups 
\[A(n,  \mathcal{V}) =\left\{ 
U(a,b)= \left[\begin{array}{ccc}
	1 & 0 & 0\\
	a & 1 & 0\\
	b & \mathcal{V}(a) & 1  
\end{array}\right]: a,b \in GF(2^n) \right\} \quad (n \geq 2),
\]
under matrix multiplication given by $U(a, b)U(a', b') = U(a + a', b + b' + a' \mathcal{V}(a))$ (here $\mathcal{V}$ be the Frobenius automorphism of $GF(2^n)$, i.e., $\mathcal{V}(x) = x^2$ for all $x \in GF(2^n)$) and 
\[A(n,p)= \left\{ 
V(a,b,c)= \left[\begin{array}{ccc}
	1 & 0 & 0\\
	a & 1 & 0\\
	b & c & 1  
\end{array}\right]: a,b,c \in GF(p^n) \right\} (p \text{ is any prime})
\]
under matrix multiplication $V (a, b, c)V (a', b', c') = V (a + a', b + b' + ca', c + c')$. 
 In this section, we compute Signless Laplacian spectrum and Signless Laplacian energy of non-commuting graph of the  groups $A(n,  \mathcal{V})$ and $A(n,p)$. Further we compare Signless Laplacian energy of $\Gamma_G$ with its predetermined energy, Laplacian energy  and look into the hyper- and hypo-properties of $\Gamma_G$ if $G$ is isomorphic to $A(n,  \mathcal{V})$ and $A(n,p)$. 
\begin{theorem}[\protect{\cite[Proposition 3.5]{DN18} and \cite[Theorem 4.1.3]{FWNT21}}]\label{H1}
If $G$ denotes the Hanaki group  $A(n,  \mathcal{V})$, 
then
\[E(\Gamma_G)=LE(\Gamma_G) = 2^{2n+1}-2^{n+2}.\]
\end{theorem}

\begin{theorem}\label{Hanaki1}
If $G$ is isomorphic to the Hanaki group $A(n,  \mathcal{V})$ then $\Q-spec(\Gamma_G)=\{(2^{2n}-2^{n+1})^{2^{2n}-2^{n+1}+1},\\(2^{2n}-3\times 2^n)^{2^n-2},(2^{2n+1}-2^{n+2})^1\}$ and $LE^+(\Gamma_G) = 2^{2n+1}-2^{n+2}$.
\end{theorem}

\begin{proof}
If $G$ is isomorphic to the Hanaki group $A(n,  \mathcal{V})$  then  $|v(\Gamma_G)|=2^{2n}-2^n$ and $\Gamma_G=K_{(2^n-1).2^n}$. Using Theorem \ref{R1}(b), we have 
\begin{align*}
Q_{\Gamma_G}(x)= & (x-(2^{2n}-2^n)+2^n)^{(2^n-1)(2^n-1)}(x-(2^{2n}-2^n)+2\times 2^n)^{(2^n-1)}\left(1-\frac{(2^n-1).2^n}{x-(2^{2n}-2^n)+2\times 2^n}\right)\\
   =\,& (x-(2^{2n}-2^{n+1}))^{(2^n-1)^2}(x-(2^{2n}-3\times 2^n))^{2^n-2}(x-(2^{2n+1}-2^{n+2})).
\end{align*}
Thus, $\Q-spec(\Gamma_G)=\{(2^{2n}-2^{n+1})^{2^{2n}-2^{n+1}+1},(2^{2n}-3\times 2^n)^{2^n-2},(2^{2n+1}-2^{n+2})^1\}$.

 Number of edges of $\Gamma_G^c$ is $2^{n-2}(2^{2n}-2^{n+1}+1)$. Thus, $|e(\Gamma_G)| = \frac{(2^{2n}-2^n)(2^{2n}-2^n-1)}{2}-2^{n-2}(2^{2n}-2^{n+1}+1)$ $=2^{4n-1}-3\times 2^{3n-1}+2^{2n}$. Now,
\[\left|2^{2n}-2^{n+1} - \frac{2|e(\Gamma_G)|}{|v(\Gamma_G)|}\right| = \left|\frac{2^{3n+1}-2^{3n+1}}{2^{2n}-2^n}\right|= 0,\quad
\left|2^{2n}-3\times 2^n - \frac{2|e(\Gamma_G)|}{|v(\Gamma_G)|}\right| = \left|\frac{2^{2n}-2^{3n}}{2^{2n}-2^n}\right|  =\frac{2^{3n}-2^{2n}}{2^{2n}-2^n}\] 
and
\[\left|2^{2n+1}-2^{n+2}- \frac{2|e(\Gamma_G)|}{|v(\Gamma_G)|}\right| =
\left|\frac{2^{4n}-3\times 2^{3n}+2^{2n+1}}{2^{2n}-2^n}\right|=\frac{2^{4n}-3\times 2^{3n}+2^{2n+1}}{2^{2n}-2^n}.\]
Therefore,
$LE^+(\Gamma_G) = (2^{2n}-2^{n+1}+1) \times 0 + (2^n-2) \times  \frac{2^{3n}-2^{2n}}{2^{2n}-2^n}+ \frac{2^{4n}-3\times 2^{3n}+2^{2n+1}}{2^{2n}-2^n}$ and the result follows on simplification.
\end{proof}

\begin{theorem}
If $G$ is isomorphic to the Hanaki group $A(n,  \mathcal{V})$ then
\begin{enumerate}
\item $E(\Gamma_G) = LE(\Gamma_G)=LE^+(\Gamma_G)$.
\item $\Gamma_G$ is non-hypoenergetic, non-hyperenergetic, not L-hyperenergetic and not Q-hyperenergetic.
\end{enumerate}
\end{theorem}
\begin{proof}
(a) Using Theorems \ref{H1} and \ref{Hanaki1}, we have $E(\Gamma_G) = LE(\Gamma_G) = LE^+(\Gamma_G) = 2^{n+1}(2^n-2)$ and hence the result follows.

\vspace{0.5cm}
\noindent (b) Here, $|v(\Gamma_G)|=2^n(2^n-1)$ and $E(K_{|v(\Gamma_G)|}) = LE(K_{|v(\Gamma_G)|})= LE^+(K_{|v(\Gamma_G)|}) = 2^{n+1}(2^n-1)-2$. Therefore, by Theorems \ref{H1} and \ref{Hanaki1}, we have $E(\Gamma_G)-|v(\Gamma_G)|=2^n(2^n-3)>0$, for $n>1$ and $E(K_{|v(\Gamma_G)|})-E(\Gamma_G)= LE(K_{|v(\Gamma_G)|})-LE(\Gamma_G)= LE^+(K_{|v(\Gamma_G)|})-LE^+(\Gamma_G)= 2(2^n-1)>0$. Hence, the results follow.
\end{proof}


\begin{theorem}[\protect{\cite[Proposition 3.6]{DN18} and \cite[Theorem 4.1.4]{FWNT21}}]\label{H2}
%
If $G$ denotes the Hanaki group $A(n, p)$ then 
\[E(\Gamma_G)=LE(\Gamma_G) = 2(p^{3n}-p^{2n}).\]
\end{theorem}

\begin{theorem}\label{Hanaki2}
If $G$ is isomorphic to the Hanaki group $A(n,p)$ then $\Q-spec(\Gamma_G)=\{(p^{3n}-p^{2n})^{(p^n+1)(p^{2n}-p^n-1)},\\(p^{3n}-2p^{2n}+p^n)^{p^n},(2p^{3n}-2p^{2n})^1\}$ and $LE^+(\Gamma_G) = 2p^{2n}(p^n-1)$.
\end{theorem}

\begin{proof}
If $G$ is isomorphic to the Hanaki group $A(n,p)$  then $|v(\Gamma_G)|=p^{3n}-p^n$ and $\Gamma_G=K_{(p^n+1).(p^{2n}-p^n)}$. Using Theorem \ref{R1}(b), we have 
\begin{align*}
Q_{\Gamma_G}(x)= & (x-(p^{3n}-p^n)+(p^{2n}-p^n))^{(p^n+1)(p^{2n}-p^n-1)}(x-(p^{3n}-p^n)+2(p^{2n}-p^n))^{(p^n+1)}\\
\,& \times \left(1-\frac{p^{3n}-p^n}{x-(p^{3n}-p^n)+2(p^{2n}-p^n)}\right)\\
   =\,& (x-(p^{3n}-p^{2n}))^{(p^n+1)(p^{2n}-p^n-1)}(x-(p^{3n}-2p^{2n}+p^n))^{p^n}(x-(2p^{3n}-2p^{2n})).
\end{align*}
Thus, $\Q-spec(\Gamma_G)=\{(p^{3n}-p^{2n})^{(p^n+1)(p^{2n}-p^n-1)},(p^{3n}-2p^{2n}+p^n)^{p^n},(2p^{3n}-2p^{2n})^1\}$.

 Number of edges of $\Gamma_G^c$ is $\frac{(p^{3n}-p^n)(p^{2n}-p^n-1)}{2}$. Thus, $|e(\Gamma_G)| = \frac{(p^{3n}-p^n)(p^{3n}-p^n-1)}{2}-\frac{(p^{3n}-p^n)(p^{2n}-p^n-1)}{2}$ $=\frac{p^{6n}-p^{5n}-p^{4n}+p^{3n}}{2}$. Now 
\[\left|p^{3n}-p^{2n} - \frac{2|e(\Gamma_G)|}{|v(\Gamma_G)|}\right| = \left|\frac{0}{p^{3n}-p^n}\right|= 0,\]
\[\left|p^{3n}-2p^{2n}+p^n - \frac{2|e(\Gamma_G)|}{|v(\Gamma_G)|}\right| = \left|\frac{-p^{5n}+p^{4n}+p^{3n}-p^{2n}}{p^{3n}-p^n}\right|  =\frac{p^{5n}-p^{4n}-p^{3n}+p^{2n}}{p^{3n}-p^n}\]
 and
\[\left|2p^{3n}-2p^{2n}- \frac{2|e(\Gamma_G)|}{|v(\Gamma_G)|}\right| =
\left|\frac{p^{6n}-p^{5n}-p^{4n}+p^{3n}}{p^{3n}-p^n}\right|=\frac{p^{6n}-p^{5n}-p^{4n}+p^{3n}}{p^{3n}-p^n}.\]
Therefore, $LE^+(\Gamma_G) = (p^n+1)(p^{2n}-p^n-1)\times 0 + p^n \times \frac{p^{5n}-p^{4n}-p^{3n}+p^{2n}}{p^{3n}-p^n}+ \frac{p^{6n}-p^{5n}-p^{4n}+p^{3n}}{p^{3n}-p^n}$ and the result follows on simplification.
\end{proof}

\begin{theorem}
If $G$ is isomorphic to the Hanaki group $A(n, p)$ then
\begin{enumerate}
\item $E(\Gamma_G) = LE(\Gamma_G)=LE^+(\Gamma_G)$.
\item $\Gamma_G$ is non-hypoenergetic, non-hyperenergetic, not L-hyperenergetic and not Q-hyperenergetic.
\end{enumerate}
\end{theorem}
\begin{proof}
(a) Using Theorems \ref{H2} and \ref{Hanaki2} we have $E(\Gamma_G) = LE(\Gamma_G) = LE^+(\Gamma_G)= 2p^{2n}(p^n-1)$ and hence the result follows.

\vspace{0.5cm}
\noindent (b) Here, $|v(\Gamma_G)|=p^{3n}-p^n$ and $E(K_{|v(\Gamma_G)|}) = LE(K_{|v(\Gamma_G)|})= LE^+(K_{|v(\Gamma_G)|})= 2(p^{2n}(p^n-1)+p^{2n}-p^n-1)$. Therefore, by Theorems \ref{H2} and \ref{Hanaki2}, we have $E(\Gamma_G)-|v(\Gamma_G)|=p^n(p^n-1)^2>0$ and $E(K_{|v(\Gamma_G)|})-E(\Gamma_G)= LE(K_{|v(\Gamma_G)|})-LE(\Gamma_G)= LE^+(K_{|v(\Gamma_G)|})-LE^+(\Gamma_G)= 2(p^n(p^n-1)-1)>0$. Hence, the results follow.
\end{proof}


\subsection{The semi-dihedral groups, $SD_{8n}$}
We consider $SD_{8n}:=\langle a,b : a^{4n} = b^2 = 1, bab^{-1} = a^{2n-1} \rangle$, the semi-dihedral groups of order $8n$ (where $n > 1$). Results regarding different energies of non-commuting graph of $SD_{8n}$ are given below.
\begin{theorem}[\protect{\cite[Theorem 4.1.1 (b) and Theorem 4.3.2]{FWNT21}}]\label{SD}
Let $G$ be isomorphic to $SD_{8n}$.
\begin{enumerate}
\item If $n$ is odd then
\[E(\Gamma_{SD_{8n}}) = 4(n-1)+4\sqrt{(n-1)(5n-1)} \, \text{ and } \,
LE(\Gamma_{SD_{8n}}) = \frac{8n(4n^2-10n+7)}{2n-1}.\]
\item If $n$ is even then
\[E(\Gamma_{SD_{8n}}) = 2(2n-1)+2\sqrt{(2n-1)(10n-1)} \, \text{ and } \,
LE(\Gamma_{SD_{8n}}) =  \frac{8n(8n^2-8n+3)}{4n-1}.\]
\end{enumerate}
\end{theorem}


\begin{theorem}\label{SD_8n1}
Let $G$ be isomorphic to $SD_{8n}$, where  $n$ is odd. Then
\begin{align*}
\Q-spec(\Gamma_{SD_{8n}})=&\left\lbrace(8n-8)^{3n},(4n)^{4n-5},(8n-12)^{n-1},\left(8n-6+2\sqrt{8n^2-16n+9}\right)^1,\right.\\
&\left.\left(8n-6-2\sqrt{8n^2-16n+9}\right)^1\right\rbrace    
\end{align*} and 
\begin{align*}
LE^+(\Gamma_{SD_{8n}}) = \begin{cases}36+4\sqrt{33}, & \mbox{if $n=3$}\vspace{.2cm}\\
\frac{32n^3-112n^2+96n-12}{2n-1}+4\sqrt{8n^2-16n+9}, & \mbox{if $n \geq 5.$}\\
\end{cases}   
\end{align*}
\end{theorem}
\begin{proof}
If $G\cong SD_{8n}$ and $n$ is odd then $|v(\Gamma_{SD_{8n}})|=8n-4$ and $\Gamma_{SD_{8n}}=K_{n.4,1.(4n-4)}$. Using Theorem \ref{R1}(b), we have 
\begin{align*}
Q_{\Gamma_{SD_{8n}}}(x)=&{\displaystyle \prod_{i=1}^2} (x-(8n-4)+p_i)^{a_i(p_i-1)}{\displaystyle \prod_{i=1}^2}(x-(8n-4)+2p_i)^{a_i}\left(1-{\displaystyle \sum_{i=1}^2}\frac{a_ip_i}{x-(8n-4)+2p_i}\right)\\
  = \,& (x-8n+8)^{3n}(x-4n)^{4n-5}(x-8n+12)^n(x-4)\left(1-\frac{4n}{x-8n+12}-\frac{4n-4}{x-4}\right)\\
   = \,& (x-(8n-8))^{3n}(x-4n)^{4n-5}(x-(8n-12))^{n-1}(x^2-(16n-12)x+32n^2-32n).
\end{align*}
Thus, $\Q-spec(\Gamma_{SD_{8n}})=\left\lbrace(8n-8)^{3n},(4n)^{4n-5},(8n-12)^{n-1},\left(8n-6+2\sqrt{8n^2-16n+9}\right)^1,\right.\\
 ~~~~~~~~~~~~~~~~~~~~~~~~~~~~~~~~~\left.\left(8n-6-2\sqrt{8n^2-16n+9}\right)^1\right\rbrace$.

 Number of edges of $\Gamma_{SD_{8n}}^c$ is $8n^2-12n+10$. Thus, $|e(\Gamma_{SD_{8n}})| = \frac{(8n-4)(8n-4-1)}{2}-(8n^2-12n+10)=24n(n-1)$. Now 
\[\left|8n-8 - \frac{2|e(\Gamma_{SD_{8n}})|}{|v(\Gamma_{SD_{8n}})|}\right| = \left|\frac{4(n-1)(n-2)}{2n-1}\right| = \frac{4(n-1)(n-2)}{2n-1},\] 
\[\left|4n - \frac{2|e(\Gamma_{SD_{8n}})|}{|v(\Gamma_{SD_{8n}})|}\right| = \left|\frac{4n(2-n)}{2n-1}\right| =
\frac{4n(n-2)}{2n-1},\]
\[\left|8n-12 - \frac{2|e(\Gamma_{SD_{8n}})|}{|v(\Gamma_{SD_{8n}})|}\right| = \left|\frac{4(n^2-5n+3)}{2n-1}\right| = \begin{cases}\frac{4(-n^2+5n-3)}{2n-1}, & \mbox{if $n \leq 3$}\vspace{.2cm}\\
\frac{4(n^2-5n+3)}{2n-1}, & \mbox{if $n \geq 5,$}
 \end{cases}\]
\begin{align*}
\left|8n-6+2\sqrt{8n^2-16n+9} - \frac{2|e(\Gamma_{SD_{8n}})|}{|v(\Gamma_{SD_{8n}})|}\right| =& \left|2\sqrt{8n^2-16n+9}+(2n-3)+\frac{3}{2n-1}\right|\\ = & 2\sqrt{8n^2-16n+9}+(2n-3)+\frac{3}{2n-1}
\end{align*} and
\begin{align*}
\left|8n-6-2\sqrt{8n^2-16n+9} - \frac{2|e(\Gamma_{SD_{8n}})|}{|v(\Gamma_{SD_{8n}})|}\right| = &\left|-2\sqrt{8n^2-16n+9}+(2n-3)+\frac{3}{2n-1}\right|\\= & 2\sqrt{8n^2-16n+9}-(2n-3)-\frac{3}{2n-1}.\end{align*}


Therefore, for $n=3$ we have $LE^+(\Gamma_{SD_{8n}})= 36+4\sqrt{33}$. For $n \geq 5$ we have 
\begin{align*}
LE^+(\Gamma_{SD_{8n}}) = & 3n\times \frac{4(n-1)(n-2)}{2n-1}+(4n-5)\times \frac{4n(n-2)}{2n-1}+(n-1)\times\frac{4(n^2-5n+3)}{2n-1}+\\ &2\sqrt{8n^2-16n+9}-2n-3+\frac{3}{2n-1}+2\sqrt{8n^2-16n+9}+2n+3-\frac{3}{2n-1}    
\end{align*} 
and the result follows on simplification.
\end{proof}

\begin{theorem}\label{SD_8n2}
Let $G$ be isomorphic to $SD_{8n}$, where $n$ is even. Then
\begin{align*}
\Q-spec(\Gamma_{SD_{8n}})=&\left\lbrace (8n-4)^{2n},(4n)^{4n-3},(8n-6)^{2n-1}, \left(8n-3+\sqrt{32n^2-32n+9}\right)^1,\right.\\
&\left.\left(8n-3-\sqrt{32n^2-32n+9}\right)^1\right\rbrace    
\end{align*} and
\begin{align*}
LE^+(\Gamma_{SD_{8n}}) = \begin{cases} \frac{134}{7}+2\sqrt{73}, & \mbox{if $n=2$}\vspace{.2cm}\\ \frac{64n^3-128n^2+64n-6}{4n-1}+2\sqrt{32n^2-32n+9}, & \mbox{if $n \geq 4$.}\\
\end{cases}  
\end{align*}
\end{theorem}

\begin{proof}
If $G\cong SD_{8n}$ and $n$ is even then  $|v(\Gamma_{SD_{8n}})|=8n-2$ and $\Gamma_{SD_{8n}}=K_{2n.2,1.(4n-2)}$. Using Theorem \ref{R1}(b), we have 
\begin{align*}
Q_{\Gamma_{SD_{8n}}}(x)=&{\displaystyle \prod_{i=1}^2} (x-(8n-2)+p_i)^{a_i(p_i-1)}{\displaystyle \prod_{i=1}^2}(x-(8n-2)+2p_i)^{a_i}\left(1-{\displaystyle \sum_{i=1}^2}\frac{a_ip_i}{x-(8n-2)+2p_i}\right)\\
  = \,& (x-8n+4)^{2n}(x-4n)^{4n-3}(x-8n+6)^{2n}(x-2)\left(1-\frac{4n}{x-8n+6}-\frac{4n-2}{x-2}\right)\\
   = \,& (x-(8n-4))^{2n}(x-4n)^{4n-3}(x-(8n-6))^{2n-1}(x^2-(16n-6)x+32n^2-16n).
\end{align*}
Thus, $\Q-spec(\Gamma_{SD_{8n}})=\left\lbrace(8n-4)^{2n},(4n)^{4n-3},(8n-6)^{2n-1}, \left(8n-3+\sqrt{32n^2-32n+9}\right)^1,\right.\\
 ~~~~~~~~~~~~~~~~~~~~~~~~~~~~~~~~~\left.\left(8n-3-\sqrt{32n^2-32n+9}\right)^1\right\rbrace$.

 Number of edges of $\Gamma_{SD_{8n}}^c$ is $8n^2-8n+3$. Therefore, $|e(\Gamma_{SD_{8n}})| = \frac{(8n-2)(8n-2-1)}{2}-(8n^2-8n+3)$ $=12n(2n-1)$. Now
\[\left|8n-4 - \frac{2|e(\Gamma_{SD_{8n}})|}{|v(\Gamma_{SD_{8n}})|}\right| = \left|\frac{(8n-4)(n-1)}{4n-1}\right| = \frac{(8n-4)(n-1)}{4n-1},\] 
\[\left|4n - \frac{2|e(\Gamma_{SD_{8n}})|}{|v(\Gamma_{SD_{8n}})|}\right| =\left|\frac{-8n(n-1)}{4n-1}\right| = \frac{8n(n-1)}{4n-1},\]
\[\left|8n-6 - \frac{2|e(\Gamma_{SD_{8n}})|}{|v(\Gamma_{SD_{8n}})|}\right|= \left|\frac{(8n^2-20n+6)}{4n-1}\right| = \begin{cases}\frac{2}{7}, & \mbox{if $n=2$} \vspace{.2cm}\\
  \frac{(8n^2-20n+6)}{4n-1}, & \mbox{if $n \geq 4$,}\end{cases}\]
\begin{align*}
 \left|8n-3+\sqrt{32n^2-32n+9} - \frac{2|e(\Gamma_{SD_{8n}})|}{|v(\Gamma_{SD_{8n}})|}\right| =&\left|\sqrt{32n^2-32n+9}+2n-\frac{3}{2}+\frac{3}{8n-2}\right|\\= & \sqrt{32n^2-32n+9}+2n-\frac{3}{2}+\frac{3}{8n-2} 
\end{align*} 
and
\begin{align*}\left|8n-3-\sqrt{32n^2-32n+9} - \frac{2|e(\Gamma_{SD_{8n}})|}{|v(\Gamma_{SD_{8n}})|}\right|  =&\left|-\sqrt{32n^2-32n+9}+2n-\frac{3}{2}+\frac{3}{8n-2}\right|\\=& \sqrt{32n^2-32n+9}-2n+\frac{3}{2}-\frac{3}{8n-2}.
\end{align*}
Therefore, for $n=2$, we have $LE^+(\Gamma_{SD_{8n}})= \frac{134}{7}+2\sqrt{73}$. For $n \geq 4$, we have
\begin{align*}
 LE^+(\Gamma_{SD_{8n}}) =& 2n\times \frac{(8n-4)(n-1)}{4n-1}+(4n-3)\times \frac{8n(n-1)}{4n-1}+(2n-1)\times\frac{(8n^2-20n+6)}{4n-1}+\\ &\sqrt{32n^2-32n+9}+2n-\frac{3}{2}+\frac{3}{8n-2}+\sqrt{32n^2-32n+9}-2n+\frac{3}{2}-\frac{3}{8n-2}   
\end{align*}and the result follows on simplification.
\end{proof}

\begin{theorem}\label{SD_{8n}}
If $G$ is isomorphic to $SD_{8n}$ then 
\begin{enumerate}
\item $E(\Gamma_{SD_{8n}})< LE^+(\Gamma_{SD_{8n}}) < LE(\Gamma_{SD_{8n}})$.
\item $\Gamma_{SD_{8n}}$ is non-hypoenergetic as well as non-hyperenergetic.
\item $\Gamma_{SD_{8n}}$ is Q-hyperenergetic and L-hyperenergetic.
\end{enumerate}
\end{theorem}

\begin{proof}
(a) \textbf{Case 1:} $n$ is odd

For $n=3$, using Theorems \ref{SD} and \ref{SD_8n1}, we have $E(\Gamma_{SD_{24}})=8+8\sqrt{7}$, $LE(\Gamma_{SD_{24}})= \frac{312}{5}$ and $LE^+(\Gamma_{SD_{24}})=36+4\sqrt{33}$. Clearly, $E(\Gamma_{SD_{24}})< LE^+(\Gamma_{SD_{24}}) < LE(\Gamma_{SD_{24}})$.

For $n \geq 5$, using Theorems \ref{SD} and \ref{SD_8n1}, we have 
\begin{equation}\label{S1}
LE(\Gamma_{SD_{8n}})-LE^+(\Gamma_{SD_{8n}})=\frac{32n^2-40n+12}{2n-1}-4\sqrt{8n^2-16n+9}
\end{equation}
and
\begin{equation}\label{S2}
LE^+(\Gamma_{SD_{8n}})-E(\Gamma_{SD_{8n}})=\frac{32n^3-116n^2+102n-14}{2n-1}+4\sqrt{8n^2-16n+9}-4\sqrt{5n^2-6n+1}.
\end{equation}
Since $32n^2-40n+12 > 0$, $4(2n-1)\sqrt{8n^2-16n+9} > 0$ and $(32n^2-40n+12)^2-\left(4\sqrt{8n^2-16n+9}\right)^2(2n-1)^2=512n^3(n-2)+128n(5n-1)>0$
we have 
$$
32n^2-40n+12 - 4(2n-1)\sqrt{8n^2-16n+9} > 0.
$$
 Therefore, by \eqref{S1}, $(2n-1)(LE(\Gamma_{SD_{8n}})-LE^+(\Gamma_{SD_{8n}})) > 0$. Hence, $LE(\Gamma_{SD_{8n}}) > LE^+(\Gamma_{SD_{8n}})$.  

Again, $\sqrt{8n^2-16n+9} > 0, \sqrt{5n^2-6n+1} > 0$ and $\left(\sqrt{8n^2-16n+9}\right)^2-\left(\sqrt{5n^2-6n+1}\right)^2=n(3n-10)+8>0$. Therefore, 
$\sqrt{8n^2-16n+9}-\sqrt{5n^2-6n+1}>0$. Since $32n^3-116n^2+102n-14 > 0$ we have $\frac{32n^3-116n^2+102n-14}{2n-1}+4\sqrt{8n^2-16n+9}-4\sqrt{5n^2-6n+1}>0$. Therefore, by \eqref{S2}, $LE^+(\Gamma_{SD_{8n}}) >  E(\Gamma_{SD_{8n}})$. Hence, 
$E(\Gamma_{SD_{8n}}) < LE^+(\Gamma_{SD_{8n}}) < LE(\Gamma_{SD_{8n}})$.


\vspace{0.5cm}
\noindent  \textbf{Case 2:} $n$ is even

For $n=2$, using Theorems \ref{SD} and \ref{SD_8n2}, we have $E(\Gamma_{SD_{16}})=6+2\sqrt{57}$, $LE(\Gamma_{SD_{16}})= \frac{304}{7}$ and $LE^+(\Gamma_{SD_{16}})=\frac{134}{7}+2\sqrt{73}$. Clearly, $E(\Gamma_{SD_{16}})< LE^+(\Gamma_{SD_{16}}) < LE(\Gamma_{SD_{16}})$.

For $n \geq 4$, using Theorems \ref{SD} and \ref{SD_8n2}, we have 
\begin{equation}\label{S3}
LE(\Gamma_{SD_{8n}})-LE^+(\Gamma_{SD_{8n}})=\frac{64n^2-40n+6}{4n-1}-2\sqrt{32n^2-32n+9}
\end{equation}
and
\begin{equation}\label{S4}
LE^+(\Gamma_{SD_{8n}})-E(\Gamma_{SD_{8n}})=\frac{16n^2(4n-9)+76n-8}{4n-1}+2\sqrt{32n^2-32n+9}-2\sqrt{20n^2-12n+1}.
\end{equation}
Since $64n^2-40n+6 > 0$, $2(4n-1)\sqrt{32n^2-32n+9} > 0$ and $(64n^2-40n+6)^2-(2\sqrt{32n^2-32n+9})^2(4n-1)^2=2048n^3(n-1)+64n(10n-1)>0$
we have $64n^2-40n+6-2(4n-1)\sqrt{32n^2-32n+9}> 0$. Therefore, by \eqref{S3}, $(4n-1)(LE^+(\Gamma_{SD_{8n}})-LE(\Gamma_{SD_{8n}})) > 0$. Hence, $LE(\Gamma_{SD_{8n}}) > LE^+(\Gamma_{SD_{8n}})$.  

Again, $\sqrt{32n^2-32n+9} > 0, \sqrt{20n^2-12n+1} > 0$ and $\left(\sqrt{32n^2-32n+9}\right)^2-\left(\sqrt{20n^2-12n+1}\right)^2=4n(3n-5)+8>0$. Therefore, 
$\sqrt{32n^2-32n+9}-\sqrt{20n^2-12n+1}>0$. Since $16n^2(4n-9)+76n-8 > 0$ we have $\frac{16n^2(4n-9)+76n-8}{2n-1}+2\sqrt{32n^2-32n+9}-2\sqrt{20n^2-12n+1}>0$. Therefore, by \eqref{S4}, $LE^+(\Gamma_{SD_{8n}}) >  E(\Gamma_{SD_{8n}})$. Hence, 
$E(\Gamma_{SD_{8n}}) < LE^+(\Gamma_{SD_{8n}}) < LE(\Gamma_{SD_{8n}})$.


\vspace{0.5cm}

\noindent (b) \textbf{Case 1:}  $n$ is odd 

Here $|v(\Gamma_{SD_{8n}})|=8n-4$ and $E(K_{|v(\Gamma_{SD_{8n}})|})= LE(K_{|v(\Gamma_{SD_{8n}})|})= LE^+(K_{|v(\Gamma_{SD_{8n}})|})= 16n-10$. Using Theorem \ref{SD}, we have
\begin{equation}\label{E35}
E(\Gamma_{SD_{8n}})-|v(\Gamma_{SD_{8n}})|= 4\sqrt{(n-1)(5n-1)}-4n     
\end{equation} and 
\begin{equation}\label{E36}
E(K_{|v(\Gamma_{SD_{8n}})|})-E(\Gamma_{SD_{8n}})= 12n-6-4\sqrt{(n-1)(5n-1)}.    
\end{equation}
Since $4\sqrt{(n-1)(5n-1)} > 0$, $4n > 0$ and $\left(4\sqrt{(n-1)(5n-1)}\right)^2-(4n)^2=16(4n^2-6n+1)>0$
we have $ 4\sqrt{(n-1)(5n-1)}-4n > 0$. Therefore, by \eqref{E35}, $E(\Gamma_{SD_{8n}}) > |v(\Gamma_{SD_{8n}})|$.

Again,  $4\sqrt{(n-1)(5n-1)} > 0$, $12n-6 > 0$ and $(12n-6)^2-\left(4\sqrt{(n-1)(5n-1)}\right)^2=4(16n^2-12n+5)>0$
and so $12n-6-4\sqrt{(n-1)(5n-1)}> 0$. Therefore, by \eqref{E36}, $E(K_{|v(\Gamma_{SD_{8n}})|}) > E(\Gamma_{SD_{8n}})$.

\vspace{0.5cm}
\noindent  \textbf{Case 2:} $n$ is even

Here $|v(\Gamma_{SD_{8n}})|=8n-2$ and $E(K_{|v(\Gamma_{SD_{8n}})|})= LE(K_{|v(\Gamma_{SD_{8n}})|})= LE^+(K_{|v(\Gamma_{SD_{8n}})|})= 16n-16$. Using Theorem \ref{SD}, we have
\begin{equation}\label{E37}
E(\Gamma_{SD_{8n}}) - |v(\Gamma_{SD_{8n}})|= 2\left(\sqrt{(2n-1)(10n-1)}-2n\right)
\end{equation} and
\begin{equation}\label{E38}
E(K_{|v(\Gamma_{SD_{8n}})|})-E(\Gamma_{SD_{8n}})= 2\left(3(2n-1)+1-\sqrt{(2n-1)(10n-1)}\right).
\end{equation}
Since $\sqrt{(2n-1)(10n-1)} > 0$, $2n > 0$ and $\left(\sqrt{(2n-1)(10n-1)}\right)^2-(2n)^2=4n(4n-3)+1>0$
we have $ \sqrt{(2n-1)(10n-1)}-2n> 0$. Therefore, by \eqref{E37}, $E(\Gamma_{SD_{8n}}) > |v(\Gamma_{SD_{8n}})|$.

Again,  $\sqrt{(2n-1)(10n-1)} > 0$, $3(2n-1)+1 > 0$ and $(3(2n-1)+1)^2-\left(\sqrt{(2n-1)(10n-1)}\right)^2=4n(4n-3)+3>0$
and so $3(2n-1)+1 - \sqrt{(2n-1)(10n-1)}> 0$. Therefore, by \eqref{E38}, $E(K_{|v(\Gamma_{SD_{8n}})|}) > E(\Gamma_{SD_{8n}})$. 

\vspace{0.5cm}
\noindent (c) \textbf{Case 1:} $n$ is odd 

Using Theorem \ref{SD_8n1}, for $n=3$ we have $LE^+(\Gamma(SD_{24}))= 36+4\sqrt{33}$ and $LE^+(K_{|v(\Gamma(SD_{24}))|})=38$. Also, for $n \geq 5$ we have $$
LE^+(\Gamma_{SD_{8n}}) - LE^+(K_{|v(\Gamma_{SD_{8n}})|})= \frac{2(8n^2(2n-9)+66n-11)}{2n-1}+4\sqrt{8n^2-14n+9}>0.
$$
  Therefore, $LE^+(\Gamma_{SD_{8n}})>LE^+(K_{|v(\Gamma_{SD_{8n}})|})$ which implies $\Gamma_{SD_{8n}}$ is Q-hyperenergetic and consequently part (a) implies $\Gamma_{SD_{8n}}$ is L-hyperenergetic.

\vspace{0.5cm}
\noindent \textbf{Case 2:} $n$ is even

For $n=2$ we have $LE^+(K_{|v(\Gamma_{SD_{16}})|})= 16$ and using Theorem \ref{SD_8n2}, $LE^+(\Gamma_{SD_{16}})=\frac{134}{7}+2\sqrt{73}$. Therefore, $\Gamma_{SD_{16}}$ is Q-hyperenergetic and consequently part (a) implies $\Gamma_{SD_{16}}$ is L-hyperenergetic.

For $n \geq 4$, using Theorem \ref{SD_8n2}, we have $$
LE^+(\Gamma_{SD_{8n}}) - LE^+(K_{|v(\Gamma_{SD_{8n}})|})= \frac{64n^2(n-3)+144n-22}{4n-1}+2\sqrt{32n^2-32n+9}>0.
$$
Therefore, $LE^+(\Gamma_{SD_{8n}})>LE^+(K_{|v(\Gamma_{SD_{8n}})|})$ which implies $\Gamma_{SD_{8n}}$ is Q-hyperenergetic and consequently part (a) implies $\Gamma_{SD_{8n}}$ is L-hyperenergetic.
\end{proof}

In Theorem \ref{SD_{8n}}, we compare $E(\Gamma_{SD_{8n}})$, $LE(\Gamma_{SD_{8n}})$ and $LE^+(\Gamma_{SD_{8n}})$. However, in the following figures, we show how close are they.

\vspace{.3cm}

\begin{minipage}[t]{.5\linewidth}
\begin{tikzpicture}
\begin{axis}
[
xlabel={$n$ $\rightarrow$},
ylabel={Energies of $\Gamma_{SD_{8n}}$ $\rightarrow$},
xmin=2, xmax=11,
ymin=0, ymax=1800,
grid = both,
minor tick num = 1,
major grid style = {lightgray},
minor grid style = {lightgray!25},
width=.7\textwidth,
height=.7\textwidth,
legend style={legend pos=north west},
 ]
\addplot[domain=5:22,samples at={3,5,7,...,21,23},mark=*,green, samples=9, mark size=.8pt]{2*(x-1)+4*(5*x*x-6*x+1)^(1/2)};
\tiny
\addlegendentry{$E$}
\addplot[domain=5:22,samples at={3,5,7,...,21,23},mark=triangle*,blue,mark size=.8pt, samples=9]{8*x*(4*x*x-10*x+7)/(2*x-1)};
\tiny
\addlegendentry{$LE$}
\addplot[domain=5:22,samples at={5,7,...,21,23},mark=square*, red, mark size=.8pt, samples=9]{(32*x*x*x-112*x*x+96*x-12)/(2*x-1)+4*(8*x*x-16*x+9)^(1/2)};
\tiny
\addlegendentry{$LE^+$}
\addplot[domain=3:4,samples at={3},mark=square*, red, mark size=.8pt, samples=9]{36+4*(33)^(1/2)};
\tiny
\end{axis}
\end{tikzpicture}
\vspace{-.2 cm}
\captionsetup{font=footnotesize}

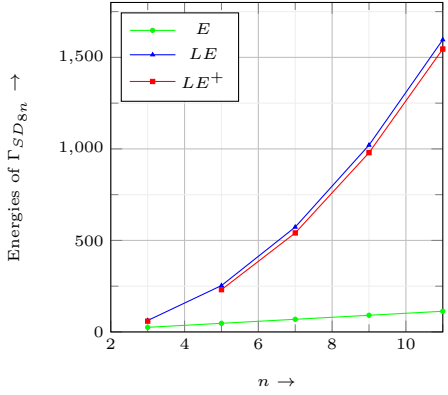
\captionof{figure}{Energies of $\Gamma_{SD_{8n}}$, where $n$ is odd}
\end{minipage}
\hspace{0.05cm}
\begin{minipage}[t]{.5\linewidth}
\begin{tikzpicture}
\begin{axis}
[
xlabel={$n$ $\rightarrow$},
ylabel={Energies of $\Gamma_{SD_{8n}}$ $\rightarrow$},
xmin=1, xmax=12,
ymin=0, ymax=2500,
grid = both,
minor tick num = 1,
major grid style = {lightgray},
minor grid style = {lightgray!25},
width=.7\textwidth,
height=.7\textwidth,
legend style={legend pos=north west},
 ]
\addplot[domain=4:22,samples at={2,4,6,8,10,12,14,16,18,20,22},mark=*,green, samples=10, mark size=.8pt]{(4*x-2)+2*(20*x*x-12*x+1)^(1/2)};
\tiny
\addlegendentry{$E$}
\addplot[domain=4:22,samples at={2,4,6,8,10,12,14,16,18,20,22},mark=triangle*,blue,mark size=.8pt, samples=10]{(64*x*x*x-64*x*x+24*x)/(4*x-1)};
\tiny
\addlegendentry{$LE$}
\addplot[domain=4:22,samples at={4,6,8,10,12,14,16,18,20,22},mark=square*, red, mark size=.8pt, samples=10]{(64*x*x*x-128*x*x+64*x-6)/(4*x-1)+2*(32*x*x-32*x+9)^(1/2)};
\tiny
\addlegendentry{$LE^+$}
\addplot[domain=2:3,samples at={2},mark=square*, red, mark size=.8pt, samples=10]{(134)/(7)+2*(73)^(1/2)};
\tiny
\end{axis}
\end{tikzpicture}
\vspace{-.2 cm}
\captionsetup{font=footnotesize}

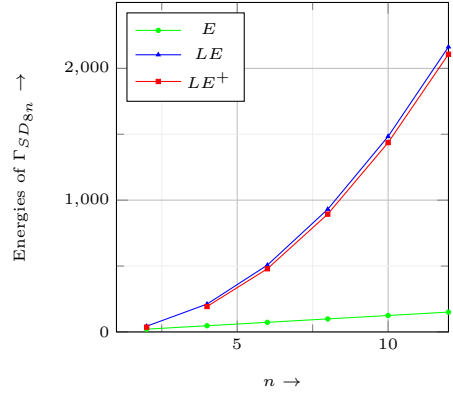
\captionof{figure}{Energies of $\Gamma_{SD_{8n}}$, where $n$ is even}
\end{minipage}

\subsection{The groups, $V_{8n}$}
We consider $V_{8n}:=\langle a,b : a^{4n} = b^4 = 1, b^{-1}ab^{-1}= bab = a^{-1} \rangle$, the groups of order $8n$ (where $n > 1$). We compute Signless Laplacian spectrum, Signless Laplacian energy, Laplacian spectrum, Laplacian energy and spectrum and energy of $\Gamma_{V_{8n}}$. If $n$ is odd then energy and Laplacian energy of $\Gamma_{V_{8n}}$ are as given in the following theorem.
\begin{theorem}[\protect{\cite[Theorem 4.1.1 (a) and Theorem 4.3.1]{FWNT21}}]\label{V8n}
Let $G$ be isomorphic to $V_{8n}$, where  $n$ is odd. Then
\[ E(\Gamma_{V_{8n}}) = 2(2n-1)+2\sqrt{(2n-1)(10n-1)} \, \text{ and } \, 
LE(\Gamma_{V_{8n}}) = \frac{8n(8n^2-8n+3)}{4n-1}.\]
\end{theorem}

\begin{theorem}\label{V8n1}
Let $G$ be isomorphic to $V_{8n}$, where $n$ is odd. Then
\begin{align*}
\Q-spec(\Gamma_{V_{8n}})=&\left\lbrace(8n-4)^{2n},(4n)^{4n-3},(8n-6)^{2n-1}, \left(8n-3+\sqrt{32n^2-32n+9}\right)^1,\right.\\
&\left.\left(8n-3-\sqrt{32n^2-32n+9}\right)^1\right\rbrace 
\end{align*} and $LE^+(\Gamma_{V_{8n}}) = \frac{64n^3-128n^2+64n-6}{4n-1}+2\sqrt{32n^2-32n+9}$. 
\end{theorem}

\begin{proof}
If $G\cong V_{8n}$ and $n$ is odd then $|v(\Gamma_{V_{8n}})|=8n-2$ and $\Gamma_{V_{8n}}=K_{2n.2,1.(4n-2)}$. Using Theorem \ref{R1}(b), we have 
\begin{align*}
Q_{\Gamma_{V_{8n}}}(x)=&{\displaystyle \prod_{i=1}^2} (x-(8n-2)+p_i)^{a_i(p_i-1)}{\displaystyle \prod_{i=1}^2}(x-(8n-2)+2p_i)^{a_i}\left(1-{\displaystyle \sum_{i=1}^2}\frac{a_ip_i}{x-(8n-2)+2p_i}\right)\\
  = \,& (x-8n+4)^{2n}(x-4n)^{4n-3}(x-8n+6)^{2n}(x-2)\left(1-\frac{4n}{x-8n+6}-\frac{4n-2}{x-2}\right)\\
   = \,& (x-(8n-4))^{2n}(x-4n)^{4n-3}(x-(8n-6))^{2n-1}(x^2-(16n-6)x+32n^2-16n).
\end{align*}
 Thus, $\Q-spec(\Gamma_{V_{8n}})=\left\lbrace(8n-4)^{2n},(4n)^{4n-3},(8n-6)^{2n-1}, \left(8n-3+\sqrt{32n^2-32n+9}\right)^1,\right.\\
 ~~~~~~~~~~~~~~~~~~~~~~~~~~~~~~~~~\left.\left(8n-3-\sqrt{32n^2-32n+9}\right)^1\right\rbrace$.

 Number of edges of $\Gamma_{V_{8n}}^c$ is $8n^2-8n+3$. Therefore, $|e(\Gamma_{V_{8n}})| = \frac{(8n-2)(8n-2-1)}{2}-(8n^2-8n+3)= $ $12n(2n-1)$. Now
\[\left|8n-4 - \frac{2|e(\Gamma_{V_{8n}})|}{|v(\Gamma_{V_{8n}})|}\right| = \left|\frac{(8n-4)(n-1)}{4n-1}\right| = \frac{(8n-4)(n-1)}{4n-1},\]
\[\left|4n - \frac{2|e(\Gamma_{V_{8n}})|}{|v(\Gamma_{V_{8n}})|}\right| =\left|\frac{-8n(n-1)}{4n-1}\right| = \frac{8n(n-1)}{4n-1},\]
\[\left|8n-6 - \frac{2|e(\Gamma_{V_{8n}})|}{|v(\Gamma_{V_{8n}})|}\right|= \left|\frac{(8n^2-20n+6)}{4n-1}\right| = \frac{8n^2-20n+6}{4n-1},\]
\begin{align*}
 \left|8n-3+\sqrt{32n^2-32n+9} - \frac{2|e(\Gamma_{SD_{8n}})|}{|v(\Gamma_{SD_{8n}})|}\right| =&\left|\sqrt{32n^2-32n+9}+2n-\frac{3}{2}+\frac{3}{8n-2}\right|\\= & \sqrt{32n^2-32n+9}+2n-\frac{3}{2}+\frac{3}{8n-2} 
\end{align*}
 and
\begin{align*}
\left|8n-3-\sqrt{32n^2-32n+9} - \frac{2|e(\Gamma_{SD_{8n}})|}{|v(\Gamma_{SD_{8n}})|}\right|  =&\left|-\sqrt{32n^2-32n+9}+2n-\frac{3}{2}+\frac{3}{8n-2}\right|\\=& \sqrt{32n^2-32n+9}-2n+\frac{3}{2}-\frac{3}{8n-2}.\end{align*}


Therefore, for $n \geq 3$, we have 
\begin{align*}
LE^+(\Gamma_{V_{8n}}) =& 2n\times \frac{(8n-4)(n-1)}{4n-1}+(4n-3)\times \frac{8n(n-1)}{4n-1}+(2n-1)\times\frac{8n^2-20n+6}{4n-1}+\\ &\sqrt{32n^2-32n+9}+2n-\frac{3}{2}+\frac{3}{8n-2}+\sqrt{32n^2-32n+9}-2n+\frac{3}{2}-\frac{3}{8n-2}    
\end{align*}
and the result follows on simplification.
\end{proof}

\begin{theorem}\label{V8n2}
Let $G$ be isomorphic to $V_{8n}$, where $n$ is even. Then 
 \begin{align*}
 \Q-spec&(\Gamma_{V_{8n}})\\ =&\left\lbrace(8n-8)^{3n},(4n)^{4n-5},(8n-12)^{n-1}, \left(8n-6+2\sqrt{8n^2-16n+9}\right)^1, \left(8n-6-2\sqrt{8n^2-16n+9}\right)^1\right\rbrace,
 \end{align*}
 \[\L-spec(\Gamma_{V_{8n}})=\{0,(4n)^{4n-5},(8n-8)^{3n},(8n-4)^{n}\} \text{ and}\]
 \[
 \spec(\Gamma_{V_{8n}})=\left\lbrace 0^{7n-5}, (-4)^{n-1}, \left(2(n-1)+2\sqrt{(n-1)(5n-1)}\right)^1, \left(2(n-1)+2\sqrt{(n-1)(5n-1)}\right)^1\right\rbrace. 
 \]
 Further
\begin{align*}
LE^+(\Gamma_{V_{8n}}) = \begin{cases}
\frac{24n^3-64n^2+32n+12}{2n-1}+4\sqrt{8n^2-16n+9}, & \mbox{if $n \leq 4$}\vspace{.2cm}\\
\frac{32n^3-112n^2+96n-12}{2n-1}+4\sqrt{8n^2-16n+9}, & \mbox{if $n \geq 6$,}\end{cases}  
\end{align*}
\[LE(\Gamma_{V_{8n}})=\frac{8n(4n^2-10n+7)}{2n-1} \, \text{ and } \,
E(\Gamma_{V_{8n}})=4(n-1)+4\sqrt{(n-1)(5n-1)}.\]
\end{theorem}
\begin{proof}
If $G\cong V_{8n}$ and $n$ is even then $|v(\Gamma_{V_{8n}})|=8n-4$ and $\Gamma_{V_{8n}}=K_{n.4,1.(4n-4)}$. Using Theorem \ref{R1}(b), we have 
\begin{align*}
Q_{\Gamma_{V_{8n}}}(x)=&{\displaystyle \prod_{i=1}^2} (x-(8n-4)+p_i)^{a_i(p_i-1)}{\displaystyle \prod_{i=1}^2}(x-(8n-4)+2p_i)^{a_i}\left(1-{\displaystyle \sum_{i=1}^2}\frac{a_ip_i}{x-(8n-4)+2p_i}\right)\\
  = \,& (x-8n+8)^{3n}(x-4n)^{4n-5}(x-8n+12)^n(x-4)\left(1-\frac{4n}{x-8n+12}-\frac{4n-4}{x-4}\right)\\
   = \,& (x-(8n-8))^{3n}(x-4n)^{4n-5}(x-(8n-12))^{n-1}(x^2-(16n-12)x+32n^2-32n).
\end{align*}
Thus, $\Q-spec(\Gamma_{V_{8n}})=\left\lbrace (8n-8)^{3n}, (4n)^{4n-5}, (8n-12)^{n-1}, \left(8n-6+2\sqrt{8n^2-16n+9}\right)^1,\right.\\
 ~~~~~~~~~~~~~~~~~~~~~~~~~~~~~~~~~\left.\left(8n-6-2\sqrt{8n^2-16n+9}\right)^1\right\rbrace$.

 Number of edges of $\Gamma_{V_{8n}}^c$ is $8n^2-12n+10$ and so $|e(\Gamma_{V_{8n}})| = \frac{(8n-4)(8n-4-1)}{2}-(8n^2-12n+10)=24n(n-1)$. Now
\[
\left|8n-8 - \frac{2|e(\Gamma_{V_{8n}})|}{|v(\Gamma_{V_{8n}})|}\right| = \left|\frac{4(n-1)(n-2)}{2n-1}\right| = \frac{4(n-1)(n-2)}{2n-1},
\] 
\[
\left|4n - \frac{2|e(\Gamma_{V_{8n}})|}{|v(\Gamma_{V_{8n}})|}\right| = \left|\frac{4n(2-n)}{2n-1}\right| = \frac{4n(n-2)}{2n-1},
\]
\[
\left|8n-12 - \frac{2|e(\Gamma_{V_{8n}})|}{|v(\Gamma_{V_{8n}})|}\right| = \left|\frac{4(n^2-5n+3)}{2n-1}\right| = \begin{cases}  \frac{4(-n^2+5n-3)}{2n-1}, & \mbox{if $n \leq 4$}\vspace{.2cm}\\
\frac{4(n^2-5n+3)}{2n-1}, & \mbox{if $n \geq 6$,}\end{cases}
\]
\begin{align*}
\left|8n-6+2\sqrt{8n^2-16n+9} - \frac{2|e(\Gamma_{V_{8n}})|}{|v(\Gamma_{V_{8n}})|}\right| =& \left|2\sqrt{8n^2-16n+9}+(2n-3)+\frac{3}{2n-1}\right|\\ =& 2\sqrt{8n^2-16n+9}+(2n-3)+\frac{3}{2n-1}    
\end{align*} 
and
\begin{align*}
\left|8n-6-2\sqrt{8n^2-16n+9} - \frac{2|e(\Gamma_{V_{8n}})|}{|v(\Gamma_{V_{8n}})|}\right| = &\left|-2\sqrt{8n^2-16n+9}+(2n-3)+\frac{3}{2n-1}\right|\\ =& 2\sqrt{8n^2-16n+9}-(2n-3)-\frac{3}{2n-1}.    
\end{align*}

 Therefore, for $n \leq 4$, we have 
\begin{align*}
LE^+(\Gamma_{V_{8n}})=&3n\times \frac{4(n-1)(n-2)}{2n-1}+(4n-5)\times \frac{4n(n-2)}{2n-1}+(n-1)\times\frac{4(-n^2+5n-3)}{2n-1}+\\ &2\sqrt{8n^2-16n+9}-2n-3+\frac{3}{2n-1}+2\sqrt{8n^2-16n+9}+2n+3-\frac{3}{2n-1}.   
\end{align*}
For $n \geq 6$, we have
\begin{align*}
LE^+(\Gamma_{V_{8n}}) =& 3n\times \frac{4(n-1)(n-2)}{2n-1}+(4n-5)\times \frac{4n(n-2)}{2n-1}+(n-1)\times\frac{4(n^2-5n+3)}{2n-1}+\\
&2\sqrt{8n^2-16n+9}-2n-3+\frac{3}{2n-1}+2\sqrt{8n^2-16n+9}+2n+3-\frac{3}{2n-1}.    
\end{align*}
Thus we get the required expressions for $LE^+(\Gamma_{V_{8n}})$ on simplification.

Since $\Gamma_{V_{8n}}^c=nK_4 \cup K_{4n-4}$,  using Theorem \ref{L1}, we have 
\begin{align*} \L-spec(\Gamma_{V_{8n}})=&\left\{(0)^1,\left({\displaystyle \sum_{i=1}^2}l_im_i-m_2\right)^{l_2(m_2-1)},\left({\displaystyle \sum_{i=1}^2}l_im_i-m_1\right)^{l_1(m_1-1)},\left({\displaystyle \sum_{i=1}^2}l_im_i\right)^{{\displaystyle \sum_{i=1}^2}l_i-1}\right\}\\
  = \,& \{(0)^1,(n.4+4n-4-4n+4)^{1(4n-4-1)},(n.4+4n-4-4)^{n(4-1)},(n.4+4n-4)^{n+1-1}\}\\
   = \,& \{(0)^1,(4n)^{4n-5},(8n-8)^{3n},(8n-4)^{n}\}.
\end{align*}

Now \[\left|0 - \frac{2|e(\Gamma_{V_{8n}})|}{|v(\Gamma_{V_{8n}})|}\right| = \left|\frac{-12n(n-1)}{2n-1}\right| = \frac{12n(n-1)}{2n-1}, \quad
\left|4n - \frac{2|e(\Gamma_{V_{8n}})|}{|v(\Gamma_{V_{8n}})|}\right| = \left|\frac{4n(2-n)}{2n-1}\right| = \frac{4n(n-2)}{2n-1},
\]
\[
\left|8n-8 - \frac{2|e(\Gamma_{V_{8n}})|}{|v(\Gamma_{V_{8n}})|}\right| = \left|\frac{4(n-1)(n-2)}{2n-1}\right| = \frac{4(n-1)(n-2)}{2n-1}
\] 
and 
\[
\left|8n-4 - \frac{2|e(\Gamma_{V_{8n}})|}{|v(\Gamma_{V_{8n}})|}\right| = \left|\frac{4(n^2-n+1)}{2n-1}\right| = \frac{4(n^2-n+1)}{2n-1}.
\]
 Therefore, $LE(\Gamma_{V_{8n}}) = 1\times \frac{12n(n-1)}{2n-1}+(4n-5)\times \frac{4n(n-2)}{2n-1}+3n \times \frac{4(n-1)(n-2)}{2n-1}+n \times \frac{4(n^2-n+1)}{2n-1}$ and we get the required expression for $LE(\Gamma_{V_{8n}})$ on simplification.

Since $\Gamma_{V_{8n}}$ is a complete $(n+1)$-partite graph with $8n-4$ vertices and $\Gamma_{V_{8n}}=K_{n.4,1.(4n-4)}$. Therefore, using Theorem \ref{R1}(a), the characteristic polynomial of $\Gamma_{V_{8n}}$ is
\begin{align*}
P_{\Gamma_{V_{8n}}}(x)&=  x^{(8n-4)-(n+1)}(x+4)^{n-1}(x+4n-4)^{1-1}(x^2-(4n-4)x-16n^2-16n)\\
&= x^{7n-5}(x+4)^{n-1}(x^2-(4n-4)x-16n^2-16n).
\end{align*}
Thus, $\spec(\Gamma_{V_{8n}})=\left\lbrace (0)^{7n-5}, (-4)^{n-1}, \left(2(n-1)+2\sqrt{(n-1)(5n-1)}\right)^1, \left(2(n-1)+2\sqrt{(n-1)(5n-1)}\right)^1\right\rbrace$.
Therefore,
 $E(\Gamma_{V_{8n}})=(7n-5)\times |0|+(n-1)\times|-4|+|2(n-1)+2\sqrt{(n-1)(5n-1)}|+|2(n-1)-2\sqrt{(n-1)(5n-1)}|$
and we get the required expression for $E(\Gamma_{V_{8n}})$ on simplification.
\end{proof}

\begin{theorem}\label{V_{8n}}
If $G$ is isomorphic to $V_{8n}$ then 
\begin{enumerate}
\item $E(\Gamma_{V_{8n}}) \leq  LE^+(\Gamma_{V_{8n}}) \leq LE(\Gamma_{V_{8n}})$; equality holds if and only if $G \cong V_{16}$.
\item $\Gamma_{V_{8n}}$ is non-hypoenergetic as well as non-hyperenergetic.
\item $\Gamma_{V_{16}}$ is not L-hyperenergetic and not Q-hyperenergetic. If $n \neq 2$ then $\Gamma_{V_{8n}}$ is Q-hyperenergetic and L-hyperenergetic. 
\end{enumerate}
\end{theorem}
\begin{proof}
(a)\textbf{Case 1:} $n$ is odd

Using Theorems \ref{V8n} and \ref{V8n1}, we have 
\begin{equation}\label{V1}
LE(\Gamma_{V_{8n}})-LE^+(\Gamma_{V_{8n}})=\frac{64n^2-40n+6}{4n-1}-2\sqrt{32n^2-32n+9}
\end{equation}
and
\begin{equation}\label{V6}
LE^+(\Gamma_{V_{8n}})-E(\Gamma_{V_{8n}})=\frac{16n^2(4n-9)+76n-8}{4n-1}+2\sqrt{32n^2-32n+9}-2\sqrt{20n^2-12n+1}.
\end{equation}
Since $64n^2-40n+6 > 0$, $2(4n-1)\sqrt{32n^2-32n+9} > 0$ and $(64n^2-40n+6)^2-(2\sqrt{32n^2-32n+9})^2(4n-1)^2=2048n^3(n-1)+64n(10n-1)>0$
we have $64n^2-40n+6-2(4n-1)\sqrt{32n^2-32n+9}> 0$. Therefore, by \eqref{V1}, $(4n-1)(LE^+(\Gamma_{V_{8n}})-LE(\Gamma_{V_{8n}})) > 0$. Hence, $LE(\Gamma_{V_{8n}}) > LE^+(\Gamma_{V_{8n}})$.  

Again, $\sqrt{32n^2-32n+9} > 0, \sqrt{20n^2-12n+1} > 0$ and $\left(\sqrt{32n^2-32n+9}\right)^2-\left(\sqrt{20n^2-12n+1}\right)^2=4n(3n-5)+8>0$. Therefore, 
$\sqrt{32n^2-32n+9}-\sqrt{20n^2-12n+1}>0$. Since $16n^2(4n-9)+76n-8 > 0$ we have $\frac{16n^2(4n-9)+76n-8}{2n-1}+2\sqrt{32n^2-32n+9}-2\sqrt{20n^2-12n+1}>0$. Therefore, by \eqref{V6}, $LE^+(\Gamma_{V_{8n}}) >  E(\Gamma_{V_{8n}})$. Hence, 
$E(\Gamma_{V_{8n}}) < LE^+(\Gamma_{V_{8n}}) < LE(\Gamma_{V_{8n}})$.



\vspace{0.5cm}
\noindent \textbf{Case 2:} $n$ is even 

Using Theorems \ref{V8n} and \ref{V8n2}, for $n \leq 4$, we have 
\begin{equation}\label{V2}
LE(\Gamma_{V_{8n}})-LE^+(\Gamma_{V_{8n}})=\frac{8n^3-16n^2+24n-12}{2n-1}-4\sqrt{8n^2-16n+9}
\end{equation}
and
\begin{equation}\label{V3}
LE^+(\Gamma_{V_{8n}})-E(\Gamma_{V_{8n}})=\frac{4(n-2)(6n^2-6n-1)}{2n-1}+4\sqrt{8n^2-16n+9}-4\sqrt{5n^2-6n+1}.
\end{equation}
Since $8n^3-16n^2+24n-12> 0$, $4(2n-1)\sqrt{8n^2-16n+9} > 0$ and
\begin{center}
 $(8n^3-16n^2+24n-12)^2-\left(4\sqrt{8n^2-16n+9}\right)^2(2n-1)^2=64n(n-2)^2(n-1)(n^2+n-1) \geq 0$
\end{center}
(equality holds if and only if  $n = 2$)
we have $8n^3-16n^2+24n-12 - 4(2n-1)\sqrt{8n^2-16n+9} \geq 0$. Therefore, by \eqref{V2}, $(2n-1)(LE(\Gamma_{V_{8n}})-LE^+(\Gamma_{V_{8n}})) \geq 0$. Hence, $LE(\Gamma_{V_{8n}}) \geq LE^+(\Gamma_{V_{8n}})$ equality holds if and only if $G \cong V_{16}$.  

Again,  $\sqrt{8n^2-16n+9} > 0, \sqrt{5n^2-6n+1} > 0$ and $\left(\sqrt{8n^2-16n+9}\right)^2-\left(\sqrt{5n^2-6n+1}\right)^2=n(3n-10)+8 \geq 0$ (equality holds if and only if $n = 2$). Therefore, 
$\sqrt{8n^2-16n+9}-\sqrt{5n^2-6n+1} \geq 0$. Since $4(n-2)(6n^2-6n-1) \geq 0$ we have $\frac{4(n-2)(6n^2-6n-1)}{2n-1}+4\sqrt{8n^2-16n+9}-4\sqrt{5n^2-6n+1} \geq 0$. Therefore, by \eqref{V3}, $LE^+(\Gamma_{V_{8n}}) \geq  E(\Gamma_{V_{8n}})$. Hence, 
$E(\Gamma_{V_{8n}}) \leq LE^+(\Gamma_{V_{8n}}) \leq LE(\Gamma_{V_{8n}})$ equality holds if and only if $G \cong V_{16}$.


Using Theorems \ref{V8n} and \ref{V8n2}, for $n \geq 6$, we have 
\begin{equation}\label{V4}
LE(\Gamma_{V_{8n}})-LE^+(\Gamma_{V_{8n}})=\frac{32n^2-40n+12}{2n-1}-4\sqrt{8n^2-16n+9}
\end{equation}
and
\begin{equation}\label{V5}
LE^+(\Gamma_{V_{8n}})-E(\Gamma_{V_{8n}})=\frac{32n^3-116n^2+102n-14}{2n-1}+4\sqrt{8n^2-16n+9}-4\sqrt{5n^2-6n+1}.
\end{equation}
Since $32n^2-40n+12 > 0$, $4(2n-1)\sqrt{8n^2-16n+9} > 0$ and $(32n^2-40n+12)^2-\left(4\sqrt{8n^2-16n+9}\right)^2(2n-1)^2=512n^3(n-2)+128n(5n-1)>0$
we have $32n^2-40n+12 - 4(2n-1)\sqrt{8n^2-16n+9} > 0$. Therefore, by \eqref{V4}, $(2n-1)(LE(\Gamma_{V_{8n}})-LE^+(\Gamma_{V_{8n}})) > 0$. Hence, $LE(\Gamma_{V_{8n}}) > LE^+(\Gamma_{V_{8n}})$.  

Again, $\sqrt{8n^2-16n+9} > 0, \sqrt{5n^2-6n+1} > 0$ and $\left(\sqrt{8n^2-16n+9}\right)^2-\left(\sqrt{5n^2-6n+1}\right)^2=n(3n-10)+8>0$. Therefore, 
$\sqrt{8n^2-16n+9}-\sqrt{5n^2-6n+1}>0$. Since $32n^3-116n^2+102n-14 > 0$ we have $\frac{32n^3-116n^2+102n-14}{2n-1}+4\sqrt{8n^2-16n+9}-4\sqrt{5n^2-6n+1}>0$. Therefore, by \eqref{V5}, $LE^+(\Gamma_{V_{8n}}) >  E(\Gamma_{V_{8n}})$. Hence, 
$E(\Gamma_{V_{8n}}) < LE^+(\Gamma_{V_{8n}}) < LE(\Gamma_{V_{8n}})$.

\vspace{0.5cm}

\noindent (b) \textbf{Case 1:} $n$ is odd

Here $|v(\Gamma_{V_{8n}})|= 8n-2$ and $E(K_{|v(\Gamma_{V_{8n}})|})= LE(K_{|v(\Gamma_{V_{8n}})|})= LE^+(K_{|v(\Gamma_{V_{8n}})|})=16n-16$.
Using Theorem \ref{V8n}, we have
\begin{equation}\label{E39}
E(\Gamma_{V_{8n}}) - |v(\Gamma_{V_{8n}})|= 2\left(\sqrt{(2n-1)(10n-1)}-2n\right)
\end{equation}
and
\begin{equation}\label{E40}
E(K_{|v(\Gamma_{V_{8n}})|})-E(\Gamma_{V_{8n}})= 2\left(3(2n-1)+1-\sqrt{(2n-1)(10n-1)}\right).
\end{equation}
Since $\sqrt{(2n-1)(10n-1)} > 0$, $2n > 0$ and $\left(\sqrt{(2n-1)(10n-1)}\right)^2-(2n)^2=4n(4n-3)+1>0$
we have $ \sqrt{(2n-1)(10n-1)}-2n> 0$. Therefore, by \eqref{E39}, $E(\Gamma_{V_{8n}}) > |v(\Gamma_{V_{8n}})|$.

Again,  $\sqrt{(2n-1)(10n-1)} > 0$, $3(2n-1)+1 > 0$ and $(3(2n-1)+1)^2-\left(\sqrt{(2n-1)(10n-1)}\right)^2=4n(4n-3)+3>0$
and so $3(2n-1)+1 - \sqrt{(2n-1)(10n-1)}> 0$. Therefore, by \eqref{E40}, $E(K_{|v(\Gamma_{V_{8n}})|}) > E(\Gamma_{V_{8n}})$.


\vspace{0.2cm}
\noindent \textbf{Case 2:} $n$ is even 

Here $|v(\Gamma_{V_{8n}})|= 8n-4$ and $E(K_{|v(\Gamma_{V_{8n}})|})= LE(K_{|v(\Gamma_{V_{8n}})|})= LE^+(K_{|v(\Gamma_{V_{8n}})|})=16n-10$.
Using Theorem \ref{V8n}, we have
\begin{equation}\label{E41}
E(\Gamma_{V_{8n}})-|v(\Gamma_{V_{8n}})|= 4\sqrt{(n-1)(5n-1)}-4n    
\end{equation}
and
\begin{equation}\label{E42}
E(K_{|v(\Gamma_{V_{8n}})|})-E(\Gamma_{V_{8n}})= 12n-6-4\sqrt{(n-1)(5n-1)}.    
\end{equation}
Since $4\sqrt{(n-1)(5n-1)} > 0$, $4n > 0$ and $\left(4\sqrt{(n-1)(5n-1)}\right)^2-(4n)^2=16(4n^2-6n+1)>0$
we have $ 4\sqrt{(n-1)(5n-1)}-4n > 0$. Therefore, by \eqref{E41}, $E(\Gamma_{V_{8n}}) > |v(\Gamma_{V_{8n}})|$.

Again,  $4\sqrt{(n-1)(5n-1)} > 0$, $12n-6 > 0$ and $(12n-6)^2-\left(4\sqrt{(n-1)(5n-1)}\right)^2=4(16n^2-12n+5)>0$
and so $12n-6-4\sqrt{(n-1)(5n-1)}> 0$. Therefore, by \eqref{E42}, $E(K_{|v(\Gamma_{V_{8n}})|}) > E(\Gamma_{V_{8n}})$.


\vspace{0.5cm}
\noindent (c) \textbf{Case 1:} $n$ is odd

Using Theorem \ref{V8n1} we have
$LE^+(\Gamma_{V_{8n}}) - LE^+(K_{|v(\Gamma_{V_{8n}})|})= \frac{64n^2(n-3)+144n-22}{4n-1}+2\sqrt{32n^2-32n+9}>0$. Therefore, $LE^+(\Gamma_{V_{8n}})>LE^+(K_{|v(\Gamma_{V_{8n}})|})$ which implies $\Gamma_{V_{8n}}$ is Q-hyperenergetic and consequently part (a) implies $\Gamma_{V_{8n}}$ is L-hyperenergetic. 

\vspace{0.2cm}
\noindent \textbf{Case 2:} $n$ is even

Using Theorem \ref{V8n2}, for $n=2$, we have $LE(\Gamma_{V_{8n}})=16$ and $LE(K_{|v(\Gamma_{V_{8n}})|})=22$. Clearly, $LE(\Gamma_{V_{8n}}) < LE(K_{|v(\Gamma_{V_{8n}})|})$. For $n \leq 4$, 
$$
LE^+(\Gamma_{V_{8n}}) - LE^+(K_{|v(\Gamma_{V_{8n}})|})= \frac{4(6n^2(n-4)+20n-1)}{2n-1}+4\sqrt{8n^2-16n+9}>0
$$
 for all $n \neq 2$. Therefore, for all $n \neq 2$, $LE^+(\Gamma_{V_{8n}})>LE^+(K_{|v(\Gamma_{V_{8n}})|})$ which implies $\Gamma_{V_{8n}}$ is Q-hyperenergetic and consequently part (a) implies $\Gamma_{V_{8n}}$ is L-hyperenergetic. For $n \geq 6$, 
 $$
 LE^+(\Gamma_{V_{8n}}) - LE^+(K_{|v(\Gamma_{V_{8n}})|})= \frac{2(8n^2(2n-9)+66n-11)}{2n-1}+4\sqrt{8n^2-16n+9}>0.
 $$
  Hence, the result holds.
\end{proof}

 In Theorem \ref{V_{8n}}, we compare $E(\Gamma_{V_{8n}})$, $LE(\Gamma_{V_{8n}})$ and $LE^+(\Gamma_{V_{8n}})$. However, in the following figures, we show how close are they.

\vspace{.3cm}

\begin{minipage}[t]{.5\linewidth}
\begin{tikzpicture}
\begin{axis}
[
xlabel={$n$ $\rightarrow$},
ylabel={Energies of $\Gamma_{V_{8n}}$ $\rightarrow$},
xmin=2, xmax=11,
ymin=0, ymax=2000,
grid = both,
minor tick num = 1,
major grid style = {lightgray},
minor grid style = {lightgray!25},
width=.7\textwidth,
height=.7\textwidth,
legend style={legend pos=north west},
 ]
\addplot[domain=3:22,samples at={3,5,7,...,21,23},mark=*,green, samples=9, mark size=.8pt]{(4*x-2)+2*(20*x*x-12*x+1)^(1/2)};
\tiny
\addlegendentry{$E$}
\addplot[domain=3:22,samples at={3,5,7,...,21,23},mark=triangle*,blue,mark size=.8pt, samples=9]{(64*x*x*x-64*x*x+24*x)/(4*x-1)};
\tiny
\addlegendentry{$LE$}
\addplot[domain=3:22,samples at={3,5,7,...,21,23},mark=square*, red, mark size=.8pt, samples=9]{(64*x*x*x-128*x*x+64*x-6)/(4*x-1)+2*(32*x*x-32*x+9)^(1/2)};
\tiny
\addlegendentry{$LE^+$}
\end{axis}
\end{tikzpicture}
\vspace{-.2 cm}
\captionsetup{font=footnotesize}

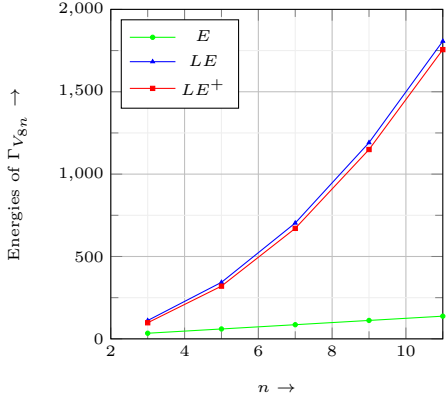
\captionof{figure}{Energies of $\Gamma_{V_{8n}}$, where $n$ is odd}
\end{minipage}
\hspace{0.05cm}
\begin{minipage}[t]{.5\linewidth}
\begin{tikzpicture}
\begin{axis}
[
xlabel={$n$ $\rightarrow$},
ylabel={Energies of $\Gamma_{V_{8n}}$  $\rightarrow$},
xmin=1, xmax=12,
ymin=0, ymax=2000,
grid = both,
minor tick num = 1,
major grid style = {lightgray},
minor grid style = {lightgray!25},
width=.7\textwidth,
height=.7\textwidth,
legend style={legend pos=north west},
 ]
\addplot[domain=4:22,samples at={2,4,6,8,10,12,14,16,18,20,22},mark=*,green, samples=10, mark size=.8pt]{2*(x-1)+4*(5*x*x-6*x+1)^(1/2)};
\tiny
\addlegendentry{$E$}
\addplot[domain=4:22,samples at={2,4,6,8,10,12,14,16,18,20,22},mark=triangle*,blue,mark size=.8pt, samples=10]{8*x*(4*x*x-10*x+7)/(2*x-1)};
\tiny
\addlegendentry{$LE$}
\addplot[domain=6:22,samples at={6,8,10,12,14,16,18,20,22},mark=square*, red, mark size=.8pt, samples=10]{(32*x*x*x-112*x*x+96*x-12)/(2*x-1)+4*(8*x*x-16*x+9)^(1/2)};
\tiny
\addlegendentry{$LE^+$}
\addplot[domain=2:4,samples at={2,4},mark=square*, red, mark size=.8pt, samples=10]{(24*x*x*x-64*x*x+32*x+12)/(2*x-1)+4*(8*x*x-16*x+9)^(1/2)};
\tiny
\end{axis}
\end{tikzpicture}
\vspace{-.2 cm}
\captionsetup{font=footnotesize}

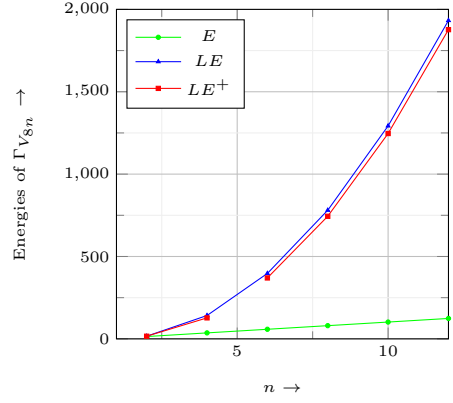
\captionof{figure}{Energies of $\Gamma_{V_{8n}}$, where $n$ is even}
\end{minipage}

\subsection{The Frobenious groups of order $pq$}
We consider $F_{p,q}: = \langle a, b : a^{p} = b^q=1; b^{-1}ab = a^{u} \rangle$, the Frobenious groups of order $pq$, where $p$ and $q$ are two primes such that $q|(p-1)$ and $u$ is an integer such that $\Bar{u} \in \mathbb{Z}_p \setminus \{\Bar{0}\}$ having order $q$. Results regarding different energies of non-commuting graph of  $F_{p,q}$ are given below.
\begin{theorem}[\protect{\cite[(4.1.e) and (4.3.b)]{FWNT21}}]\label{pq}
	Let $G$ be isomorphic to $F_{p,q}$ . Then
	\[E(\Gamma_G)= \alpha + \sqrt{\alpha^2+4p\alpha} \, \text{ and } \, 
	LE(\Gamma_{G})=\frac{2p^2\alpha+2p(q-1)^2}{pq-1}, \text{ where } \alpha=(p-1)(q-1).\] 
\end{theorem}

\begin{theorem}\label{pq1}
	Let $G$ be isomorphic to $F_{p,q}$. Then
	\begin{align*}
		\Q-spec(\Gamma_G)=&\left\{(pq-p)^{p-2},(pq-q)^{pq-2p},(pq-2q+1)^{p-1}, \left(\frac{A}{2}\right)^1, \left(\frac{B}{2}\right)^1 \right\},
	\end{align*} where $A=3pq-2p-2q+1+\sqrt{pq(pq-2)+4(p-q)(pq-p-q+1)+1}$ and $B= 3pq-2p-2q+1-\sqrt{pq(pq-2)+4(p-q)(pq-p-q+1)+1}$ and 
	\begin{align*}
		LE^+(\Gamma_G) = \frac{2p^3q-p^2q^2-2pq^2-6pq-4p^3+6p^2+2q-1}{pq-1}+\sqrt{pq(pq-2)+4(p-q)(pq-p-q+1)+1}.  
	\end{align*}
\end{theorem}
\begin{proof}
	If $G \cong F_{p,q}$  then $|v(\Gamma_G)|=pq-1$ and $\Gamma_G=K_{1.(p-1), p.(q-1)}$. Using Theorem \ref{R1}(b), we have 
	\begin{align*}
		Q_{\Gamma_G}(x)= & {\displaystyle \prod_{i=1}^2} (x-(pq-1)+p_i)^{a_i(p_i-1)}{\displaystyle \prod_{i=1}^2}(x-(pq-1)+2p_i)^{a_i}\left(1-{\displaystyle \sum_{i=1}^2}\frac{a_ip_i}{x-(pq-1)+2p_i}\right)\\
		=& (x-pq+p)^{p-2}(x-pq+q)^{pq-2p}(x-pq+2p-1)(x-pq+2q-1)^p\\
		&\times \left(1-\frac{p-1}{x-pq+2p-1}-\frac{pq-p}{x-pq+2q-1}\right)\\
		=& (x-(pq-p))^{p-2}(x-(pq-q))^{pq-2p}(x-(pq-2q+1))^p\\
		&\times (x^2-(3pq-2p-2q+1)x+2p^2q^2-4p^2q+2pq-2pq^2+2p^2-2p+2q-1).
	\end{align*}
	Thus, $\Q-spec(\Gamma_G)=\left\{(pq-p)^{p-2},(pq-q)^{pq-2p},(pq-2q+1)^{p-1}, \left(\frac{A}{2}\right)^1, \left(\frac{B}{2}\right)^1 \right\}$.
	
	Number of edges of $\Gamma_G^c$ is $\frac{p^2-3p+2+pq^2-3pq+2p}{2}$. Thus, $|e(\Gamma_G)| = \frac{(pq-1)(pq-2)}{2}-\frac{p^2-3p+2+pq^2-3pq+2p}{2}=\frac{(p^2-p)(q^1-1)}{2}$. Now,
	\[\left|pq-p - \frac{2|e(\Gamma_G)|}{|v(\Gamma_G)|}\right| = \left|\frac{(q-1)(pq-p^2)}{pq-1}\right|= \frac{(q-1)(p^2-pq)}{pq-1},\] 
	\[\left|pq-q - \frac{2|e(\Gamma_G)|}{|v(\Gamma_G)|}\right| = \left|\frac{(p-q)(p-1)}{pq-1}\right|  =\frac{(p-q)(p-1)}{pq-1},\]
	\[\left|pq-2q+1 - \frac{2|e(\Gamma_G)|}{|v(\Gamma_G)|}\right| = \left|\frac{(p^2+2q)-(pq^2+p+1)}{pq-1}\right|  =\frac{-(p^2+2q)+(pq^2+p+1)}{pq-1},\]
	\[\left|A- \frac{2|e(\Gamma_G)|}{|v(\Gamma_G)|}\right| =
	\left|\frac{p^2q^2+2p^2+p+2q-2p^2q-2pq-1}{2(pq-1)}+\frac{C}{2}\right|=\frac{p^2q^2+2p^2+p+2q-2p^2q-2pq-1}{2(pq-1)}+\frac{C}{2}
	\] 
		and 
	\[\left|B- \frac{2|e(\Gamma_G)|}{|v(\Gamma_G)|}\right| =
	\left|\frac{p^2q^2+2p^2+p+2q-2p^2q-2pq-1}{2(pq-1)}-\frac{C}{2}\right|=\frac{-(p^2q^2+2p^2+p+2q-2p^2q-2pq-1)}{2(pq-1)}+\frac{C}{2},\]
	where $C= \sqrt{pq(pq-2)+4(p-q)(pq-p-q+1)+1}$. 
	Therefore, 
	
	\begin{align*}
		LE^+(\Gamma_G) =& (p-2) \times \frac{(q-1)(p^2-pq)}{pq-1} + (pq-2p) \times  \frac{(p-q)(p-1)}{pq-1}+ (p-1) \times \frac{(pq^2+p+1)-(p^2+2q)}{pq-1} +\\ 
		& \frac{p^2q^2+2p^2+p+2q-2p^2q-2pq-1}{2(pq-1)}+\frac{C}{2}-\left(\frac{p^2q^2+2p^2+p+2q-2p^2q-2pq-1}{2(pq-1)}\right)+\frac{C}{2}     
	\end{align*} and the result follows on simplification.
\end{proof}

\section{Some implications of the preceding findings}
It was shown in \cite{DDN18} that the non-commuting graphs of the groups considered in this paper are L-integral. Also, in \cite[Chapter 4]{FWNT21}, several conditions were obtained such that the non-commuting graphs of these groups are integral. In view of Theorems \ref{U6n} and \ref{Suzuki}, it follows that $\Gamma_G$ is not Q-integral if $G \cong U_{6n}$ or  $\frac{G}{Z(G)} \cong Sz(2)$. However, $\Gamma_G$ is Q-integral if  $G \cong A(n,  \mathcal{V})$, $A(n,p)$ or $\frac{G}{Z(G)} \cong \mathbb{Z}_p \times \mathbb{Z}_p$ (follows from Theorems \ref{Z_p*Z_p}, \ref{Hanaki1} and \ref{Hanaki2}).  As a consequence of our results we also have the following theorem related to Question \ref{Q1-Q-integral}. 

\begin{theorem}
 $\Gamma_G$ is Q-integral if	
\begin{enumerate}
\item $G \cong D_{2r}, M_{2rs}$ or $SD_{8r}$,  $r$ is odd  and  $8r^2 -16r + 9$ is a perfect square. 		

\item $G \cong V_{8n}$, $n$ is even and $8n^2 -16n + 9$ is a perfect square.	

\item  $G \cong Q_{4n}$ or $\frac{G}{Z(G)} \cong D_{2n}$ and  $8n^2 -16n + 9$ is a perfect square.

\item  $G \cong D_{2r}$ or $M_{2rs}$,  $r$ is even and $2r^2 -8r + 9$ is a perfect square. 

\item $G \cong V_{8n}$, $n$ is odd and $32n^2 -32n + 9$ is a perfect square.

\item $G \cong SD_{8n}$, $n$ is even and $32n^2 -32n + 9$ is a perfect square.

\item $G \cong QD_{2^n}$ and $2^{2n-1}-2^{n+2}+9$ is a perfect square.
%
	\end{enumerate}
\end{theorem}

In the following table we give some positive integers such that $8n^2-16n + 9$, $2n^2 -8n + 9$ and $32n^2 -32n + 9$ 
are perfect squares. It may be interesting to obtain general terms of such sequences of positive integers.  
\begin{table}[h]
	\centering
	\begin{tabular}{|c|c|c|c|c|c|}
		\hline
		&  &  &  &  &   \\
		$n$ & $\sqrt{8n^2-16n + 9}$ & $n$  &$\sqrt{2n^2-8n + 9}$  &$n$   & $\sqrt{32n^2-32n + 9}$   \\
		\hline
		1 & 1 & 2 &1 & 1 & 3  \\
		\hline
		2 & 3 & 4 & 3 &18  & 99  \\
		\hline
		7 & 17 & 14 & 17 &595  &3363   \\
		\hline
		36 & 99 & 72  &99  &20196  & 114243   \\
		\hline
		205 & 577 & 410 & 577 & 686053 & 3880899   \\
		\hline
		1190 & 3363 & 2380 & 3363 &23305590 &131836323    \\
		\hline
		6931 & 19601 &13862  &19601  & &   \\
		\hline
		40392 & 114243 & 80784 &114243  &  &    \\
		\hline
		235417 & 665857 &470834  &665857  &  &    \\
		\hline
		1372106 & 3880899 &2744212  &3880899  &  &    \\
		\hline
		7997215 & 22619537 &15994430  &22619537  &  &   \\
		\hline
		46611180 & 131836323 &93222360  &131836323  &  &   \\
		\hline
		271669861 & 768398401  &543339722  &768398401  &  &    \\
		\hline
	\end{tabular}
\caption{ }
\end{table}

%

\newpage

As  consequences of our results we also have the following theorems related to Questions \ref{Q2-hyperenergetic} and \ref{Q3-energies}. 
\begin{theorem}
	Let $G$ be a finite non-abelian group. Then








\begin{enumerate}
\item  $E(\Gamma_G)= LE(\Gamma_G)= LE^+(\Gamma_G)$ if $G \cong D_8, Q_8, M_{8s}, A(n,  \mathcal{V}), A(n, p)$, $V_{16}$ a non-abelian group of order $p^3$.	Also, if $\frac{G}{Z(G)} \cong \mathbb{Z}_p \times \mathbb{Z}_p$, where $p$ is a prime, then $E(\Gamma_G)= LE(\Gamma_G)= LE^+(\Gamma_G)$.

\item $E(\Gamma_G) < LE^+(\Gamma_G) < LE(\Gamma_G)$ if $G \cong D_{2m} (m \ne 4)$, $QD_{2^n}$, $M_{2rs} (r \ne 4)$, $Q_{4n} (n \ne 2)$, $U_{6n}$, $SD_{8n}$ and $V_{8n}$ ($n \ne 2$). Also, if   $\frac{G}{Z(G)} \cong D_{2m} (m \geq 3)$ and  $Sz(2)$ then $E(\Gamma_G) < LE^+(\Gamma_G) < LE(\Gamma_G)$. 

		
		
		
		
		
		
		
		
		
		

		

		\item 
		
		
		
		

		
		
		
		
		
		
		
		$\Gamma_G$ is non-hypoenergetic as well as non-hyperenergetic if $G \cong D_{2m}, QD_{2^n}, M_{2rs},  Q_{4n}, U_{6n}, A(n,  \mathcal{V})$ , $A(n, p), SD_{8n}$ and $V_{8n}$. Also, if $\frac{G}{Z(G)} \cong D_{2m},  \mathbb{Z}_p \times \mathbb{Z}_p$ and $Sz(2)$, where $m \geq 3$ and $p$ is a prime, then  	$\Gamma_G$ is non-hypoenergetic as well as non-hyperenergetic.

\item $\Gamma_G$ is L-hyperenergetic but not Q-hyperenergetic if $G \cong D_6, M_6$ and $Sz(2)$. 

		
		
		
\item $\Gamma_G$ is neither L-hyperenergetic nor Q-hyperenergetic if $G \cong D_8, M_{8s}, Q_8, A(n,  \mathcal{V}), A(n, p), V_{16}$. Also, if $\frac{G}{Z(G)} \cong \mathbb{Z}_p \times \mathbb{Z}_p$ then $\Gamma_G$ is neither L-hyperenergetic nor Q-hyperenergetic.

		
		
		
		
		
		

\item $\Gamma_G$ is L-hyperenergetic as well as Q-hyperenergetic if $G \cong D_{2m} (m \ne 3, 4)$, $QD_{2^n}$, $M_{2rs} (2rs \ne 6, 8s)$, $Q_{4n} (n \neq 2)$, $U_{6n}$,  $SD_{8n}$ and  $V_{8n} (n \ne 2)$. Also, if $\frac{G}{Z(G)} \cong Sz(2)$ ($G \ncong 
Sz(2)$) and $D_{2m}$ ($m = 3, 4$ and $|Z(G)|\ne 1$ or $m \geq 5$ and $|Z(G)| \geq 1$) then $\Gamma_G$ is L-hyperenergetic as well as Q-hyperenergetic.

		
		
				
		
		
		
		
		

	\end{enumerate}
\end{theorem}

In the following theorem we get an  example of a graph (non-commuting graph of the symmetric group of degree $4$) disproving Conjecture \ref{con-hyper}.
\begin{theorem}\label{com}
Let the commuting graph of a finite group $G$ be planar. Then
\begin{enumerate}
\item $\Gamma_G$ is Q-integral if 	$G \cong D_8, Q_8, \mathbb{Z}_2 \times D_8, \mathbb{Z}_2 \times Q_8, \mathcal{M}_{16}, \mathbb{Z}_4 \rtimes \mathbb{Z}_4, D_8 \ast \mathbb{Z}_4$ or $SG(16,3)$, otherwise not  Q-integral.
\item $\Gamma_G$ is non-hypoenergetic.
\item $\Gamma_G$ is hyperenergetic if $G \cong S_4$, otherwise non-hyperenergetic.
\item $\Gamma_G$ is Q-hyperenergetic if  $G \cong D_{10}, D_{12}, Q_{12}, A_4, A_5, S_4 \text{ or } SL(2,3)$, otherwise 
 not Q-hyperenergetic
\item $\Gamma_G$ is L-hyperenergetic if \, $G \, \cong \,  D_6, \, D_{10}, \, D_{12}, \, Q_{12},\, A_4, \, A_5, \, S_4, \, SL(2,3) \text{ or } Sz(2)$, otherwise not\\ L-hyperenergetic.
\end{enumerate}
\end{theorem}
\begin{proof}
If the commuting graph of $G$ is planar then, by  \cite[Theorem 2.2]{AFK15}, $G \cong D_6, D_8, D_{10}, D_{12}, Q_8, Q_{12}, \mathbb{Z}_2 \times D_8, \mathbb{Z}_2 \times Q_8, \mathcal{M}_{16}, \mathbb{Z}_4 \rtimes \mathbb{Z}_4, D_8 \ast \mathbb{Z}_4, SG(16,3), A_4, A_5, S_4, SL(2,3), Sz(2)$.

By Theorems \ref{Dihedral1}, \ref{Dihedral2} and \ref{Quarternion}, $\Gamma_G$ is Q-integral if 	$G \cong D_8, Q_8$ and not Q-integral if 	$G \cong D_6, D_{10}, D_{12}, Q_{12}$. If $G \cong D_6$, then from Theorem \ref{D_{2m}}, $\Gamma_G$ is non-hypoenergetic, non-hyperenergetic, not Q-hyperenergetic but is L-hyperenergetic. If $G \cong D_{10}, D_{12}$, then from Theorem \ref{D_{2m}}, $\Gamma_G$ is non-hypoenergetic, non-hyperenergetic but is  Q-hyperenergetic and L-hyperenergetic. If $G \cong D_8$, then from Theorem \ref{D_{2m}}, $\Gamma_G$ is non-hypoenergetic, non-hyperenergetic, not Q-hyperenergetic and not L-hyperenergetic.

 If $G \cong Q_8$, then from Theorem \ref{Q_{4m}}, $\Gamma_G$ is  non-hypoenergetic, non-hyperenergetic, not Q-hyperenergetic and not L-hyperenergetic. If $G \cong Q_{12}$, then from Theorem \ref{Q_{4m}} $\Gamma_G$ is  non-hypoenergetic, non-hyperenergetic but is  Q-hyperenergetic and L-hyperenergetic.

 If $G \cong \mathbb{Z}_2 \times D_8, \mathbb{Z}_2 \times Q_8, \mathcal{M}_{16}, \mathbb{Z}_4 \rtimes \mathbb{Z}_4, D_8 \ast \mathbb{Z}_4, SG(16,3)$, then $\frac{G}{Z(G)} \cong \mathbb{Z}_2 \times \mathbb{Z}_2$. Using Theorems \ref{Z_p*Z_p} and \ref{Zp-energetic}, for $p=2$, we get $\Gamma_G$ is Q-integral but not hypoenergetic, hyperenergetic,  Q-hyperenergetic as well as L-hyperenergetic.

 If $G \cong A_4$ then from Proposition 4.3.13 of \cite{FWNT21}, we have $E(\Gamma_{A_4})=6+2\sqrt{33}$ and $LE(\Gamma_{A_4})=\frac{224}{11}$. Now, $|v(\Gamma_{A_4})|=11$ so $E(K_{|v(\Gamma_{A_4})|})=LE^+(K_{|v(\Gamma_{A_4})|})= LE(K_{|v(\Gamma_{A_4})|})= 20$. Here, $\Gamma_{A_4}=K_{4.2,1.3}$ so using Theorem \ref{R1}(b) we get 
 \[
 \Q-spec(\Gamma_{A_{4}})=\left\lbrace(9)^4, (8)^2, (7)^3, \left(\frac{23+\sqrt{145}}{2}\right)^1, \left(\frac{23-\sqrt{145}}{2}\right)^1\right\rbrace.
 \]
 It follows that $\Gamma_{A_4}$ is not Q-integral. We have $|e(\Gamma_{A_{4}})|=48$ and so $\frac{2|e(\Gamma_{A_{4}})|}{|v(\Gamma_{A_{4}}|}=\frac{96}{11}$. Therefore, $\left|9-\frac{96}{11}\right|=\frac{3}{11}$, $\left|8-\frac{96}{11}\right|=\frac{8}{11}$, $\left|7-\frac{96}{11}\right|=\frac{19}{11}$, $\left|\frac{23+\sqrt{145}}{2}-\frac{96}{11}\right|=\frac{61}{22}+\frac{\sqrt{145}}{2}$ and $\left|\frac{23-\sqrt{145}}{2}-\frac{96}{11}\right|=-\frac{61}{22}+\frac{\sqrt{145}}{2}$. Thus, 
 \[
 LE^+(\Gamma_{A_{4}})=4 \times \frac{3}{11}+2 \times \frac{8}{11}+3 \times \frac{19}{11}+\frac{61}{22}+\frac{\sqrt{145}}{2}-\frac{61}{22}+\frac{\sqrt{145}}{2}=\frac{85}{11}+\sqrt{145} > 20.
 \]
  Hence, $\Gamma_{A_4}$ is  non-hypoenergetic, non-hyperenergetic but is Q-hyperenergetic and L-hyperenergetic.

 If $G \cong A_5$ then from Proposition 4.3.13 of \cite{FWNT21},  we have $E(\Gamma_{A_5})=111.89$ and $LE(\Gamma_{A_5})=\frac{8580}{59}$. Now, $|v(\Gamma_{A_5})|=59$ so $E(K_{|v(\Gamma_{A_5})|})=LE^+(K_{|v(\Gamma_{A_5})|})=LE(K_{|v(\Gamma_{A_5})|})= 116$. Here, $\Gamma_{A_5}=K_{5.3,10.2,6.4}$ so using Theorem \ref{R1}(b) we get 
 \[
 \Q-spec(\Gamma_{A_{5}})=\{(56)^{10}, (57)^{10}, (55)^{27}, (53)^4, (51)^5, (x_1)^1, (x_2)^1, (x_3)^1\},
 \]
  where $x_1 \approx 52.03252, x_2\approx 54.05266$ and  $x_3 \approx 111.91482$ are the roots of the equation $x^3-218x^2+14685x-314760=0$. It follows that $\Gamma_{A_5}$ is not Q-integral. We have $|e(\Gamma_{A_{5}})|=1650$ and so $\frac{2|e(\Gamma_{A_5})|}{|v(\Gamma_{A_5}|}=\frac{3300}{59}$. Therefore, $\left|57-\frac{3300}{59}\right|=\frac{63}{59}$, $\left|56-\frac{3300}{59}\right|=\frac{4}{59}$, $\left|55-\frac{3300}{59}\right|=\frac{55}{59}$, $\left|53-\frac{3300}{59}\right|=\frac{173}{59}$, $\left|51-\frac{3300}{59}\right|=-\frac{291}{59}+\frac{\sqrt{145}}{2}$,
 $\left|x_1-\frac{3300}{59}\right| = -(x_1-\frac{3300}{59})$,
 $\left|x_2-\frac{3300}{59}\right| = -(x_2-\frac{3300}{59})$
 and $\left|x_3-\frac{3300}{59}\right| = x_3-\frac{3300}{59}$.
 Thus, 
 $$
 LE^+(\Gamma_{A_{5}}) = \frac{7602}{59} - x_1 - x_2 + x_3 > 116.
 $$
 Hence, $\Gamma_{A_5}$ is  non-hypoenergetic, non-hyperenergetic but is  Q-hyperenergetic and L-hyperenergetic.

 If $G \cong S_4$ then from Proposition 4.3.13 of \cite{FWNT21}, we have $E(\Gamma_{S_4})=35.866+4\sqrt{5}$ and $LE(\Gamma_{S_4})=\frac{1072}{23}+4\sqrt{13}$. Now, $|v(\Gamma_{S_4})|=23$ so $E(K_{|v(\Gamma_{S_4})|})=LE^+(K_{|v(\Gamma_{S_4})|})=LE(K_{|v(\Gamma_{S_4})|})=44$. Using GAP \cite{GAP2017}, the characteristic polynomial of $Q(\Gamma_{S_4})$ is 
 $$
 Q_{\Gamma_{S_4}}(x) = x(x+20)^4(x+21)^7(x+23)^7(x^2+40x+394)^2
 $$
  and so 
$$
\Q-spec(\Gamma_{S_{4}})=\left\lbrace (0)^1, (-20)^4, (-21)^7, (-23)^7, \left(-20+\sqrt{6}\right)^2, \left(-20-\sqrt{6}\right)^2\right\rbrace.
$$
It follows that $\Gamma_{S_4}$ is not Q-integral. We have $|e(\Gamma_{S_{4}})|=228$ and so $\frac{2|e(\Gamma_{S_4})|}{|v(\Gamma_{S_4}|}=\frac{456}{23}$. Therefore, $\left|0-\frac{456}{23}\right|=\frac{456}{23}$, $\left|20-\frac{456}{23}\right|=\frac{4}{23}$, $\left|21-\frac{456}{23}\right|=\frac{27}{23}$,  $\left|23-\frac{456}{23}\right|=\frac{73}{23}$,  $\left|-20+\sqrt{6}-\frac{456}{23}\right|=\frac{916}{23}-\sqrt{6}$ and $\left|-20-\sqrt{6}-\frac{456}{23}\right|=\frac{916}{23}+\sqrt{6}$. Thus, 
 $$
 LE^+(\Gamma_{S_{4}})=\frac{456}{23} + 4 \times \frac{4}{23}+7 \times \frac{27}{23}+7 \times \frac{73}{23}+2 \times \left(\frac{916}{23}-\sqrt{6}\right)+2 \times \left(\frac{916}{23}+\sqrt{6}\right)=\frac{4836}{23}.
 $$
   Hence, $\Gamma_{S_4}$ is  non-hypoenergetic, hyperenergetic but is Q-hyperenergetic and  L-hyperenergetic.

 If $G \cong SL(2,3)$ then from Proposition 4.3.13 of \cite{FWNT21}, we have $E(\Gamma_{SL(2,3)})=16+8\sqrt{7}$ and $LE(\Gamma_{SL(2,3)}))=\frac{552}{11}$. Now, $|v(\Gamma_{SL(2,3)}))|=22$ so $E(K_{|v(\Gamma_{SL(2,3)}))|})=LE^+(K_{|v(\Gamma_{SL(2,3)})|})= LE(K_{|v(\Gamma_{SL(2,3)})|})= 42$. Here, $\Gamma_{SL(2,3)}=K_{3.2,4.4}$ so using Theorem \ref{R1}(b) we get
 $$
 \Q-spec(\Gamma_{SL(2,3)})=\left\lbrace(20)^3, (18)^{14}, (14)^3, \left(\frac{54+\sqrt{420}}{2}\right)^1, \left(\frac{54-\sqrt{420}}{2}\right)^1\right\rbrace.
 $$
It follows that $\Gamma_{SL(2, 3)}$ is not Q-integral.  We have $|e(\Gamma_{SL(2,3)})|=204$ and so $\frac{2|e(\Gamma_{SL(2,3)})|}{|v(\Gamma_{SL(2,3)}|}=\frac{204}{11}$. Therefore, $\left|20-\frac{204}{11}\right|=\frac{16}{11}$, $\left|18-\frac{204}{11}\right|=\frac{6}{11}$, $\left|14-\frac{204}{11}\right|=\frac{50}{11}$, $\left|\frac{54+\sqrt{420}}{2}-\frac{204}{11}\right|=\frac{93}{11}+\frac{\sqrt{420}}{2}$ and $\left|\frac{54-\sqrt{420}}{2}-\frac{204}{11}\right|=-\frac{93}{22}+\frac{\sqrt{420}}{2}$. Thus, 
  $$
  LE^+(\Gamma_{SL(2,3)})=3 \times \frac{16}{11}+14 \times \frac{6}{11}+3 \times \frac{50}{11}+\frac{93}{11}+\frac{\sqrt{420}}{2}-\frac{93}{22}+\frac{\sqrt{420}}{2}=\frac{282}{11}+\sqrt{420}.
  $$
   Hence, $\Gamma_{SL(2,3)})$ is  non-hypoenergetic, non-hyperenergetic but is Q-hyperenergetic and L-hyperenergetic.

 If $G \cong Sz(2)$ then, by Theorem \ref{Suzuki}, we have $\Gamma_G$ is not Q-integral. Also,  Theorem \ref{SSzz} gives that $\Gamma_G$ is  non-hypoenergetic, non-hyperenergetic, not Q-hyperenergetic but is L-hyperenergetic. 
\end{proof}

\begin{theorem}
Let $G$ be a finite group and the commuting graph of $G$ is toroidal. Then
\begin{enumerate}
\item  $\Gamma_G$ is Q-integral if $G \cong D_{14}$ or $A_4 \times \mathbb{Z}_2$, otherwise not Q-integral.  	
\item  $\Gamma_G$ is non-hypoenergetic, non-hyperenergetic but is Q-hyperenergetic and L-hyperenergetic.
\end{enumerate}
\end{theorem}
\begin{proof}
If commuting graph of $G$ is toroidal then, by \cite[Theorem 3.3]{DN-17} $G \cong D_{14}, D_{16}, Q_{16}, QD_{16}, \mathbb{Z}_7 \rtimes \mathbb{Z}_3, D_6 \times \mathbb{Z}_3, A_4 \times \mathbb{Z}_2$. 
By Theorems \ref{Dihedral1}, \ref{Quasidihedral} and \ref{Quarternion}, $\Gamma_G$ is Q-integral if 	$G \cong D_{14}$ and not Q-integral if 	$G \cong D_{16}, Q_{16}$ or $QD_{16}$. If $G \cong D_{14}, D_{16}$, then from Theorem \ref{D_{2m}}, $\Gamma_G$ is  non-hypoenergetic, non-hyperenergetic but is Q-hyperenergetic and L-hyperenergetic.
 If $G \cong Q_{16}$, then from Theorem \ref{Q_{4m}}, $\Gamma_G$ is  non-hypoenergetic, non-hyperenergetic but is Q-hyperenergetic and L-hyperenergetic.
 If $G \cong QD_{16}$, then from Theorem \ref{QDn}, $\Gamma_G$ is  non-hypoenergetic, non-hyperenergetic but is Q-hyperenergetic and L-hyperenergetic. 

 If $G \cong \mathbb{Z}_7 \rtimes \mathbb{Z}_3$ then, by  Theorem \ref{pq1}, we have 
\[
		\Q-spec(\Gamma_G)= \left\{(14)^{5},(18)^{7},(16)^{6}, \left(22 + 2\sqrt{37}\right)^1, \left(22 - 2\sqrt{37}\right)^1 \right\}.
\] 
Thus $\Gamma_G$ is not Q-integral. By Theorem \ref{pq}, we also have $E(\Gamma_G)=12+4\sqrt{30}$  and $LE(\Gamma_G)=\frac{308}{5}$ and from Theorem \ref{pq1}, we have $LE^+(\Gamma_G)=\frac{292}{20}+4\sqrt{37}$. Now, $|v(\Gamma_G)|=20$ so $E(K_{|v(\Gamma_G)|})=LE^+(K_{|v(\Gamma_G)|})=LE(K_{|v(\Gamma_G)|})=38$. Hence, $\Gamma_G$ is  non-hypoenergetic, non-hyperenergetic but is Q-hyperenergetic and L-hyperenergetic.


  If $G \cong D_6 \times \mathbb{Z}_3$, then from Proposition 4.3.14 of \cite{FWNT21}, we have $E(\Gamma_G)=6+6\sqrt{7}$ and $LE(\Gamma_G)=\frac{594}{15}$. Now, $|v(\Gamma_G)|=15$ so $E(K_{|v(\Gamma_G)|})=LE^+(K_{|v(\Gamma_G)|})= LE(K_{|v(\Gamma_{G})|})= 28$. Here, $\Gamma_{G}=K_{3.3,1.6}$ so using Theorem \ref{R1}(b) we get
 $$
  \Q-spec(\Gamma_{G})=\left\lbrace(12)^6, (9)^7, \left(\frac{27+\sqrt{297}}{2}\right)^1, \left(\frac{27-\sqrt{297}}{2}\right)^1\right\rbrace.
  $$
 It follows that $\Gamma_G$ is not Q-integral.   We have $|e(\Gamma_{G})|=81$ and so $\frac{2|e(\Gamma_{G})|}{|v(\Gamma_{G}|}=\frac{162}{15}$. Therefore, $\left|12-\frac{162}{15}\right|=\frac{18}{15}$, $\left|9-\frac{162}{15}\right|=\frac{27}{15}$, $\left|\frac{27+\sqrt{297}}{2}-\frac{162}{15}\right|=\frac{81}{30}+\frac{\sqrt{297}}{2}$ and $\left|\frac{27-\sqrt{145}}{2}-\frac{162}{15}\right|=-\frac{81}{30}+\frac{\sqrt{297}}{2}$. Thus, 
 $$
 LE^+(\Gamma_{G})=6 \times \frac{18}{15}+7 \times \frac{189}{11}+\frac{81}{30}+\frac{\sqrt{297}}{2}-\frac{81}{30}+\frac{\sqrt{297}}{2}=\frac{99}{5}+3\sqrt{33}.
 $$
  Hence, $\Gamma_G$ is  non-hypoenergetic, non-hyperenergetic but is Q-hyperenergetic and L-hyperenergetic.

 If $G \cong A_4 \times \mathbb{Z}_2$, then from Proposition 4.3.14 of \cite{FWNT21}, we have $E(\Gamma_G)=12+4\sqrt{33}$ and $LE(\Gamma_G)=\frac{544}{11}$. Now, $|v(\Gamma_G)|=22$ so $E(K_{|v(\Gamma_G)|})=LE^+(K_{|v(\Gamma_G)|})=LE(K_{|v(\Gamma_{G})|})= 42$. Here, $\Gamma_{G}=K_{4.4,1.6}$ so using Theorem \ref{R1}(b) we get 
 $$
 \Q-spec(\Gamma_{G})=\{(18)^{12}, (16)^5, (14)^3, (36)^1, (10)^1\}.
 $$
Clearly, $\Gamma_G$ is  Q-integral.  We have $|e(\Gamma_{G})|=192$ and so $\frac{2|e(\Gamma_{G})|}{|v(\Gamma_{G}|}=\frac{192}{11}$. Therefore, $\left|18-\frac{192}{11}\right|=\frac{6}{11}$, $\left|16-\frac{192}{11}\right|=\frac{16}{11}$, $\left|14-\frac{192}{11}\right|=\frac{38}{11}$, $\left|36-\frac{192}{11}\right|=\frac{204}{11}$ and $\left|10-\frac{192}{11}\right|=\frac{82}{11}$. Thus, 
$$
LE^+(\Gamma_{G})=12 \times \frac{6}{11}+5 \times \frac{16}{11}+ 3 \times \frac{38}{11}+\frac{204}{11}+\frac{82}{11}=\frac{552}{11}.
$$
 Hence, $\Gamma_G$ is  non-hypoenergetic, non-hyperenergetic but is Q-hyperenergetic and L-hyperenergetic.
\end{proof}

If non-commuting graph of $G$ is planar then, by \cite[Proposition 2.3]{AAM06} $G \cong D_6, D_8, Q_8$. Therefore, we have the following theorem.
\begin{theorem}
Let $G$ be a finite group whose non-commuting graph is planar. Then 
\begin{enumerate}
\item $\Gamma_G$ is non-hypoenergetic, non-hyperenergetic and not Q-hyperenergetic.
\item $\Gamma_G$ is not L-hyperenergetic but Q-integral if $G \ncong D_{6}$.
\end{enumerate}
\end{theorem}

\section*{Acknowledgements}
The authors would like to thank Mr. Nabin  Pokhrel and Mr. Kallol Ray for their help in drawing Figures and computing perfect squares in Table 1.  
Ms. M. Sharma expresses gratitude to DST  for the INSPIRE fellowship.

\section*{Funding} No funding was received by the authors.
\section*{Data Availability} No data was used in the preparation of this manuscript.
\section*{Conflict of interest} The authors declare that they have no conflict of interest.

\end{document}